\mag=\magstephalf
%\mag=\magstep1
\pageno=1
\input amstex
\documentstyle{amsppt}
\interlinepenalty=1000
\NoRunningHeads
\pagewidth{17 truecm}
\pageheight{24 truecm}
\vcorrection{-0.5cm}
\hcorrection{-0.7cm}
\nologo
\NoBlackBoxes
\font\twobf=cmbx12
\define \ii{\roman i}
\define \ee{\roman e}
\define \dd{\roman d}
\define \tr{\roman {tr}}
\define \bx{\overline{x}}
\define \bby{\overline{y}}
\define \CC{{\Bbb C}}
\define \CCtwo{{\CC^2\setminus \{0\}}}
\define \EE{{\roman{E} }}
\define \RR{{\Bbb R}}
\define \PP{{\Bbb P}}
\define \ZZ{{\Bbb Z}}
\define \NN{{\Bbb N}}
\define \QQ{{\Bbb Q}}
\define \CCinf{{\Bbb C\cup\{\infty\}}}
\define \PSL{\roman{PSL}_2}
\define \SL{\roman{SL}_2}

\define \SD{\roman{SD}}
\define \elas{\roman{elas}}

\define \cls{{\roman{cls}}}
\define \KdV{{\roman{KdV}}}
\define \hyp{{\roman{hyp}}}

\define \ev{\roman{ev}}
\define \finite{{\roman{finite}}}
\define \Pic{ {\roman{Pic}} }
\define \Cinf{{\Cal C^{\infty}}}
\define \SCinf{{\underline{\Cal C^{\infty}}}}
\define \SF{{\underline{\Cal F }}}
\define \Diff{ {\roman{Diff}^+} }
\define \Isom{ {\roman{Isom}} }
\define \res{ {\roman{res}} }
\define \DDs{{\frak D^s}}
\define \EEs{{\frak E^s}}

\define \DDf{{\frak D^f}}
\define \DDc{{\frak D^c}}
\define \EEf{{\frak E^f}}
\define \VVf{{\frak V^f}}
\define \AAc{{\frak A^c}}
\define \WWc{{\frak W^c}}
\define \WWf{{\frak W^f}}
\define \WWt{{\frak W^t}}

\define \DDt{{\frak D^t}}
\define \EEt{{\frak E^t}}
\define \EEc{{\frak E^c}}

\define \AAf{{\frak A^f}}
\define \aac{{\Cal A^c}}
\define \aaf{{\Cal A^f}}
\define \aat{{\Cal A^t}}
\define \val{{\roman{val}}}
\define \CMeP{{\Cal M_\elas^{\Bbb P}}}
\define \CMPP{{\Cal M^{\Bbb P}}}
\define \CCMeP{{\Cal C\Cal M_\elas^{\Bbb P}}}
\define \DCMeP{{{\Cal D\Cal M_\elas^{\Bbb P}}}}
\define \CMeCtwo{{ \Cal M_\elas^{\CC^2\setminus \{0\}}}}

\define \CMeC{{\Cal M_\elas^{\Bbb C}}}

\define \CMeCpi{{\Cal M_\elas^{\Bbb C, 2\pi}}}

\define \BMeP{{\Bbb M_\elas^{\Bbb P}}}
\define \BMP{{\Bbb M^{\Bbb P}}}
\define \fMeP{{\frak M_\elas^{\Bbb P}}}
\define \BMeCtwo{{\Bbb M_\elas^{\CC^2\setminus \{0\}}}}
\define \BMCtwo{{\Bbb M^{\CC^2\setminus \{0\}}}}
\define \BMeC{{\Bbb M_\elas^{\Bbb C}}}
\define \BMeCpi{\Bbb M_\elas^{\Bbb C, 2\pi}}
\define \BMeinf{\Bbb M_\elas^{\infty}}
\define \BMeinfc{\Bbb M_\elas^{\infty, {\roman cvtr}}}

\define \CMKdV{{\Cal M_\KdV}}
\define \BMKdV{{\Bbb M_\KdV}}
\define \fMKdV{{\frak M_\KdV}}
\define \fMhyp{{\frak M_\hyp}}
\define \Top{{\bold{Top}}}
\define \DGeom{{\bold{DGeom}}}
\define \Spect{ {\roman{Spect}} }
\define \HH{{\roman{H}}}
\define \Hc{{\roman{H}_{\roman{c}}}}

\define \simKdVH{{\underset{\roman{KdVHf}}\to{\sim}}}
\define \Proof{ {\it Proof} }
\define \GCD{\roman{GCD}}
\define \gr{\roman{gr}}
\define \ft{{\frak t}}
\define \pt{{\roman{pt}}}
\define\hw{\widehat{\omega}}

%\define \tvskip{\vskip 0.5 cm}
\define \tvskip{\vskip 1.0 cm}
\define\ce#1{\lceil#1\rceil}
\define\dg#1{(d^{\circ}\geq#1)}
\define\Dg#1#2{(d^{\circ}(#1)\geq#2)}
\define\dint{\dsize\int}

\define\s#1{\sigma_{#1}}
\define\tp#1{\negthinspace\left.\ ^t#1\right.}
\define\mrm#1{\text{\rm#1}}
\define\lr#1{^{\sssize\left(#1\right)}}
\redefine\ord{\text{ord}}

\redefine\qed{\hbox{\vrule height6pt width3pt depth0pt}}

%DRAFT: 98/8/21 version
%\baselineskip  0.6 cm

{\centerline{\bf{ On the Moduli of a Quantized Elastica in $\PP$
and KdV Flows:}}}
{\centerline{\bf{Study of Hyperelliptic Curves as an
Extension of Euler's Perspective of Elastica I}}}

\author
Shigeki MATSUTANI${}^*$ and
 Yoshihiro \^ONISHI${}^\dagger$
\endauthor
\affil
${}^*$8-21-1 Higashi-Linkan Sagamihara, 228-0811 Japan \\
${}^\dagger$Faculty of Humanities and Social Sciences,
 Iwate University,\\
  Ueda, Morioka, Iwate, 020-8550 Japan
\endaffil \endtopmatter

\centerline{\twobf Abstract }\tvskip

Quantization needs evaluation of all of states of a quantized object
rather than its stationary states with respect to its energy.
In this paper, we have investigated  moduli $\CMeP$ of a quantized
elastica, a quantized loop  with an energy functional associated
with the Schwarz derivative, on a Riemann sphere $\PP$.
Then it is proved that its moduli space
is decomposed to a set of equivalent
classes determined by flows obeying the Korteweg-de Vries (KdV)
 hierarchy
which conserve the energy. Since the  flow obeying the KdV hierarchy
has a natural topology, it induces topology in the moduli space
$\CMeP$.  Using the topology,  $\CMeP$ is classified.

Studies on a loop space in  the category of topological spaces
$\Top$
are well-established and its cohomological properties are
well-known. As the moduli space of
 a quantized elastica can be regarded
as a loop space in the category of differential geometry $\DGeom$,
we also proved an existence of a functor between a triangle
category related to a loop space in {\bf  Top}  and  that in
$\DGeom$ using the induced topology.

As Euler investigated the elliptic integrals and its moduli
by observing a shape of classical elastica on $\CC$,
this paper devotes relations between hyperelliptic
curves and a quantized elastica on $\PP$ as an extension of
Euler's perspective of elastica.

\document
\tvskip
%\baselineskip  0.6 cm
%\newpage
\centerline{\twobf \S 1. Introduction }\tvskip

History of investigations of elastica was opened by James Bernoulli in
1691 according to Truesdell's inquiry [T1, 2, L].
He named a shape of a thin non-stretching elastic rod elastica
and proposed the elastica problem:
what shape does elastica take for a given boundary
condition? It should be, further, noted that he also proposed
the lemniscate problem and discovered an elliptic integral corresponding
to the lemniscate function by investigation of elastica.
He considered a smooth curve with the arc-length in a plane $\CC$,
$$
  \tilde \gamma: [0,l] \hookrightarrow \CC,
   \quad ( s \mapsto \tilde \gamma(s) ).
$$
Following his studies, his nephew Daniel Bernoulli discovered that
the elastica obeys the minimal principle
that shape of the elastica is realized as a stationary point of an
energy functional, which is called Euler-Bernoulli functional nowadays,
$$   E[\tilde\gamma] = \int k^2 d s,
$$
where $k$ is the curvature of the curve $\tilde\gamma$ in $\CC$,
$k = -\sqrt{-1} \partial_s^2 \tilde \gamma/ \partial_s \tilde \gamma$,
$\partial_s := d / d s$, and $s$ is the arc-length of the curve using
the induced metric in $\CC$. (It should be noted that this functional
{\it differs} from that of a {\lq\lq}string" in the literature of the
string theory in the elementary particle physics:
although an elastica is a model of a string of the chord {\it e.g},
the guitar, {\lq\lq}string"  in the string theory can not be realized in
the classical mechanical regime.)

Since the curvature $k$ is expressed as
$k = \partial_s \phi$ where $\phi$
is the tangential angle, and the energy is given by
$E = \int | \partial_s \phi |^2 d s$,
the elastica problem could be interpreted as the oldest problem of
a harmonic map into a target space U(1); if we write
 $\partial_s \phi ds = g^{-1} dg$, for U(1) valued function $g$
over $\CC$,
 then the Hodge-star dual $*g^{-1} d g = \partial_s
\phi$ and $E = \int < g^{-1} d g \wedge *g^{-1} d g>$.

The elastica problem is to investigate moduli
$\tilde{\Cal M}_{\elas,\cls}$,
$$   \tilde \Cal M_{\elas,\cls} := \{
   \tilde \gamma: [0,1] \hookrightarrow \CC\ |\
   \delta E[\tilde\gamma]/\delta \tilde\gamma =0 \}/ \sim .
$$
Here \lq\lq$\sim$" means modulo Euclidean
move in $\Bbb  C$ and dilatation.
We sometimes call this space {\it moduli space of the classical elasticas}.
The classification of this moduli space $\tilde \Cal M_{\elas,\cls}$
was essentially done by Euler in 1744 by means of
numerical computations [E]. The moduli space
$\tilde \Cal M_{\elas,\cls}$
is classified by the moduli of the elliptic curves [T1, 2, L, We].
It is noted that before Euler refereed
to Fagnano's paper on his discovery of an algebraic properties
of the lemniscate function (an elliptic function of a special modulus)
at December 31 1751, the elliptic integrals for more general modulus
was investigated in the study of this classical harmonic map
problem. (It is known that Jacobi recognized that the day is
the birthday of the elliptic function. Thus we think that elastica
is a kind of the movements of  the fetus of algebraic
curves.) We also emphasize that from the beginning,
the harmonic map problem
(classical field theory in
physics) is closely related to algebraic varieties.
Recently Mumford investigated this elastica problem from
a viewpoint of
applied mathematics and gave simple and deep expressions of
the shape of elastica, which show the depth, importance
and beauty of this problem [Mum3].

Especially for a closed elastica, Euler showed that
its moduli space,
$$      \Cal M_{\elas,\cls} := \{
   \gamma: S^1 \hookrightarrow \CC|
   \delta E[\gamma]/\delta \gamma =0 \}/\sim ,
$$
consists of two disjoint points: the corresponding moduli
 $\tau$ of the elliptic curves consist of two points
$\tau =0$ and $\tau= 0.70946\cdots$ [E, T1, 2, L].

Recently a loop space is one of the most concerned objects in mathematics
and there have been so many efforts to investigate it
[Br, G, LP, Se, SW and reference therein].
Further it is well-known that soliton equations are
closely related to the loop spaces, loop groups and loop algebras [G, SW].
However these studies are sometimes
too abstract to be related to physical problems, except  problems in the
elementary particle physics; for example, the embedded space is often
a group manifold, {\it e.g.}, U($N$).  Further the energy function is
paid little attention in these studies.

On the other hand, our concerned object is a non-stretching
 elastica, which
 is related to a large polymer, such as the deoxyribonucleic acid
(DNA) as a physical model [Mat1, 2, 4, KV].  Elastica has an
energy functional as we described above.
Thus our problem, basically, differs from the arguments in an
ordinary loop space in [G, Se, SW] except [Br, LP]
though it is closely related to them.

One of these authors (S.M.) considered the quantization of a closed
elastica (precisely speaking, statistical mechanics of elasticas)
[Mat2]. He defined the moduli space of the closed quantized elastica,
which is an isometric immersion of $S^1$ into $\CC$ module the
Euclidean motion and dilatation,
$$
   \Cal M_{\elas}^{\CC}  := \{  \gamma : S^1
           \hookrightarrow \CC| \text{ isometric immersion } \}/\sim   .
$$
He investigated the partition function from a physical point of view,
which has not been mathematically justified:
$$   Z : \Cal M_{\elas}^{\CC}  \times \Bbb R_{>0}
                    \longrightarrow \Bbb R,
$$ with
$$   Z[\beta] = \int_{\Cal M_{\elas}^{\CC} } D \gamma
   \exp( -\beta E[\gamma]),
$$
where $\beta \in \Bbb R_{>0}:=\{x \in \Bbb R\ | \ x>0 \}$ and
$D \gamma$ is the Feynman measure.
On the quantization of an elastica, we need more information
of the moduli space of curves besides those around its
stationary points.
To evaluate this map $Z$, he classified the moduli space
of a quantized
closed elastica $\CMeC$ and attempted to redefine
the Feynman measure by replacing it with the series of Riemann
integral over $\CMeC$.
His quantization is somewhat novel for an elastica.
He physically proved that the moduli space of the quantized
elastica is given as a subspace of the moduli space
of the modified Korteweg-de Vries (MKdV) equation [Mat2].

Here we should emphasize that it is very surprising that a physical
system is completely described by a soliton equation as mentioned
in [Mat2].  Even in  physical phenomena which are known as systems
represented by soliton equations, like shallow waves, plasma waves,
charge density waves and so on, the higher soliton solutions are,
in general,  out of their approximation regions; of course one or
two soliton solutions do represent these phenomena well.
On the other hand, in the quantized elastica problem,
its functional space is completely expressed by the MKdV
hierarchy, even though problems in polymer physics are, in general,
too complex to be solved exactly [KV].

In this paper, we will rewrite the physical theorem in [Mat2] from
a mathematical point of view and extend it.
Pedit gave a lecture
on a loop space over a Riemann sphere $\PP$
at  Tokyo Metropolitan University in 1998 \cite{Ped}.
There he showed that the loop space is related to the
 Korteweg-de Vries (KdV) flow by considering
a loop in $\CCtwo$.
As his treatment is given in the framework of
pure mathematics,
we will follow the expressions of Pedit and
deal with the KdV flow instead of the MKdV flow here.
Due to the Miura map (the Ricatti type
differential equation), the MKdV flow and the KdV flow can be regarded
as different aspects of the same object; this choice is not significant.
Mathematical investigations on the KdV flow
leads us to our main results, Theorems 3-4,  4-2 and 7-4.

As we will show later, our investigation of  a quantized elastica leads
us to study  the hyperelliptic curves and their moduli space as Euler
encountered the elliptic integrals and studied of the moduli of
the elliptic functions by observing a shape of classical elastica on $\CC$.
One of our purposes of this study is to know  the hyperelliptic functions
and its moduli by investigating a quantized
elastica in $\PP$ as an extension of Euler's perspective of elastica.
After we submitted the first version of this paper,
these works progressed [Mat7-10]. Hence in this revised version,
we also rewrite the related parts.

Contents of this paper is as follows.

\S 2 shows an expression of a real curve immersed in a Riemann
sphere $\PP$ according to the lecture of Pedit [Ped]. Using his
expressions, we define the moduli space
of a real smooth curve immersed in
$\PP$ and an energy functional of the curve whose integrand is the
Schwarz derivative along the curve. When we regard $\PP$ as complex
plane with the infinity point, $\CCinf$,  the energy functional is
identified with the Euler-Bernoulli  energy functional around the origin
$\{0\}$ of $\CCinf$ and
the curve with the energy  is reduced to a quantized elastica which
was studied by one of these authors [Mat2].
Thus we  continue to
refer such a curve in $\PP$ \lq\lq
quantized elastica in $\PP$".
In order to consider a quantum effect, we should get knowledge of a
set of curves with different energies instead of investigation of only a
stationary point of the energy functional even though we are dealing
with a single elastica.
Thus we will call, in this paper, the moduli
$\BMeP$ defined in
Definition 2-10 and 2-12,
{\lq\lq}moduli of a quantized elastica" rather than moduli of loops.
In \S 2, we will give an equivalence between a loop space over
$\CCtwo$ and $\PP$ in a certain sense. Further following
MacLaughlin and Beylinski [McLau, Br], we will introduce a
natural topology of the loop space which is induced from
the topology of the base space.

In \S 3, we introduce infinite dimensional parameters $t=(t_1,t_2,t_3,
\cdots)$ which deform a given curve and  define a flow obeying the KdV
hierarchy  along $t$, which is called KdVH flow.
First we give our first main
Theorem 3-4 in this paper. Since the energy functional of a curve turns
out to be the first integral with respect to the parameter $t$,
we  prove that using the KdVH flow we can classify
the moduli $\BMeP$ of a quantized elastica  in $\PP$.
In other words, the moduli space $\BMeP$ is decomposed to a set of the
equivalent classes with respect to the KdVH flow.
As \S 3 gives the differential geometrical and dynamical properties
of the quantized elastica, we will attempt to express the theorem
in terms of the words of the differential geometers.
Remark 3-10 is a key of the study in \S 3.

Primary considerations
 leads the fact that the moduli space
of a quantized elastica $\BMeP$
is a subspace of the moduli space of the KdVH flow $\BMKdV$
as shown in Proposition 4-29.
The system of the KdV hierarchy has
a natural topology, which essentially determines the algebraic
properties of the KdV hierarchy [D, S, SN, SS, SW].
Using results on these studies of the KdV hierarchy,
we give finer classification of the moduli space  $\BMeP$ in
Theorem 4-2 and Proposition 4-33, which is our second main theorem.
There a dense subspace in $\BMeP$ is decomposed by a subspace
characterized by a natural number.
As we defined below Lemma 4-1, we encounter a
finite type of the KdVH flow, which corresponds to the finite type
solutions of the KdV equation and are related to a hyperelliptic
curve.  The natural number is related to genus of the hyperelliptic curve.

In order that we mention our second main statements, Theorem 4-2 and
Proposition 4-33,
\S 4 reviews the algebro-geometrical properties of the KdV hierarchy
based upon the so-called Sato-Mulase theory [Mul, SS, SN].
As the completion of set of finite type solutions is equal to $\CMKdV$,
we concentrate our attention on the finite solution of the KdV flow
and consider $\CMKdV$ algebro-geometrically.
As Sato-Mulase theory is of the algebraic analysis and
is based upon the formal power series ring,
we replace the base ring of smooth functions by
the formal power series.
There we find that a commutative differential ring is
connected with geometry of a commutative ring, {\it i.e.},
a hyperelliptic curve.
Using the inclusion $\CMeP\subset\CMKdV$, we will
introduce the  relative topology in $\CMeP$ induced from the
topology of $\CMKdV$.

In \S 5, we will show another algorithm of explicit computation of
solutions of the KdV flow.
There we will reconsider the KdV equation in the framework of
inverse scattering method and comment the meanings of Theorem 4-2
again at Proposition 5-20. In other words, we will rewrite our second
result more analytically.
So readers can skip this section except Example 5-21.
There we will also review Krichever's construction of
algebro-geometrical solutions [Kr, BBEIM]
and Baker's original method given about one hundred years ago
 [Ba2, Ma7].
 Using it we showed that there is an
injection from the moduli space $\frak M_\hyp$ of hyperelliptic
curves to the moduli space
$\frak M_\KdV$ of the KdV equation up to an
ambiguity; this correspondence enables us to  determine function
forms of hyperelliptic $\wp$ functions  as solutions of the
 KdV equations
for any algebraically given
hyperelliptic curves including  degenerate curves.

\S 6 is digression and we will review a result of a loop space over
$S^2$ in the category of topological space $\Top$, whose morphism is
a continuous map, following the arguments in the textbook of
Bott and Tu [BT]. Studies
on a loop space in $\Top$ are well-established and its cohomological
properties are well-known. On the other hand, the moduli space of a
quantized elastica in $\PP$ can be regarded as a loop space in the
category of the differential geometry $\DGeom$. Thus by loosening the
properties in $\DGeom$ and regarding them
as those in $\Top$,  it is expected that the
moduli of a quantized elastica $\CMeP$ in $\PP$ are
topologically related to those of a loop space in $\Top$.
Thus in \S 6, we will
review  a loop space in $\Top$ and show its cohomological properties.

In \S 7, we will mention the topological properties of the moduli of a
quantized elastica $\CMeP$ and give our third main theorem.
 As loop spaces in both $\Top$ and $\DGeom$ are not finite dimensional
spaces when we regard them as  manifolds in an appropriate sense,
 it is not known that de Rham's
theorem can be applicable to them.
However it is expected that cohomological
sequences in both categories should correspond to each other. In other
words, it is important to argue existence of functor between triangle
categories related to them, {\it i.e.}, quasi-isomorphism.
Precisely speaking, though the closed condition and the reality
condition in the moduli $\CMeP$
make its topological properties difficult to treat,  we
will tune the low dimensional parts of chain complex of $\CMeP$ and
consider a complex of a quotient spaces
$\CCMeP$.  Then we will show
existence of a functor between the triangle categories in
loop spaces in both $\Top$ and
$\DGeom$ as our third main theorem at Theorem 7-4.
The existence of the functor means the our theory in $\DGeom$
is justified in topological investigation.
We believe that this result is  meaningful to  the investigations
of the loop space.

\S 8 gives the remarks and comments upon our results.
First we will comment upon sequences of homotopy of loop spaces in both
$\Top$ and $\DGeom$. Next, we will give a possibility of computations
of the partition function of a quantized elastica in $\CC$. Even in the
quantized system, we will show that the orbit space is meaningful,
whereas it is well-known that in noncommutative space, concepts of orbit
and geometry are sometimes nonsense [C]. So we will comment upon the
fact. Further we will remark the relations between our system and
Painlev\`e equation of the first kind [Mat2, In], and between our system
and conformal field theory. Finally we will comment upon our results
from the a point of view recent progress of Dirac operator related to
immersion object based upon [Mat1-6]. We will also mention possibility
of higher dimensional case of our consideration there.

\tvskip
\centerline{\twobf Acknowledgment}
\tvskip
One of us (S.M.) would like to thank
Prof.~F.~Pedit and Prof.~K.~Tamano for
critical discussions and drawing his attention to this problem. It is
acknowledged that Prof.~K.~Tamano have taught him algebraic topology
and differential geometry based upon [BT] and
[Br] for over this decade and critically
read this manuscript. He also thanks Prof.~S.~Saito, Prof.T.~Tokihiro,
W.~Kawase and H.~Mitsuhashi for helpful discussions and comments in
early stage of this study.
Prof.~K.~Sogo privately suggested him that soliton equations should be
expressed in a projective space before starting this study
and thus this study is one of answers to his suggestions.

He thanks to Prof.~A.~Koholodnko for telling him the reference [Br]
and so many encouragements and discussions by using e-mails and
Prof.~B.~L.~Konopelchenko for kind letters to encourage his works.

He is also grateful to Prof.~Y.~Ohnita, Prof.~M.~Guest,
Dr.~R.~Aiyama and Prof.~K.~Akutagawa
for inviting him their seminars and for critical
discussions and especially Prof.~M.~Guest for sending him the reference
[Se].
Further we thank to Prof. J. McKay for his interest on this article;
his kind comment encouraged us revising the manuscript.
Finally we would like to express our
sincere thanks to the referee for appropriate suggestions, which
improved this article.

%\newpage
\tvskip
\centerline{\twobf Notations}
\tvskip

$\RR$ and $\CC$ are real and complex number fields respectively.
$\RR_{\ge0}$ is the set of the non-negative real numbers.
$\ZZ$ are the set of integers and $\NN$ is the set of natural numbers
$1, 2, 3, \cdots$.
$\ZZ_{\ge0}$ is the set of the non-negative integers.
$\Cinf(A,B)$ means the set of $B$-valued smooth functions
over $A$.
$R[x_1, \cdots, x_n]$
is the set of polynomial of $x_1, \cdots, x_n$ with $R$ valued coefficients
and
$R[[x_1, \cdots, x_n]]$ is the set of formal power series of
$x_1, \cdots, x_n$ with $R$ valued coefficients.
Others important quantities are listed as follows.

$$
\vbox{
\halign{\hfil#\hfil & \quad#\hfil &\quad #\hfil \cr
$\BMP$     :& {Moduli of Loops in $\PP$}&{Defintion 2-4}\cr
$\BMCtwo$ :&{Moduli of Loops in $\CCtwo$,
             $\varpi:\BMCtwo\to\BMP$}&{Defintion 2-4, Remark 2-5}\cr
$\{\gamma,s\}_{\SD}$ :&{Schwarz derivative }&{Defintion 2-6}\cr
$\BMeP$ :&{Moduli of quantized elastica in $\PP$}
               &{Defintion 2-10}\cr
$\BMeC$ :&{Moduli of quantized elastica in $\CC$}
               &{Defintion 2-10}\cr
$\BMeCtwo$ :&{$\varpi:\BMeCtwo\to\BMeP$}&
              {Defintion 2-10, Remark 2-11}\cr
$\CMeP$ :&$\pi_{\elas}^{\PP} :\BMeP \to {\CMeP}$
               &{Defintion 2-12}\cr
$\CMeC$:&$\pi_{\elas}^{\CC} :\BMeP \to {\CMeC}$
               &{Defintion 2-12}\cr
$\CMeCtwo$ :&$\pi_{\elas}^{\CCtwo} :\BMeP \to {\CMeCtwo}$&
                {Defintion 2-12}\cr
$\Cal E[\gamma]$:&{energy of elastica in $\PP$}&
{Defintion 2-18}\cr
$\DDs$:&{Differential ring over $\Cinf(S^1,\CC)$}&
                {Defintion 3-1}\cr
$\EEs$ :& {Micro differential ring to $\DDs$}&
                {Defintion 3-1}\cr
${\Cal V}^\infty$ :& $S^1 \times (\prod_{n=1}^\infty \Bbb R)$&
                {Defintion 3-2}\cr
$\overline\phi_{\partial_s u,t},\quad
 \overline\varphi_{\partial_s u,t}$
                  :& the KdVH flow &
                {Defintion 3-2, Proposition 3-11}\cr
$\Omega, \quad \underline\Omega, $
                  :&Recursion differenital operator &
                {Defintion 3-2, Lemma 3-6}\cr
$\simKdVH $
               :& Equivlent relation related to the KdVH flow&
                {Defintion 3-2}\cr
$\phi_{A,t}$
                  :& a flow for $A$&
                {Defintion 3-7}\cr
$\fMeP$ : & $\pi_\elas : \BMeP\to \fMeP:=\BMeP/\simKdVH$
                   & Definition 3-18\cr
}}
$$
$$
\vbox{
\halign{\hfil#\hfil & \quad#\hfil &\quad #\hfil \cr
$\BMeP_{, \finite}, \quad \BMeP_g$
                :& finite type of the KdV flow and finite $g$-type flow &
                 Theorem 4-2 \cr
$\DDf$:&{Differential ring over $\CC[[t_1]]$}&
                {Definition 4-4}\cr
$\EEf$ :& {Micro differential ring to $\DDf$}&
                {Definition 4-4}\cr
$\EEf$ :& {Micro differential ring with coefficient $\CC$}&
                {Definition 4-4}\cr
$\WWf,\quad \WWc$ :& {Subsets of $\EEf$ and $\EEc$}&
                {Definition 4-4}\cr
$\Cal L$ :& {Subset of $\EEf$}&
                {Lemma 4-8}\cr
$\aaf, \quad \aac$ :& {commutative subrings of $\EEf$ and $\EEc$}&
                {Lemma 4-8}\cr
$\AAc$ :& {set of the commutative subrings in $\EEc$}&
                {Definition 4-12}\cr
$\DDt, \EEt, \WWt$:
         &{Differential rings over $\CC[[t_1, t_2,\cdots]]$}&
                {Definition 4-19}\cr
$\BMKdV,\quad\BMKdV^\infty$
          :& {Moduli of the KdV hierarchy}&
       {Definition 4-20, Proposition 4-32}\cr
$\fMKdV,g$
          :&$\pi_{\KdV}^{g} :\BMKdV \to {\fMKdV,g}$&
       {Above Proposition 4-27}\cr
$F_g\BMKdV,\quad\BMKdV_g$
          :& {Filter of Moduli of the KdV hierarchy}&
       {Definition 4-26}\cr
$\EEf$ :& {Micro differential ring with coefficient $\CC$}&
                {Definition 4-4}\cr
$\WWf_g,\quad \WWf_{0,1}$
          :& {Gauge freedom}&
       {Lemma 4-28}\cr
$F_g\BMeP$
          :& {Filter of Moduli of quantized elastica}&
       {Proposition 4-29}\cr
$\frak M_{\hyp,g}$
          :&  Moduli of hyperelliptic curves of genus $g$&
       {Proposition 5-4}\cr
$P(X)$   :& {Path space over $X$ in $\Top$}&
       {Proposition 6-2}\cr
$\Omega X$     :& {Loop space over $X$ in $\Top$}&
       {Proposition 6-2}\cr
$\Cal D\CMeP, \quad \Cal C\CMeP$&
Complex related to quantized elastica&
       {Proposition 7-1}\cr
}}
$$

%\newpage
\tvskip
\centerline{\twobf \S 2. A Loop in $\PP$}\tvskip

In this section we will give an expression of a real curve immersed
in a Riemann sphere  $\PP$ following one of Pedit [Ped].
His expression is based upon the oldest theory of a complex curve
embedded in a complex plane $\CC$ or an upper half plane $\Bbb H$,
which was found in ending of the nineteenth
 century and studied by
Klein, Schwarz, Fuchs, Poincar\'e and so on [Po]. Using the expression,
 we will define the moduli space  $\BMeP$ ($\BMeCtwo$)
and $\CMeP$ ($\CMeCtwo$) of smooth curves in $\PP$ ($\CCtwo$)
in Definition 2-10 and 2-12 and an energy functional of a curve
in Definition 2-18,
whose integrand is the Schwarz derivative along the curve. As mentioned
in Introduction, we will call $\CMeP$ a moduli space of a
quantized elastica.

\vskip 1.0 cm

Let us consider a smooth immersion of a circle into two dimensional
complex plane without origin,
$$   \psi : S^1 (:=\RR/\ZZ)\hookrightarrow \CC^2\setminus \{0\},
   \quad  \left(s \mapsto \psi(s) =
\pmatrix \psi_1(s) \\ \psi_2(s) \endpmatrix \right).
$$
Using this map and the natural projection of $\CC^2\setminus \{0\}$
to the
complex projective space (Riemann sphere) $\PP$, we can define the
immersion of a loop in $\PP$:

\tvskip
\proclaim{2-1. Definition}
{\rm
 We define an immersion $\gamma: S^1 \hookrightarrow \PP$ by the
 commutative diagram as $\gamma=\varpi \circ \psi$,
$$\CD           @.  \CC^2\setminus \{0\}   \\
                @.    @VV{\varpi}V \\
                S^1 @>\gamma>> \PP
\endCD .
$$
}
\endproclaim
\tvskip

For a chart around $\psi_2 \neq 0$,
$ s \mapsto\gamma(s)= \dfrac{\psi_1}{\psi_2}(s)$.

\proclaim{2-2. Definition}
{\rm
\roster
\item
The {\it special linear map} $\SL(\CC):\CC^2\setminus \{0\}
\longrightarrow \CC^2\setminus \{0\}$,
$$   m \in \SL(\CC)
   :=\left\{\pmatrix a& b \\ c & d\endpmatrix \ \Big| \
   a, b, c, d \in \CC\quad a d - b c =1 \right\} ,
$$
acts on  $\PP$ through the M{\"o}bius transformation as a symmetric group of
$\PP$: $g_m : \varpi \circ \psi\mapsto \varpi \circ m \psi$ for $ m \in
\SL(\CC) $ and for a point $\gamma \in \PP$,
$$  g_m :\gamma \mapsto \frac{ a \gamma + b}{ c \gamma + d },
    \quad \roman{for }\ m =\pmatrix a& b \\ c & d\endpmatrix.
$$
Let $\PSL(\CC)$ denote  this group including the  group action.

\item
Let  $\Gamma_0(\CC)$ denote the subgroup which is characterized by
vanishing condition of $(2,1)$-component,
$$   \Gamma_0(\CC)
   :=\left\{ \pmatrix a& b \\ 0 & d\endpmatrix \in \SL(\CC) \ \right\} ,
$$
and $\EE_0(\CC)$ denote
 its action to $\CCinf$ using the M{\"o}bius transformation.

\item
Let $\Gamma_1(\CC)$  denote the other subgroup which is characterized by
$$   \Gamma_1(\CC)
   :=\left\{ \pmatrix a& b \\ 0 & d\endpmatrix \in \SL(\CC) \ \Big| \
    |a|=1 \right\} ,
$$
and $\EE_1(\CC)$ denote
 its action to $\CCinf$ using the M{\"o}bius transformation.

\endroster
}
\endproclaim

\tvskip
\proclaim{2-3. Remark}{\rm
$\PSL(\CC)$ has following properties:
\roster
\item Translation, rotation and global dilatation:
$\pmatrix a & b \\ 0 & d\endpmatrix
 \in \Gamma_0(\CC)$, ($b=1/a$)
$$   z \longrightarrow \frac{a z + b}{d} = a^2 z  + a b .
$$
If we restrict the action into $\Gamma_1(\CC)$, it generates
a Euclidean motion induced from $\CC=\PP\setminus \{\infty\}$.

\item Coordinate transformation from chart around $0$ to chart
around $\infty$:
$$
z \longrightarrow -1/z.
$$
\endroster}
\endproclaim

\tvskip

In Definition 2-10,
 we give the definitions of moduli spaces of a quantized elastica
in $\PP$, which are our main objects in this article.
However as Proposition 2-8 is correct for a more complicate
system, we will give provisional moduli spaces of loops.

\proclaim{2-4.  Definition}
{\rm
We define the moduli spaces  of loops as sets as follows:
\roster
\item
$$
   \BMP := \{ \gamma: S^1 \hookrightarrow \PP\ |\
    \gamma \text{ is smooth immersion} \}/\PSL(\CC) .
$$

\item
$$
 \BMCtwo :=\{ \psi:
S^1 \hookrightarrow \CC^2\setminus \{0\}\ |\
     \psi \text{ is smooth immersion},
    \det(\psi(s), \partial_s \psi(s)) =1\ \  \}/ \SL(\CC) .
$$
Here $\partial_s :=d/d s$.

\endroster}
\endproclaim

\tvskip
\proclaim{2-5. Remark}{\rm
\roster

\item
Let $[\gamma]$ denote an element in $\BMP$ for a
representative  element $\gamma \in \PP$ and
 an element in $\BMCtwo$ by $[\psi]$ for a
representative element $\psi \in \CCtwo$.

\item
For a free loop space $\Bbb M$ over a base space M,
$$
	\Bbb M :=\{\delta:S^1 \to M\ | \
               \delta \text{ is smooth immersion} \},
$$
we can define an evaluation map
$\ev$ from $S^1 \times \Bbb M$ to $M$ by
$\ev( s, \delta ) = \delta(s)$ [Br].
For $\Bbb M^\circ$
 ($\circ$ is $\PP$ or $\CCtwo$), we have the evaluation
map whose image is a little bit complicate.

\item
For loops $\psi_1$ and $\psi_2$ in $\CCtwo$ such that
$[\psi_1]=[\psi_2] \in \BMCtwo$, we obviously obtain
 $[\varpi \psi_1] = [\varpi \psi_2]$ in $\BMP$.
Thus we also use the notation of $\varpi$ as
the map,
$$
	\varpi: \BMCtwo \to \BMP.
$$

\endroster
}
\endproclaim

\tvskip

\proclaim{2-6 Definition}{\rm (Schwarz derivative) [Po]
$\{\gamma(s),s\}_{\SD}$ is called {\it Schwarz derivative},
which is defined for
a smooth map $\gamma: S^1 \to \PP$ equipped with a parameter
$s\in S^1$ by,
$$    \{\gamma(s),s\}_{\SD}:= \partial_s \left( \frac{\partial_s^2 \gamma(s)}
      {\partial_s \gamma(s)} \right)-\frac{1}{2}
      \left( \frac{\partial_s^2 \gamma(s)}{\partial_s \gamma(s)} \right)^2   .
$$
We write it by $\{\gamma,s\}_{\SD}$ or $\{\gamma,s\}_{\SD}(s)$
for brevity.
}
\endproclaim

\tvskip
By elementally computations,
the Schwarz derivative is also expressed by
$$ \{\gamma,s\}_{\SD}=
   \left( \frac{\partial_s^3 \gamma(s)}
   {\partial_s \gamma(s)} \right) -\frac{3}{2}
   \left( \frac{\partial_s^2 \gamma(s)}
   {\partial_s \gamma(s)} \right)^2.
$$

Straightforward computations give following lemma.

\vskip 0.5 cm\proclaim{2-7. Lemma} {\rm [Po]}
{\it \roster\item For the action of $g \in \PSL(\CC)$, the  Schwarz
derivative $\{\gamma,s\}_{\SD}$ is invariant:
$$   \{\gamma,s\}_{\SD}=\{g(\gamma),s\}_{\SD} .
$$

\item For a diffeomorphism $s'\in \Diff(S^1)$
$$   \{\gamma,s\}_{\SD}=\left(\partial_s s'\right)^2
   (\{\gamma,s'\}_{\SD}-\{s,s'\}_{\SD})
$$
and for $\roman{U}(1)$ action on $S^1$, {\it i.e.},
$s'=s +\alpha$,
$$
\{\gamma(s),s\}_{\SD}=\{\gamma(s'-\alpha),s'\}_{\SD}
$$
\endroster}
\endproclaim

\tvskip\proclaim{2-8. Definition/Proposition}{\rm [Po]}
{\it
 There is a natural one-to-one correspondence between
$\BMCtwo $ and $\BMP$
 with
the following properties.

 \roster
\item
If $[\gamma]$ is an element of $\BMP$,
there exists a unique lifted curve
 $[\psi]$ as an inverse of the map $\varpi$
$(\varpi [\psi] =[\gamma])$.
Let the correspondence be denoted by $\widetilde\sigma$,
{\it i.e.},
$\widetilde\sigma:\BMP  \longrightarrow
\BMCtwo$, $([\psi]=\widetilde\sigma([\gamma]))$.
Then we have
 $\varpi\circ\widetilde \sigma ([\gamma]) = [\gamma]$
and $\widetilde \sigma\circ \varpi([\psi]) = [\psi]$.

\item For  a map $\gamma:S^1 \to  \PP$
representing a point of $\gamma(s)\in \PP$,
there is a curve $\psi$ in  $\CCtwo$
as a solution of the differential equation,
$$   (-\partial_s^2 - \frac{1}{2}\{\gamma,s\}_{\SD}(s) )\psi(s) =0 ,
$$
so that $\psi$ defines an element $[\psi]\in \BMCtwo$
and  $[\varpi \psi] = [\gamma]$.
This algorithm is a realization of $\widetilde\sigma$.

\endroster}
\endproclaim

\vskip 0.5 cm
\demo {Proof}
In this proof, we will deal only with representative elements $\gamma$ and
$\psi$ of $\BMP$ and $\BMCtwo$.
First we will check the well-definedness of
$\widetilde\sigma$ in (2).
Without loss of its generality, we use the chart of $\psi_2 \neq 0$
a loop $\psi:= \pmatrix \psi_1 \\ \psi_2 \endpmatrix \in
\CCtwo$. Noting the Remark 2-5 (3),
the well-definedness means that the lift of the loop
$\gamma(S^1) :=\varpi \psi(S^1)$ is uniquely $\psi$ up to $\SL (\CC)$.
By differentiating $\det(\psi,
\partial_s \psi) =1$ in $s$,
$(\partial_s^2 \psi_2)/\psi_2 =(\partial_s^2
\psi_1)/\psi_1 $. After straightforward computation, for
$\gamma=\psi_1/\psi_2$, we obtain the relation, $(\partial_s^2
\psi_2)/\psi_2 = -\{\gamma,s\}_{\SD}/2$.
Up to $\SL(\CC)$, $\psi$ is identified with a solution of (2).
Hence well-definedness is asserted.
Further existence of a solution of this equation in (2) is
guaranteed by a special solution,
$ \psi = \pmatrix \sqrt{-1} \gamma/\sqrt{\partial_s
\gamma}\\ \sqrt{-1}/\sqrt{\partial_s \gamma} \endpmatrix$, whose
$\det(\psi,\partial_s \psi)$ is unit. The property of Wronskian
$\det(\psi, \partial_s\psi) =1$ and the uniqueness of the solutions
of a second order differential
equation confirms uniqueness of the solution of (2) up to $\SL(\CC)$.
Further due to the construction of the solutions,
we will consider the effect of $\Diff(S^1)$;
for $s'(s)$, the Schwarz derivative changes as in
Lemma 2-7, and $\partial_s'$ given by the chain rule,
 $\psi(s'):=\psi(s)/\sqrt{\partial_{s'} s}$. Then
$\psi_1(s'):\psi_2(s')=\gamma(s')$ and
$\gamma(S^1)$ and $\psi(S^1)$ do not depend on the parameterization.
Thus (1) and (2) are completely proved.
\qed \enddemo

\vskip 0.5 cm
\proclaim{2-9. Remark} {\rm (Poincar\'e and Schwarz)[Po, Ba1]
By the analytical continuation of $s \in S^1$, $\gamma$ can be
complexfied to $\gamma_c$. If $\gamma_c$ is also embedded in $\CC$,
$\gamma_c^{-1}$ is automorphic function. (In general, even though
$\gamma$ is immersed or embedded in $\PP$, $\gamma_c$ can not be
immersed in $\PP$.) For example the case $s = \wp(\xi)$ ($\xi =
\wp^{-1}(s) \in X_1$) for $s \in \PP$, $\xi \in \CC$,
$\{\xi,s\}_{\SD}$ is a meromorphic function of $s$, where $\wp(\xi)$ is
the Weierstrass elliptic function and $X_1=\CC/\Lambda$ is an elliptic
curve.
These studies are by Klein, Riemann, Poincar\'e, Schwarz and so on.
In this article, we will not restrict ourselves to deal with only with
meromorphic function. We will consider transcendental functions of $s$
because our problem is related to a physical problem or an elastica
problem as the catenary, another physical curve, is also given by the
transcendental function.
}
\endproclaim

\tvskip
In this article,
we are concerned with a loop with an energy functional
in Definition 2-18. However the integrand $\{ \gamma, s\}_{\SD}$
in the energy integration depends upon the parameterization of $S^1$
or $\Diff(S^1)$ from Lemma 2-7. Hence we must fix the
parameterizations of the loop in order to treat
a loop with the energy functional.
Even in $\PP$, we can locally define the metric because
its tangent space $T\PP$ is isomorphic to $\CC$ but
the action of $\PSL(\CC)$ prevents that the metric
becomes global. Hence we restrict ourselves to consider
an action of the subgroup $\Gamma_0(\CC)$
instead  of $\PSL(\CC)$.

Let us introduce our main objects in this article.

\proclaim{2-10.  Definition}{\rm
We define the moduli spaces of loops,
which are called {\it moduli of a quantized elastica}
or {\it moduli spaces of a quantized elastica},
as follows:

 \roster

\item $$
   \BMeP  := \{  \gamma : S^1
           \hookrightarrow \PP \ |\  \text{ smooth immersion},
\ |\partial_s \gamma(s)|=1\  \}/\EE_0(\CC)  .
$$

\item
$$
   \BMeC := \{ \gamma : S^1
           \hookrightarrow \CC \ |  \text{ smooth immersion},
\ |\partial_s \gamma(s)|=1\  \}/\EE_0(\CC)  .
$$

\item
 $$
\split
{\BMeCtwo}
:= \{ \psi:
S^1 \hookrightarrow \CC^2\setminus \{0\}\ |\ &
 \text{s smooth immersion},\ \det(\psi(s), \partial_s \psi(s)) =1,\\
    & \ |\psi_a(s)|=1\
    (a = 1 \text{ or }2)
 \  \}/ \Gamma_0(\CC). \\
\endsplit
$$

\endroster}
\endproclaim

\proclaim{2-11. Remark}
{\rm
\roster
\item The condition $|\partial_s \gamma(s)|=1$ means that we will
treat only loops  with the arc-length parameter $s$
in $\CC$ or $\PP=\CC\cup\{\infty\}$ equipped
with the standard flat metric hereafter.
We call the condition $|\partial_s \gamma(s)|=1$
{\it reality condition} or {\it arc-length condition}.

\item
We continue to express the elements in $\BMeP$, $\BMeCtwo$ by
$[\gamma]$ and $[\psi]$ for loops $\gamma \in \PP$
 and $\psi\in \CCtwo$ satisfying appropriate conditions
respectively for a while.

\item
Further similar to Remark 2-5 (3), we can define
the map from $\BMeCtwo$ to $\BMeP$ by $\varpi$
noting that the reality condition
$|\partial_s \gamma(s)|=1$ means $|\psi_a(s)|=1$ ($a=1$ or $2$)
under the condition $\det(\psi(s),\partial_s \psi(s))=1$
since $\partial_s \gamma(s) =$
$\det(\psi(s),\partial_s \psi(s))/\psi_2(s)^2$ for
$\psi_2(s) \neq 0$.

\item
We can find a representative element by tuning the dilation of
$\EE_0(\CC)$. By letting $\oint |d \gamma| =2\pi$
for a curve with finite length in $\CCinf$
 we have a natural isomorphism,
$$
\split
   \BMeC \approx \BMeCpi:=\{ \gamma : S^1
           \hookrightarrow \CCinf \ |\ &\text{ smooth immersion, }\\
 &| \partial_s \gamma(s)|=1, \
           \oint |d \gamma| =2 \pi\ \}/
\EE_1(\CC)  .
\endsplit
$$
Here $ |d \gamma| =|\partial_s \gamma(s)| ds$.

\item Using (1), we have a decomposition whether the length is
finite or not, {\it i.e.},
$$
\BMeP \approx \BMeCpi\coprod \BMeinf.
$$
This picture is also asserted if one considers a smooth loop
in a two-sphere $S^2$ and its stereographic projection.

For $ \BMeinf$, we have another representation element,
$$
\split
   \BMeinf \approx {\BMeinfc}:=\{ \gamma : S^1
           \hookrightarrow \CCinf \ |\
 & \text{ smooth immersion, } \\
 & | \partial_s \gamma(s)|=1,\
          \sup_{s \in S^1}
|\partial_s \log \partial_s \gamma(s)| =1 \ \}/
\EE_1(\CC)  .\\
\endsplit
$$
Using the equivalences and  such a representative element, we can
introduce scale in our system.

\endroster}
\endproclaim

\tvskip

Next let us introduce smaller moduli spaces for later convenience.
Even under the reality condition $|\partial \gamma(s)|=1$,
there is a freedom to choose its origin of the loop, which is
 denoted by $ \Isom(S^1)=$U(1).
Thus let us define smaller sets with projections to these
sets of moduli spaces,
{\it e.g.}, $\pi_{\elas}^{\PP} :$ $\BMeP
\to {\BMeP}/\Isom(S^1)$.

\tvskip
\proclaim{2-12.  Definition}
{\rm (Moduli of a quantized elastica)
 We define  moduli spaces of loops,
which are also called {\it moduli of a quantized elastica},
or {\it moduli spaces of a quantized elastica},
as follows:

\roster

\item  $\pi_{\elas}^{\PP} :\BMeP \to {\CMeP}:=\BMeP/\Isom(S^1)$.
\tvskip

\item $\pi_{\elas}^{\CC} :\BMeC\to {\CMeC}:=\BMeC/\Isom(S^1)$.

\tvskip
\item $\pi_{\elas}^{\CC^2\setminus \{0\}} :
\BMeCtwo\to
 \CMeCtwo:=\BMeCtwo/\Isom(S^1)$.

\endroster }
\endproclaim

 In physics, we are concerned only
with the shape of elastica. $\Cal M_{\elas}^\circ$ is
more important than $\Bbb M_{\elas}^\circ$.
Further we remark here that we have a natural isomorphism
$\CMeP \approx \BMP/\Diff(S^1)$ as a connection between
$\BMeP$ and $\BMP$.
\tvskip

Next we give our correspondence between $\BMeP$ and $\BMeCtwo$
based on Proposition 2-8

\proclaim{2-13. Definition/Proposition}
{\it \roster
\item"(1)\ \ \ " There is a natural one-to-one and continuous
correspondence between ${\BMeCtwo} $ and ${\BMeP}$
%up to $\roman{U}(1)$ action on $S^1$
 with the following properties.

\item"(1-1)" If $[\gamma]$ is an element of $\BMeP$,
there exists a unique lifted curve
 $[\psi]$ in ${\BMeCtwo}$
as an inverse of the map $\varpi$,
$(\varpi[\psi]=[\gamma])$.
Let the correspondence be denoted by
$$
\sigma: {\BMeP}  \longrightarrow
 {\BMeCtwo},  \quad([\psi]=\sigma([\gamma])).
$$
Then we have
 $\varpi\circ \sigma ([\gamma]) = [\gamma]$
and $ \sigma\circ \varpi([\psi]) = [\psi]$.

\item"(1-2)" For a  curve $\gamma(s) \in \PP$ representing
$[\gamma]\in {\BMeP}$,
there is a curve $\psi$ in  $ {\CCtwo}$
as a solution of the differential equation,
$$   (-\partial_s^2 - \frac{1}{2}\{\gamma,s\}_{\SD}(s) )\psi(s) =0 ,
$$
so that $\psi$ uniquely defines an element $[\psi]$ in
$\BMeCtwo$ and  $\varpi [\psi] =[\gamma]$.
This algorithm is a realization of $\sigma$.

\item"(2)\ \ \ " There is a natural one-to-one correspondence between
$ {\CMeCtwo} $ and ${\CMeP}$
induced from $\sigma$ and $\varpi$.

\item"(3)\ \ \ " $\BMeC$ and $\BMeCtwo$ ( $\CMeC$ and $\CMeCtwo$)
are connected with
$$ \CMeC = \varpi \CMeCtwo, \quad
 \BMeC = \varpi \BMeCtwo.
$$

\endroster
\endproclaim

\demo{Proof}
(1) and (2) are essentially the same as the proves in Proposition 2-8
if we check the compatibility between $|\partial_s \gamma(s)|=1$ and
$|\psi_2(s)|=1$, and continuity of the map.
Due to Remark 2-11 (3), the condition $|\psi_2(s) | = 1$
in $\BMeCtwo$ essentially means the reality
 condition $|\partial_s \gamma(s)|=1$
of the map $\gamma$ on the chart around $\psi_2 \neq 0$.

 Let us consider the continuity.
For a loop in $\PP$, $\gamma'(s):=\gamma(s)+\epsilon v(s)$
with small number $\epsilon$ and an element $v$ of $\Cinf(S^1, \CC)$,
the Schwarz derivative changes
$$
\{\gamma',s\}_{\SD} =\{\gamma,s\}_{\SD}
          +\epsilon\left[
         \frac{ \partial_s^3 v}{\partial_s \gamma}
        -\frac{ \partial_s^3 \gamma}{(\partial_s \gamma)^2}
              \partial_s v
        - 3 \left(\frac{\partial_s^2}{\partial_s \gamma}\right)
         \left(\frac{\partial_s^2 v}{\partial_s \gamma}
         -\frac{\partial_s^2 \gamma}{(\partial_s \gamma)^2}
             \partial_s v\right)\right].
$$
From the proof of Proposition 2-8,
a solution of the its differential equation, $\psi'$, is
periodic and thus is a loop in $\CCtwo$.
For sufficiently small $\epsilon$, the second term becomes small enough.
Then using the perturbation theory, we have $\psi'=\psi+\epsilon \eta$
as a solutions of
$$
\partial_s^2 \psi'+ \frac{1}{2}\{\gamma',s\}_{\SD} \psi' =0.
$$
We note $\{\gamma',s\}_{\SD}$ is also invariant for $\PSL(\CC)$.
For the condition $|\psi_2|=1$, we replace the parameter $s$
by $s'$ using the fact in the proof of Proposition 2-8.
On the other hand, for $\psi'=\psi+\epsilon \eta'$
we can find $v \in \Cinf(S^1,\CC)$ such that
$\gamma':= \psi_1'/\psi_2' = \psi_1/\psi_2 + \epsilon v'$.
Hence both maps are continuous. \qed\enddemo

Here we will consider the natural neighborhood in the
moduli space $\BMeP$.

\proclaim{2-14. Corollary}
\roster
\item
There is a continuous injective map from $\BMeP$
 to the function space $\Cinf(S^1, \CC)$.

\item $\BMeP$ has a natural topology
generated by neighborhood in
the function space $\Cinf(S^1, \CC)$.

\endroster
\endproclaim

\demo{Proof}
(2) is obvious if (1) is proved.
For a given $u \in \Cinf(S^1, \RR)$
the functions $\psi \in \Cinf( [0, 2\pi), \CC^2)$
satisfying
$$
	(-\partial_s^2 - u)\psi =0
$$
is uniquely determined up to $\SL(\CC)$. In general,
even though $u$ is periodic and a function over $S^1$,
$\psi$ is not periodic due to Floquet theorem [MM].
 However if $\psi$ is periodic
for some $u$, by letting $|\psi_2|=1$,
$\gamma \in \PP$ is uniquely determined by
$\gamma=\psi_1/\psi_2$ up to $\EE_0(\CC)$.
We note that for such $\gamma$ and $u$, $u$ is
given as $u = \{\gamma, s\}_{\SD}/2$ and
$\{\gamma, s\}_{\SD}/2$ is invariant for the action of $\EE_1(\CC)$.
(For an action $\Diff(S^1)$ to the reparameterization
of the coordinate $s$, $\{\gamma, s\}_{\SD}/2$
is not invariant and $\{\gamma, s\}_{\SD}/2$
 changes its value. However we have considered
only the arc-length parameterization of $s$ as
 $|\partial_s \gamma(s)|=1$.)
Hence if $\psi$ is periodic
for some $u$, by letting $|\psi_2|=1$,
it determines $[\gamma] \in \BMeP$.
The continuity is obvious from the previous
proposition. \qed
\enddemo

\tvskip
As the injective map in the Corollary
2-14 is a continuous map,
the above neighborhood can be geometrically interpreted
as a neighborhood of $\gamma$.

\proclaim{2-15. Remark}
{\rm
It is known that the free loop space can be a metric space
and has natural topology if the base space is a Riemannian
manifold [Br, McLau].
As we can regard that an element in  $\BMeP$ with finite
length $\int |d \gamma(s)| < \infty$
 can be represented by a loop with
$2\pi$ length whose gravity center exists at the origin of
$\CC$, it is not
difficult to treat the quotient by $\EE_0$ or $\EE_1$.
Consider the image $C$ of the map $\gamma: S^1 \to \PP$;
$C:=\gamma(S^1)$.
For such a loop $C \subset \PP$ whose represents a point
$[C]\in\CMeP$, there is a normal bundle
characterized by the exact sequence of tangent bundle
of $C$ and $\PP$,
$$
	0 \to T C \to T \PP|_C \to N_C \to 0.
$$
Any elements in $T\PP|_C$ are decomposed to
 $T\PP|_C \equiv N_C\oplus T C$.
For a smooth section $v\in \Cinf(C, T \PP|_C)$ of
$T\PP|_C$ over $C$,
and for an infinitesimal real parameter $\epsilon$,
we have $C + \epsilon v$ as a loop in $\PP$.
Here $+$ means the natural addition in the local chart $\CC$
of $\PP$ in the sense of euclidean geometry.
Of course, it is important to check whether such a loop is in
a different point in $\BMeP$ or not but if it is, we can find
an infinitesimal path from $[\gamma]$ to
$[\gamma + \epsilon v_\gamma]$ in $\BMeP$ by letting
$v_\gamma\in\Cinf(S^1, T \PP|_{\gamma(S^1)})$,
$v_\gamma:=v\circ \gamma$.
$|\partial_s (\gamma(s)+ \epsilon v_\gamma(s)|=1$
is not difficult to be treated
by reparameterizing $s$ by $s'$ in primitive sense.
Further even for the case $[\gamma]=[\gamma + \epsilon
v_\gamma]$,
we can regard it as a trivial path.
If $[\gamma]$ and $[\gamma + \epsilon v_\gamma]$ are
different points for an infinitesimal small $\epsilon$,
 we can regard such $[v_\gamma]$ as an element in  a set of
 smooth sections of tangent bundle of $\BMeP$,
$\Cinf(T_{[\gamma]}\BMeP)\equiv \Cinf(\BMeP, T_{[\gamma]}\BMeP)$.

We show that there exist such different points.
From Corollary 2-14,  $\Cinf(T_{[\gamma]}\BMeP)$ is not the
empty set. Let us find an element in
$\Cinf(T_{[\gamma]}\BMeP)$ for each $[\gamma]\in\BMeP$.
As the fiber of $T_p\PP$ is isomorphic to $\CC$, we
define a norm in  $v\in \Cinf(S^1, T \PP|_{\gamma(S^1)})$ by sup-norm.
(In our article, our argument does not strongly
depend upon the norm in $\Cinf(S^1, T \PP|_{\gamma(S^1)})$.)
As it is difficult to define a length in scaleless space,
we might consider an element in $\BMeCpi$
rather than $\BMeP$ due to Remark 2-11.
Consider $[\gamma] \in \BMeCpi$ for $\gamma:S^1 \to\PP$
satisfying the reality condition
$|\partial_s \gamma(s)|=1$ and $\int d \gamma = 2\pi$,
and $v \in\Cinf(S^1, T\PP|_{\gamma(S^1)})$.
Suppose that $\gamma(S^1)+\epsilon v(S^1)$ preserves local and total
length of $\gamma(S^1)$ for sufficiently small $\epsilon$,
{\it i.e.},
$|\partial_s (\gamma(S^1)+\epsilon v(S^1))|=1$ and
$\int d (\gamma+\epsilon v) = 2\pi$ are satisfied.
We  call the deformation as {\it isometric}.
Then we regard $[\gamma+\epsilon v]$ as an element in $\BMeCpi$.
If the vector field $v ( \not\in\{$Euclidean move$\}$)
is the isometric deformation,
$[\gamma+\epsilon v]$ is a different point from $[\gamma]$
in $\BMeCpi$ and thus they are different points in $\BMeP$.
Then we can naturally define a neighborhood around
a loop with  finite length in our moduli space $\BMeP$.

Similarly for an element in $\BMeinf$, we can define the
neighborhood.
For an element $[\gamma]$ in $\BMeinfc$,  we define
its tangent space and velocity $\Cinf(S^1, T \PP|_{\gamma(S^1)})$.
 We can constraint the velocity as an isometric
path which locally preserves the length.
However since  $\BMeinfc$ is defined by
sup-norm in Remark 2-11, for an element $[\gamma] \in \BMeinfc$ and
an isometric path for $v \in\Cinf(S^1, T \PP|_{\gamma(S^1)})$,
$[\gamma+\epsilon v]$ generally does not belong to $\BMeinfc$
even with a sufficiently small $\epsilon$.
However $[\gamma+\epsilon v]$ is in $\BMeP$
and $[\gamma+\epsilon v]$ generates a point in $\BMeinf$ again.
Hence the path is
well-defined by local argument.

Accordingly we can naturally define a neighborhood in
our moduli space $\BMeP$.

Further as we also define a neighborhood in $\BMeCtwo$
using Proposition 2-13, we  define a topology in
our moduli spaces $\BMeP$ and $\BMeCtwo$.
As the topology comes from that in the loop space
we call it {\it topology of loop space} [Br, McLau].

}
\endproclaim

\tvskip
As the topology of loop space is generated by
$\Cinf(T_\gamma\BMeP)$, we will consider an
infinitesimal deformation parameterized by
$t\in [0, \epsilon]$ for a sufficiently small $\epsilon$
in detail.

\proclaim{2-16.  Remark}
{\rm Due to the arguments in Remark 2-15,
we wish to find  one parameter family $[\gamma_t]$ in $\BMeP$ such that
$\partial_t [\gamma_t] $ belongs to $\Cinf(T_\gamma\BMeP)$.

\roster
\item First, we will consider
an  {\it isometric deformation}
 which locally preserves the arc-length of one parameter family of
loops immersed in $\PP$: $\gamma_\circ: S^1\times [0, \epsilon] \to \PP$,
($\gamma_t(s):=\gamma(s,t) \in \PP$) satisfying
$$
  [\partial_t, \partial_s]\gamma_t(s) = 0.
$$
Here $\partial_s :=\partial/\partial s$ and
$\partial_t :=\partial/\partial t$.
We call this condition {\it isometric condition}.
Then if $|\partial_s \gamma_{t=0}(s)|=1$ for $s\in S^1$,
$|\partial_s \gamma_{t}(s)|=1$ for $(s,t) \in S^1\times[0,\epsilon]$.

For $(s,[\gamma_t]) \in S^1 \times \BMeP$, $\partial_s$ acts only on
$S^1$ whereas $\partial_t$ acts only on $[\gamma_t] \in \BMeP$.
Of course the relation $ [\partial_t, \partial_s](s,[\gamma_t]) = 0$
trivially holds. On the other hand,
 for the evaluation map $\ev(s,[\gamma_t])$
as Remark 2-5 (2), the action of $[\partial_t, \partial_s]$ is not
trivial. However by dealing only with the isometric deformation,
we can avoid the noncommutativity between a deformation and
the evaluation map.

\item Let us consider one parameter family of loops immersed in
$\PP$, $\gamma_\circ:S^1\times [0,\epsilon] \to \PP$,
given by a differential equation, which the right hand side
depends upon $\gamma_t$ itself.
First assume that the differential equation is
 $\partial_t \gamma_t = f(\gamma_t)$ for a given functional $f$.
In this case, the deformation depends upon the affine
coordinate $\gamma_t$ in $\PP$
and it is not invariant for the action of $\EE_1$.
Further we note that $\partial_s \gamma_t(s)$ is the tangential
vector of the circle $\gamma_t(S^1)$ and
$\phi:= \log \partial_s \gamma_t(s)$ denotes its
tangential angle if $|\partial_s \gamma_t(s)|=1$.
Provided that the deformation $ \partial_t \gamma_t$
is governed by a function of $\partial_s \gamma_t(s)$ itself,
the deformation must
depend upon the angle of $\CC$ in $\PP$ and a euclidean move.
Hence they can not be  deformations in $\BMeP$
and a deformation in $\BMeP$ does not include $\gamma(s)$
and $\phi$.

\item From Lemma 2-8 (1),
$u_t(s)\equiv$ $u(s,t):= \{\gamma_t(s), s\}_{\SD}/2$
 is a function of $\BMeP$. Further
$u_t(s)$ depends only on $\partial_s \log(\partial_s \gamma_t(s))$
and $\partial_s^2 \log(\partial_s \gamma(s))$ due to Definition 2-6.
We might consider the deformation in an element $\gamma_t$ in
$\BMeP$ through an equation,
$$
\partial_t u_t =f(u_t, \partial_s u_t, \partial_s^2 u_t, \cdots, A),
$$
for  appropriate functional $f$ and function $A \in
\Cinf( S^1\times [0,\epsilon] , \CC)$;
the function $A$ must be invariant
for the action of $\EE_0(\CC)$.
If $u_t$ is determined at a time $t$, $\gamma_t$ can be
reconstructed by Proposition 2-13.
If there does not appear $\gamma_t(s)$ or $\partial_s \gamma_t(s)$
themselves in right hand side, the deformation is invariant
for the action of $\EE^1(\CC)$ for an appropriate $A$.
\endroster

Due to the above consideration, we can consider an infinitesimal
deformation in $\BMeP$ by $\partial_t u
=f(u, \partial_s u, \partial_s^2 u, \cdots, A)$ and
$[\partial_t, \partial_s]\gamma(s,t)=0$.

}
\endproclaim
\tvskip

As we prepared the tools, it is not difficult to deal with
the quotient space $\BMeP$ and $\BMeCtwo$.
From here, let $\gamma$  itself
denote an element of $\BMeP$
 instead of $[\gamma]$
 for a loop $\gamma(S^1)\in \PP$ satisfying the
conditions. Similarly we write $\psi\in \BMeCtwo$
for the sake of simplicity.
Further we will consider a {\it flow} in $\BMeP$.

\proclaim{2-17. Remark}{\rm
Let us consider the situation that
for a point $\gamma \in \BMeP$ and its neighborhood $U_\gamma$
in terms of  the loop topology,
 we can find another point $\gamma' \in U_\gamma$
such that $\gamma' = \gamma+ \epsilon v$ for a sufficiently
small $\epsilon$ and some velocity
$v \in \Cinf(T_\gamma\BMeP)$.
Suppose that by sequentially finding  such points, we construct a curve
$\gamma_t$, $t\in [0,1]$ in $\BMeP$ connecting between
 a starting point $\gamma_0=\gamma$ and a terminal point
$\gamma_1 =\gamma''$ for some $\gamma'' \in U_\gamma$.
Then we may write the velocity as
$\partial_t \gamma_t$ at $\gamma_t$.
In this way, for each point $\gamma$ in $\BMeP$, we can define
an immersion of $[0,1]$ in $\BMeP$
 for a smooth section  $v \in \Cinf(T_\gamma\BMeP)$ if
it is well-defined.
We call such a immersion {\it flow} in $\BMeP$.

Further for each a point $\gamma_t$  in a flow,
$t\in[0,1]$ with $v_t\in \Cinf(T_{\gamma_t}\BMeP)$, let us
assume that we can choose another
element $v_t'\in \Cinf(T_{\gamma_t}\BMeP)$ and find
a point $\gamma_{t, \epsilon}$ in the neighborhood of $
\gamma_{t}$ such that
$\gamma_{t, \epsilon}=\gamma_t+\epsilon u'$ for
a sufficiently small parameter $\epsilon$. Then we
can consider duplex flow such as $\gamma_{t,t'}$
for $[0,1]^2$ in $\BMeP$.
Similarly we can deal with an immersion of
$[0,1]^m$ in $\BMeP$.
For the case, we  call the immersion $\gamma_{t}$ of $t\in
[0,1]^m$
{\it multiple flow}.
Further for a certain case, $[0,1]^m \in \BMeP$ can be
extended to $\RR^m \in \BMeP$ where $m$ is a positive
integer or the infinite number.

Similarly we can deal with flow in $\BMeCtwo$.
We  define the KdV and KdVH flow as an
extension of $[0,1]^m$ immersion to $\RR^m$
immersion in \S 3.
}
\endproclaim

\tvskip
\proclaim{2-18. Definition}{\rm ( Energy of a quantized elastica)
We introduce an energy functional
 of  $\gamma \in \CMeCpi\approx
\CMeC $, called
{\it Euler-Bernoulli energy functional}, by
$$
 {\Cal E[\gamma]}:=\frac{1}{2\pi}\int_{S^1} \{\gamma(s), s\}_{\SD} d s .
$$

}
\endproclaim

\proclaim{2-19. Lemma}
{\it  For $\gamma\in \CMeP$,
 the energy ${\Cal E[\gamma]}$ is non-negative real.}
\endproclaim

\vskip 0.5 cm
\demo {Proof} The Schwarz derivative can be expressed by
$$   \{\gamma, s\}_{\SD} =
   \partial_s^2 \log(\partial_s \gamma) - \frac{1}{2}
   (\partial_s \log(\partial_s \gamma) )^2 .
$$
Due to Definition 2-10, the reality condition
 $|\partial_s \gamma(s)|=1$, we let
$\partial_s \gamma= \exp(\sqrt{-1} \phi)$, $\phi$ is a real
smooth function over $S^1$, $\phi(0)=\phi(2 \pi)$.
Hence
$$ \int_{S^1} \{\gamma, s\}_{\SD}d s= \int_{S^1} d s
    \frac{1}{2}(\partial_s \phi )^2,
$$
which is real.
\qed \enddemo

\tvskip\proclaim{2-20. Remark}{\rm  [Mat2]
\roster

\item
By Lemma 2-7, the integrand in the energy $\Cal E$ is invariant
for the action of $\PSL(\CC)$.
However the diffeomorphism of $S^1$, $\Diff(S^1)$, changes
the energy. Hence we cannot find a well-defined energy
over $\CMPP$.

Further for $\gamma \in \pi_\elas^\PP\BMeinfc$, we can also consider
a correspondence
$$
\int_{S^1} \{\gamma, s\}_{\SD} d s,
$$
by giving up to considering dilatation symmetry.
However it is neither well-defined for the dilatation.

As we wish to neglect the problem for $\pi_\elas^\PP\BMeinf$,
 we restrict ourselves to deal with $\CMeC$.
In other words, we will consider the energy functional
only for $\CMeC$.

\item We regard the energy function as a section of line bundle over
$\CMeC$,
$$\CD       \RR  @>>> \text{Energy}(\CMeC)  \\@.  @VVV \\
                   @.     \CMeC   \\.
\endCD$$

\item
As mentioned in Introduction, for $\gamma \in \CMeC$,
this energy functional  $\Cal E=\int_{S^1} \{\gamma, s\}_{\SD} d s$
is identified with $\int_{S^1} (\partial_s^2
\gamma/ \partial  \gamma)^2 d s$; thus we call it
 Euler-Bernoulli energy functional.
The  stationary points of $\Cal E$ in  $\CMeP$ in the meaning
(1) were investigated by Euler [E, L, T1, 2].
Even though we will not touch the problem in this paper,
we are implicitly considering the partition function of an elastica
as a problem of physics [Mat2],
$$   Z = \int_{\CMeC} D \gamma \exp\left( -\beta \int_{S^1}
   \{\gamma, s\}_{\SD} \dd s \right) .
$$
In order to know this partition function
(which is not mathematically still well-defined),
we must investigate the moduli space of curve
$\CMeC \subset \CMeP$
and we will do in this paper.

\endroster
}
\endproclaim

%\newpage
\tvskip
\centerline{\twobf \S 3. KdVH flow }
\tvskip

Our studies are based upon the discovery of
Goldstein and Pertich [GP1, 2]
on the MKdV flow for a loop in $\CC$ and that of Langer and Perline [LP]
on the nonlinear Schr{\"o}dinger flow for a loop in $\Bbb R^3$.
Using their results, one of authors studied the moduli of loops in $\CC$
[Mat2] and loops in $\Bbb R^3$ [Mat4].
Our purpose is to give mathematical implications of these works
[Mat2, 4] using results of Pedit [Ped].
In this section,  we will give our main theorem 3-4 and its proof,
which are of a relation between the moduli of a quantized elastica
in $\PP$ and the KdV flow.

\tvskip

In order to express the system of the KdV equation,
we will introduce the differential algebra and its division
algebra before our main arguments in this section.

\proclaim{3-1. Definition}{\rm [SN]
\roster
\item The {\it differential ring} $\DDs$ is defined by,
$$  \DDs:=\{ \sum_{k\ge 0}^N a_k(s)\partial_s^k \ |
    \ N < \infty, \  a_k(s)\in \Cinf(S^1,\CC), s \in S^1 \}.
$$

\item The {\it degree} of a differential operator, $D \in \DDs$,
is denoted by $\roman{deg} D$,
$$
\roman{deg}: \DDs \longrightarrow \Bbb Z_{\ge0},
$$
 where $ \Bbb Z_{\ge0}$ is the set of non-negative integers.

\item The {\it micro-differential ring}   $\EEs$ to  $\DDs$ is
defined by
$$
    \EEs:=\{ \sum_{k= -\infty}^N a_k(s)\partial_s^k \ |
   \ N < \infty, \  a_k(s) \in \Cinf(S^1,\CC),\  s \in S^1 \},
$$
where  $\roman{deg}: \EEs \longrightarrow \Bbb Z$ and
the product  is defined by the extended Leibniz rule,
$$
  \partial_s^n a  = \sum_{r=0}^\infty \pmatrix n \\ r \endpmatrix
         (\partial_s^r a)\partial_s^{n-r} ,
          \quad
 \pmatrix n \\ r \endpmatrix:=\frac{1}{r !} n (n-1)\cdots(n-r+1) .
$$

\item The projections $+$ and $-$ are defined by
$$   +: \EEs \longrightarrow \DDs , \quad ( L \mapsto L_+),
   \quad   -:\EEs \longrightarrow \EEs\setminus \DDs,
   \quad (L \mapsto L_-, \quad L = L_++L_-) .
$$

\endroster}
\endproclaim

\tvskip
Hereafter we will write a map from $S^1$ to $\PP$,
$\BMeP$ and $\CMeP$ by the same $\gamma$.
Noting Remark 2-17,
let us define the KdV and KdVH flows,
which satisfy the isometric condition as
in Proposition 3-11.

\proclaim{3-2.  Definition}{\rm (KdV flow and KdVH flow)
\roster
\item The {\it KdV flow} is
defined as the immersion

$ \gamma_\circ:\RR \hookrightarrow \BMeP$ and
$\psi_\circ: \RR \hookrightarrow
\BMeCtwo$, $(t \mapsto (\gamma_t, \psi_t) )$,

 which satisfies the following properties:

1-1) $  \gamma_t(s)=\varpi \circ \psi_t(s)$,  for each $t\in \RR$.

1-2)
$u(s,t):=\{\gamma_t(s), s\}_{\SD}/2$ obeys the KdV equation,
$$     \partial_t u + 6 u \partial_s u + \partial_s^3 u =0 .
$$

If for a point $\gamma\in \BMeP$ and its corresponding point
$\psi\in \BMeCtwo$, there is one of the  KdV flows such that
$\gamma_t = \gamma$ and $\psi_t =\psi$
for some $t \in \RR$, we say that $\gamma$ or $\psi$
{\it belongs} to  the KdV flow $\gamma_\circ$ or $\psi_\circ$.

\item
Let us introduce a formal infinite dimensional parameter space,
$$   {\Cal V}^\infty := S^1 \times (\prod_{n=1}^\infty \Bbb R),
  \quad t=( t_1,  t_2, t_3, \cdots) \in {\Cal V}^\infty, \ t_1 \in S^1.
$$
Then the {\it KdVH flow} is defined as the immersion
$$
(\gamma_\circ, \psi_\circ)
\equiv \overline{\phi_{\partial_s u,t}}
    : \Cal V^{\infty} \hookrightarrow
      \BMeP\times\BMeCtwo,\quad (t \mapsto (\gamma_t, \psi_t) ),
$$
which satisfy the  following conditions:

1-1) $ \gamma_t(s)=\varpi \circ \psi_t(s)$,

1-2)
 $\overline{\phi_{\partial_s u,t}}$ is given by
$$   \gamma(s,t) \longrightarrow   \gamma(s,t+\delta t)=
   \exp( \sum_{n=1}\delta t_{n} \partial_{t_n} )\gamma(s,t) ,
$$
whose each $t_n$ deformation obeys the $n$-th KdV equation ($n \ge 1$),
$$   \partial_{t_{n}}u=-\Omega^{n-1}\partial_s u ,
$$where  $\Omega$ is a micro-differential operator,
$$   \Omega =  ( \partial_s^2 + 2 u  + 2\partial_s u \partial_s^{-1})
 \in \EEs.
$$

If for a point $\gamma\in \BMeP$ and its corresponding point
$\psi\in \BMeCtwo$,  there is  one of the  KdV flows such that
$\gamma_t = \gamma$ and $\psi_t =\psi$
for some $t \in{\Cal V}^\infty $, we say that $\gamma$ or $\psi$
{\it belongs} to the KdVH flow $\gamma_\circ$ or $\psi_\circ$.

\item We define a  relation,
$$   \gamma \simKdVH \gamma',
$$
for two points $\gamma, \gamma' \in  \CMeP$
if these $\gamma$ and $\gamma'$ are on an orbit of the projection of the
KdVH flow $\pi^{\PP}_{\elas} \circ \overline{\phi_{\partial_s u, t}}$,
{\it i.e.},   every points in the fibers ${\pi^{\PP}_{\elas}}^{-1}\gamma$
and ${\pi^{\PP}_{\elas}}^{-1}\gamma'$ belongs the same KdVH flow.
\endroster}
\endproclaim

For convenience, let $\gamma_t$ ($\psi_t$)
denote the KdV flow or KdVH flow
instead of $\gamma_\circ$ ($\psi_\circ$) from this.

\tvskip \proclaim{3-3}
{\rm \roster

\item Though the well-definedness of the above definition is
  later investigated, these flows satisfy the isometric
  condition as in
Proposition 3-11.

\item We will note that the space ${\Cal
 V}^\infty$ has an algebro-analytic
 structure induced from the equations,
$$   \partial_{t_{n+1}}u=-\Omega \partial_{t_{n}}u.
$$

\item The $n=2$ KdVH flow obeying
$   \partial_{t_{2}}u=\Omega \partial_{s}u  $
is identified with the KdV flow in (1) of
Definition 3-2. \endroster
\endproclaim
\tvskip

Our first main theorem is as follows:

\proclaim{3-4. Theorem}
{\it \roster
\item The relation $\simKdVH$ in the Definition 3-2 becomes an
equivalent relation; for arbitrary $\gamma$ in $\CMeP$ there is
one of the KdVH flows to which $\gamma$ belongs,
and for $ \gamma \simKdVH \gamma'$ and $\gamma'
\simKdVH \gamma^{\prime\prime}$, we have  a relation $\gamma
\simKdVH \gamma^{\prime\prime}$. By this relation, we can define an
equivalent class
$$   \frak C[\gamma]:=\{\gamma' \in \CMeP \ | \
   \gamma' \simKdVH \gamma\},\quad
   \CMeP =  \coprod_{\gamma} \frak C[\gamma].
$$
Similarly we can define
$$  \frak C_{\CC}[\gamma]:=\{\gamma' \in \CMeC \ | \
   \gamma' \simKdVH \gamma\},\quad
   \CMeC =  \coprod_{\gamma} \frak C_{\CC}[\gamma].
$$

\item The KdVH flow  conserves the energy
$\Cal E$.   In other words, for the subspace of $\CMeC$,
$$ \Cal M_{\elas,E}^{\CC} :=\left\{ \gamma \in \CMeC
   \Big|\ \Cal E[\gamma]-E=0\ \right\},
$$
and a curve $\gamma \in \CMeC$,
the following relation holds
$$
      \Cal M_{\elas, \Cal E[\gamma]}^{\CC}
   \supset \frak C_{\CC}[\gamma] .
$$
\item The moduli space
of a quantized elastica $\Cal M_{\elas}^{\CC}$ is
decomposed as
$$  \Cal M_{\elas}^{\CC} =  \coprod_E \Cal M_{\elas,E}^{\CC} ,
\quad \Cal M_{\elas,E}^{\CC} =  \coprod_{\gamma, \Cal E[\gamma]=E}
 \frak C_{\CC}[\gamma] .
$$

\endroster}
\endproclaim

\tvskip

Noting Remark 2-16, we will investigate the moduli spaces
$\BMeP$ and $\BMeCtwo$ by considering flows over there
and prove our theorem.
Here we  mention the strategies of the
proof of the theorem.

\proclaim{3-5}{\rm
We plane to investigate $\BMeP$ by dealing with a group
which is generated by a Lie algebra associated with $T\BMeP$.
By the correspondence between $\BMeP$ and
$\BMeCtwo$ in Proposition 2-8, we can identify
$\gamma(s)$ with $(\psi_1, \psi_2)(s)$.
We  firstly deal
with wider class of flows $\phi_{A,t} $
in Lemma 3-6, which is characterized by
a smooth function $A$ over $S^1\times[0,1]$.
In Lemma 3-9, we  find that an arbitrary flow $\phi_{A,t} $
approximately preserves the energy of elastica in Definition
2-18. Due to the argument in Remark 3-10, we
choose a special $A$ as $A=\partial_s \{\gamma, s\}_{\SD}$
and then the  flow is identified with
the KdV flow in Proposition 3-11.
 As shown in Proposition 3-15, 16, and 17,
we use the regular properties of the KdV hierarchy
and prove the theorem.
}
\endproclaim

\tvskip

Noting Remark 2-16, we have the following lemma.
\proclaim{3-6. Lemma}{\rm (Goldstein-Pertich, Pedit)[GP1, GP2, Ped]}
{\it
Let us consider a flow of $[0, \epsilon]$ for
a real number $\epsilon >0$:
$$   [0,\epsilon] \to \BMeP, \quad
   ( t \mapsto \gamma_t) ,
$$
{\it i.e.}, it is realized by an isometric deformation,
$$   [\partial_s, \partial_t ] \gamma_t(s)=0 .
$$

\roster
\item
Every isometric deformation $\gamma_t(s)$ locally
obeys the equation of motion,
$$   \partial_t u = -\Omega A(s,t) ,
$$
where  $u= \{\gamma, s\}_{\SD}/2$ and $A(s,t)$ is an
appropriate smooth function over
$(s,t) \in S^1\times [0,\epsilon]$.

\item For the  function $A(s,t)$, there exists a smooth function
$B(s,t)$ such that $A(s,t)=-\partial_s B(s,t)/2$ and this equation of
motion is locally rewritten by,
$$   \partial_t u = \frac{1}{2}\underline \Omega B(s,t),
$$
where $\underline \Omega:=\Omega \partial_s$,
$$        \underline \Omega =  ( \partial_s^3
   + 2 u \partial_s + 2\partial_s u ) .
$$
\endroster}
\endproclaim

\demo {Proof} Using the one-to-one correspondence between $\BMeP$
and ${\BMeCtwo}$, we lift the flow
$\gamma_t$ to $\psi_t := \sigma \gamma_t$.
In this proof, we  consider representative elements
of the image of its evaluation map, $\gamma_t(s)$ and $\psi_t(s)$.
 Due to the linear
independence given by $\det(\partial_s \psi_t(s),\psi_t(s))=1$,
we  express the deformation in terms of
$\psi_t(s)$ and $\partial_s \psi_t(s)$;
$$   \partial_t \psi_t(s)
          =( A(s,t) + B(s,t) \partial_s) \psi_t(s),
$$
where $A(s,t)$ and $B(s,t)$ are smooth functions over $(s,t)$. However
from $\partial_t \det(\psi_t(s), \partial_s \psi_t(s) ) =0$,
we have the constraint,
$$   \partial_s B(s,t) = -2 A(s,t) ,\tag 3.1
$$
using $[\partial_s, \partial_t ] \psi_t(s)=0$.
Noting $u(s,t)=-(\partial_s^2 \psi_2(s))/\psi_2(s)$,
we perform a straightforward
computations of $\partial_t u(t, s)$, we obtain the equation in (1).
On the other hand,
if the equation is satisfied, we can reduce the
equation to $[\partial_t, \partial_s]\gamma_t(s)=0$.
Similarly we obtain (2). \qed \enddemo

 \tvskip

Let us introduce another formal infinite dimensional parameter spaces, $
\quad t=( t_1,  t_2, t_3, \cdots) \in   [0,\epsilon]^\infty$
and  a formal multiple flow $\phi_{A, t}$ with the infinite dimensional
parameters, which is locally defined.

\vskip 0.5 cm
\proclaim{3-7 Definition}
{\rm
For $t \in [0,\epsilon]^\infty$ for a sufficiently small parameter
$\epsilon$, we will
define an infinitesimal multiple flow,
$$     \phi_{A,t}:[0,\epsilon]^{\infty} \longrightarrow \BMeP,
     \quad ( t \mapsto \gamma_t) ,
$$
induced from the formal variation for a sufficiently small $\delta t$,
$( \delta t_i < \epsilon < N \delta t_i $, a small natural
number $N$) and  image of evaluation map
$ \gamma(s,t):=\gamma_t(s)$,
$$     \gamma(s,t) \mapsto \gamma(s,t+\delta t)=\exp(\sum_{n=1}
       \delta t_{n} \partial_{t_n} )\gamma(s,t) :=(1+\sum_{n=1}
       \delta t_{n} \partial_{t_n} ) \gamma(s,t) + \Cal O( \delta t^2),
$$
with local relations,
$$   [\partial_s, \partial_{t_n}]\gamma(s,t)=0, \quad (n \ge 1),
$$ $$   \partial_{t_{n}}u=-\Omega^{n-1} A(s,t), \quad (n \ge 1),
$$
where  $u(s,t)= \{\gamma_t(s), s\}_{\SD}/2$,
$A(s,t)$ and $B(s,t)$ are appropriate
 smooth functions over $S^1 \times [0,\epsilon]^\infty$
such that $2 A = - \partial_s B$.
}
\endproclaim

\vskip 0.5 cm
\proclaim{3-8 Remark}
{\rm\roster

\item In terms of the definition of the exponential
 function to the base $e$,
$$
	\exp(O) = \lim_{N'\to \infty} \left(1+\frac{O}{N'}\right)^{N'},
$$
the development of $\delta t_n$ generates $[0,\epsilon]^\infty$.
By tuning $N'$ compatible to $N$ in Definition 3-7, we can define
the exponent action to $\gamma(s,t)$.

\item If $\Omega^{n-1} A(s,t)$ vanishes for $n > M$ for a natural
number $M$, the deformation is of finite dimensional. Then
the flow $\phi_{A, t}$ is well-defined for a sufficiently small
$\epsilon$.

\item In general, the above flow $\phi_{A, t}$ is a formal one and
its well-definedness is not guaranteed.
 However if it is well-defined, it gives an isometric
deformation of a curve $\gamma(s) \in \BMeP$. In fact due
to the relation $\partial_{t_{n+1}} = \Omega\partial_{t_{n}}$, we have
the flow
$$   \partial_{t_n} \psi_t(s) = (A_n + B_n \partial_s )\psi_t(s) ,
$$
where $A_2 = A= -\partial_s B/2$, $B_2 = B$, $A_1=\Omega^{-1} A$, and
$$   A_n = \Omega A_{n-1}, \quad
       \quad \partial_s B_n = \Omega B_{n-1}, \quad (n \ge 2) .
$$
Then the above relation $\partial_{t_{n}}u=-\Omega^{n-1} A(s,t)$
turns out
to be the standard type of the flow for $A_n$ in Lemma 3-6.
}
\endproclaim

\tvskip \proclaim{3-9. Lemma}
{\it\roster
For $\gamma\in\BMeC$ and
$A\in \Cinf(S^1\times [0,\epsilon]^\infty, \CC)$,
 the infinitesimal flow $\phi_{A,t}$
preserves the energy functional modulo $(\delta t)^2$;
$$  \frac{1}{2\pi} \int_{S^1}\{\gamma_t, s\}_{\SD}\ d s =
     \frac{1}{2\pi}
 \int_{S^1}\{\gamma_{t+\delta t}, s\}_{\SD}\ d s +
      \Cal O((\delta t)^2) .
$$
\endroster}
\endproclaim

\vskip 0.5 cm
\demo {Proof} Noting Remark 3-8 and by Eq.(3.1)
in the proof of Lemma 3-6, we have the relations,
$$   \partial_s B_n = -2 A_n = -2 \Omega A_{n-1},
        \quad   \partial_s B_n = \Omega B_{n-1}.
$$
When we will apply the relation to the right hand side of the lemma,
($u(s,t):=\{\gamma_t, s\}_{\SD}/2$),
$$ \split
   \int_{S^1}u(s, t+\delta t) d s
   &=\int_{S^1}u(s,t) d s+\sum_{n=1}\delta
   t_n\int_{S^1}\partial_{t_n} u(s,t) d s +
   \Cal O((\delta t)^2)\\
   &=\int_{S^1}u(s,t) d s-\sum_{n=2}\delta t_n\int_{S^1}
   \Omega A_n d s +\frac{1}{2}\int_{S^1}\partial_s B d s
   +\Cal O((\delta t)^2)\\
   &=\int_{S^1}u(s,t) d s+
   \frac{1}{2} \sum_n\delta t_n\int_{S^1}
    \partial_s B_{n+1}(s,t) d s +
   \Cal O((\delta t)^2)\\
   &=\int_{S^1}u(s,t) d s+
   \Cal O((\delta t)^2).
\endsplit $$
We completely prove the lemma. \qed \enddemo

\tvskip
\proclaim{3-10. Remark}
{\rm
\roster
\item This flow $\phi_{A,t}$ could be regarded as
an infinitesimal action of a diffeomorphism of $\BMeP$,
which is a (infinite dimensional) Lie group $G_A$
if it can be well-defined.

\item
We can regard $S^1$ as a Riemannian manifold with a metric $d s^2$.
Then $\partial_s$ is the Killing vector
and $ \exp(\sqrt{-1}s) $ is a geodesic flow. They are
a generator and an element of the $\Isom(S^1)=$U(1) group  respectively;
$$   \text{\rm{U}}(1):S^1 \longrightarrow S^1, \quad
   ( \exp(\sqrt{-1}s) \mapsto \exp(\sqrt{-1}(s+s_0) ) ,
$$
For $g_0 \in $ U(1), $g_0$ gives a natural automorphism of
$\BMeP$.

\item
Since there is the natural projection $\pi_{\elas}^{\PP}:
 \BMeP \to \CMeP$,
the U(1) action on $\gamma\in \CMeP$
must be trivial $g_0\gamma = \gamma$
for $g_0 \in $U(1) and we have the relation
$g_0 \circ \pi_{\elas}^{\PP} = \pi_{\elas}^{\PP}\circ g_0$.
It implies that the immersion of the  loop $S^1$ is
consistent with U(1) action.

\item
For a curve  $\gamma(s) \in \BMeP$,
we can locally express the  U(1) action,
$$   (\partial_s -\partial_{s_0})\{ \gamma, s\}_\SD (s,s_0)=0, \quad
    (\partial_s -\partial_{s_0})\gamma(s,s_0) =0.
$$
These equations faithfully represent the U(1) symmetry or translation,
$\gamma(s) \to \gamma(s-s_0)$ in $\BMeP$.

\item
Due to the above remarks, if exists,
$G_A$ should include $G_0=$U(1) as
its normal subgroup.
Accordingly it is natural that $A$ in Definition 3-7 starts with
the internal symmetry: $A=\partial_s u$ and
$\partial_{t_1}u = \partial_s u$
for $u = \{ \gamma, s\}_\SD$.

\item When we  consider the multiple flow generated by
$\phi_{\partial_s u,t}$ ($A = \partial_s u$), it means that
we deal with the variation,
$$   \gamma(s,t) \longrightarrow \gamma(s,t+\delta t)=
   \exp( \sum_n \delta t_{n} \partial_{t_n} ) \gamma(s,t) ,
$$
which obeys
$$  \partial_{t_{n}}u=-\Omega^{n-1}\partial_s u .
$$
Following Definition 3-2, they are locally
 identified with the KdVH flow.

\item Due to the Remark 2-16, this multiple flow is
locally well-defined in
$\BMeP$.

\item
Physically speaking for the above arguments, we are implicitly
investigating a partition
function of a {\lq\lq}elastic" curve in $\PP$. We require that
the partition function must naturally include classical shapes
 whose
have the above trivial translation symmetry as the Goldstone
 bosons or the Jacobi fields [R].
This requirement  makes the group structure
acting $\CMeP$ (if exists) contain this trivial
symmetry [Mat2].
\endroster
}\endproclaim

We will summarize the above results as a proposition.
\tvskip
\proclaim{3-11. Proposition}
{\it \roster

\item The multiple flow $\phi_{\partial_s u,t}$  contains a subflow
$\phi_{\partial_s u,t_1}$ generated by
 $$   (\partial_s -\partial_{s_0})\{ \gamma, s\}_\SD =0.
 $$
This domain of $t_1\in [0,\epsilon]$ is extended to $S^1$
and  is consistent with the projection $ \pi_{\elas}^{\PP}:\BMeP
\to \CMeP$, {\it i.e.}, there exists $\varphi_{\partial_s u,t}$
such that
$ \pi_{\elas}^{\PP}\circ\phi_{\partial_s u,t}
=\varphi_{\partial_s u,t}\circ  \pi_{\elas}^{\PP}$.

\item
By choosing $A=\partial_s u$ for $u= \{\gamma, s\}_{\SD}/2$,
the flow $\phi_{\partial_s u,t}$ defined in Definition 3-7
 is well-defined as a flow in $\BMeP$
and can extend the domain of the flow
$[0,\epsilon]^\infty \to \Cal V^\infty$.

\item
$\phi_{\partial_s u,t}$
 is identified with the KdVH flow
$\overline{\phi_{\partial_s u,t}}$ by extending
$[0,\epsilon]^\infty$ to $ \Cal V^\infty$.

\item  There exists a flow $\overline{\varphi_{\partial_s u,t}}$
 in $\CMeP$ such that
$\pi_{\elas}^{\PP} \circ \overline{ \phi_{\partial_s u,t}}
 =\overline{ \varphi_{\partial_s u,t}}\circ \pi_{\elas}^{\PP} $,
we also call it  KdVH flow.

\item For the KdVH flow, we have algebraic relations among multi-times
$t_n$ as $ \partial_{t_{n+1}}u=\Omega \partial_{t_{n}}u$.

\item The KdVH flow preserves the decomposition in Remark 2-11.

\item The restricted flow of the KdVH flow
to $\BMeC$  preserves the energy
functional exactly.

\endroster }
\endproclaim

\demo{Proof}
(1) is obvious from Remark 3-10.
If (2) is satisfied,
(3), (4) and (5) are naturally given from Remark 3-10.
Since the KdVH flow consists of isometric deformations, (6) is obvious.
(2) and (7) will be asserted by Propositions
 3-15 and 3-16.
\qed
\enddemo

Firstly we note that (7) should be compared with the Lemma 3-9.
Next we also note that in order to prove (2), we should check
1) the well-definedness of the KdVH flow locally
and 2) the extension of the domain to $\Cal V^\infty$.
If the well-definedness of the KdVH flow is guaranteed,
we can find the neighborhood of a point $\gamma \in \BMeP$
by the KdVH flow to $\gamma$
as its initial state, because the KdVH flow
consists of the isometric deformations.
We can consider the process in
$\BMeCpi$ as mentioned in Remark 2-15.

\tvskip Here we will introduce the words of
a dynamic system here apart
from our notations in main subject [AM].

\proclaim{3-12. Definition}{\rm [AM, Br]}
{\rm We will consider a manifold $M$ equipped with a closed real
2-form $\omega$. We will use the notations: $i_Y v$ is the
interior product
of a vector field $Y$ and a differential form $v$.
\roster
\item A vector field $Y$ is called
         {\it symplectic} if $i_Y \omega$ is closed.

\item A vector field $Y$ is called a
{\it Hamiltonian vector field} if there exists a
function $f$ such that $i_Y \omega= d f$. \endroster}
\endproclaim

\vskip 0.5 cm
Corresponding to Definition 3-12, we will define quantities  in the KdV
flow in Definition 3-13 and give Proposition 3-15 by assuming
$\BMeP$ as a (infinite dimensional) manifold [AM].

\vskip 0.5 cm \proclaim{3-13. Definition}{\rm [AM]}
{\rm\roster \item In our KdVH flow,  we define a 2-form $\omega$
for vectors $Y_1$ and $Y_2$ over $\BMeP$,
$$   \omega(Y_1, Y_2) := \frac{1}{2} \int_{S^1}
   \left(\int_0^s (Y_2(s)Y_1(s')-Y_1(s) Y_2(s') ) d s'\right)
   d s .
$$

\item We define the quantities $X_n$ and $h_n$ and variation
$\overline
\delta/\overline \delta u$ for the KdVH flow: $h_0 =u/2$, $X_0=0$ and
$$   X_n(u) := \Omega^{n-1}\partial_s u, \quad
   X_n(u)=\partial_s \frac{\overline \delta h_n}{\overline \delta u},
$$ where $$
      \frac{\overline \delta h_n}{\overline \delta u}
     =\frac{\delta h_n}{\delta u}
   - \partial_s \frac{\delta h_n}{\delta (\partial_s u)}
   + \partial_s^2 \frac{\delta h_n}{\delta (\partial_s^2 u)}
   - \partial_s^3 \frac{\delta h_n}{\delta (\partial_s^3 u)}+\cdots .
$$
\endroster}
\endproclaim

\tvskip
The existence of such $h_n$ will be guaranteed in Proposition 4-18.
Noting $\underline\Omega =\Omega\partial_s$,
and from the definition we have a recursion relation,
$$
	\partial_s \frac{\overline\delta h_n}{\overline\delta u}=
	\underline\Omega
          \frac{\overline\delta h_{n-1}}{\overline\delta u},
$$
if $h_n$ exists.
In Proposition 4-18, we show existence of a set of functionals
$\overline h_n = \res \dfrac{2^{2n}}{2n-1} L^{(2n-1)/2}$,
$u$ satisfying
$
\partial_s \overline h_n = \Omega \partial_s \overline h_{n-1}
$
with $\overline h_1=u/2$. Here $L = \partial_s^2 + u$ and
{\lq\lq}res" means the coefficient of the $\partial_s^{-1}$
in the notations in \S 4.
In other words,
$
\partial_s \overline h_n = \Omega^n \partial_s \overline h_{1}
 =\Omega^n \partial_s u/2
$.
Further from the definition, we have [D]
$$
	\int \frac{\overline\delta \overline h_n}
              {\overline \delta u} d s \equiv
	2(2n-3)\int \overline h_{n-1} ds,
$$
since
$\dint \dfrac{\overline\delta\overline h_n}{\overline \delta u}
 \equiv
\dint \dfrac{\delta\overline h_n}{\delta u}$ due to periodicity and
$$
\split
	\int \frac{\overline \delta }
             {\overline \delta u} \res (L^{r/2}) ds
       &= \int
    \res \left(   \sum_{i=1}^{r-1}(L^{1/2})^i
	\frac{\overline \delta L^{1/2}}{\overline \delta u}
             ( L^{1/2})^{r-i-1}  \right) ds \\
	&= r\int
  \res(   L^{(r-1)/2}
  \frac{\overline \delta L^{1/2}}{\overline \delta u} ) ds \\
	&= \frac{r}{2}\int
    \res(   L^{(r-2)/2}
   \frac{\overline \delta L}{\overline \delta u} ) ds \\
	&= \frac{r}{2}\int
    \res(   L^{(r-2)/2} ) ds.
\endsplit
$$
Let $h_n \equiv 2^n \overline h_{n-1} /(2n+1)$
 modulo periodic functions
and $X_n = \partial_s \overline h_{n}$ with $\overline h_0 = 0$.
Hence Definition 3-13 is guaranteed by Proposition 4-18.

Here we  give the vector fields $X_n$ and quantities
$h_n$ explicitly:
\proclaim{3-14. Example} {\rm( KdVH flow)}
$$     \matrix
  n= 0:& \quad X_0(u)=0, & h_0 = \frac{1}{2}u \\
   n= 1:& \quad X_1(u)=\partial_s u, & h_1 = \frac{1}{2}u^2 \\
   n= 2:& \quad X_2(u)=\partial_s
               (3 u^2 +\partial_s^2 u) ,
               &h_2 = u^3+\frac{1}{2}(\partial_s u)^2 \\
   n= 3:& \quad X_3(u)= \partial_s(10 u^3 + 5 (\partial_s u)^2
         +10 u \partial_s^2 u +\partial_s^4 u),
     &  h_3 = \frac{5}{2} u^4 + 10 u (\partial_s u)^3
     + (\partial_s ^2 u)^2.\\       \endmatrix
$$ $$     \align
   n= 1:& \quad \partial_{t_1} u + \partial_s u =0 ,    \\
   n= 2:& \quad \partial_{t_2}u
               +6 u\partial_s u +\partial_s^3 u =0 , \\
   n= 3:& \quad
    \partial_{t_3}u + 30 u^2 \partial_s u+20\partial_s u \partial_s^2 u
      +10u \partial_s^3 u + \partial_s^5 u=0 .\\
 \endalign$$
\endproclaim

\tvskip
\proclaim{3-15. Proposition}{\rm [AM]}
\roster

\item $\omega$ is a cocycle 2-form.

\item The KdVH flow has the Hamiltonian structures with
their Hamiltonian,
$$   H_n := \int_{S^1}  h_n \ d s, \quad (n \ge 0),
$$
with involutive relations for the Poisson bracket,
$\{H_n,H_m\}:=\omega(X_n,X_m)$,
$$
   \{H_n, H_m\} =0, \quad \text{for all } n,m.
$$

\item The $n$-th KdV flow has infinite conserved quantities $H_m$ $n \in
\Bbb Z_{\ge 0}$.

\item We have the relation,
$$   [\partial_{t_n}, \partial_{t_m}]u = 0, \quad \text{for all } n,m.
$$

\item For an arbitrary curve $\gamma$,
the $n$-th $(n \ge 1)$ KdV flow is uniquely determined. \endroster}
\endproclaim

\tvskip
\demo {Proof} We will prove these following to the arguments in [AM].
First we will show that $i_X \omega$ is exact: For all $n>0$,
we have the relation,
$$   i_{X_n} \omega( v) =\omega(X_n(u),v)
   = \int_{S^1} d s \frac{\overline\delta h_n}{\overline\delta u} v
   = (d H_n)(v) , \quad \text{for } n \ge 1.
$$
Hence $X_n(u)$ is a Hamiltonian vector field from the Definition
3-12 (2).
Our system is a Hamiltonian system and the  $n$-th KdV equation is given
by,
$$   u_{t_n} = X_n(u) .
$$
Next we will show that the KdVH flow is involutive.
As the time $t_m$ development of $H_n$ is given by
$$
	\partial_{t_m} H_n = \int
            \frac{\overline\delta h_n}{\overline\delta u}
                          \partial_{t_m} u
                          = \int \{H_n, H_m \},
$$
the involution relations are important.
From the Definition
3-12, we have  relations for $n \ge 1$,
$$ \split   X_n&=\partial_s
\frac{\overline\delta h_n} {\overline\delta u}\\
   &=\underline \Omega
   \frac{\overline\delta h_{n-1}} {\overline\delta u}.
   \endsplit
$$
Since in terms of $\omega$ in Definition 3-13 (1), the Poisson bracket
between $H_n$'s are given by $\{H_n,H_m\}=\omega(X_n,X_m)$, we obtain
the following relation for $n,m>0$:
$$\split  \{H_n,H_m\}&=\int_{S^1} ds
       \frac{\overline\delta h_n}{\overline\delta u} X_m(u)\\
   &=\int_{S^1} ds
   \frac{\overline\delta h_n}{\overline\delta u}\underline\Omega
            \frac{\overline\delta h_{m-1}}{\overline\delta u} \\
          &=\int_{S^1} ds
   \underline\Omega\frac{\overline\delta h_{n-1}}{\overline\delta u}
   \frac{\overline\delta h_{m}}{\overline\delta u}\\
 &=\{H_{n+1},H_{m-1}\} .
\endsplit $$
Using this relations and noting  $\{H_n,H_m\}=-\{H_m,H_n\}$, we will
prove the involutive relation. When both $n$ and $m$ are even or both
$n$ and $m$ are odd,
$$   \{H_{n},H_{m}\}=\{H_{(n+m)/2},H_{(n+m)/2}\}=0 .
$$
On the other hand, when $n$ is odd and $m$ is even,
$$   \{H_{n},H_{m}\}=\{H_{(n+m-1)/2},H_{(n+m-1)/2+1}\}=
   \{H_{(n+m-1)/2+1},H_{(n+m-1)/2}\}=0 .
$$
Hence $H_n$'s are involutive and the KdVH flow has infinite conserved
quantities.

We can express  the relation $\{H_{n},H_{m}\}=0$ by using a vector
representation for $n,m >0$,
$$   [X_n,X_m]=0.
$$
In the solution of the KdV hierarchy, we can identify $\partial_{t_n}$ with
$X_n$ itself: $\partial_{t_n}\equiv X_n$. Hence we obtain (4).

Further (5) can be proved as follows. For a given curve $\gamma$, we
uniquely have the data, $u$,$\partial_s u$,
$\partial_s^2 u$, $\cdots$. The KdV equations are given by
$$   \partial_t u = f(u,\partial_s u,
   \partial_s^2 u,\cdots) .
$$
Hence for an arbitrary curve $\gamma \in \CMeP$, the KdVH
flow is uniquely determined by the KdV hierarchy. Due to the
integrability, the {\lq\lq}time" development of the $\gamma$ is stably
determined. \qed \enddemo

\tvskip
Since the KdVH flow is a Hamiltonian system with infinite time parameters,
we can find a group $g \in G$
such that $\gamma_{t+t'}= g_{t'}\gamma_t$.
The multiplication is given as
$ g_{t'} g_{t}=g_{t'+t}$. $g_{0} $ is unit  and $g_{-t}$
is the inverse of $g_t$.
Further Proposition 3-15 (4) means that
 $[\partial_{t_1},\partial_{t_n}]u =0$
and the projection of $ \pi_{\elas}^{\PP}:\BMeP
\to \CMeP$ consists with the KdV flow.

Further as solving the KdV hierarchy is an initial problem with the first
derivative with respect to the time, for an arbitrary $\gamma\in\BMeP$
we can find the KdVH flow to which $\gamma$ belongs as an initial
state.

We will give a proposition as a summary of the  above arguments.

\proclaim{3-16. Proposition}
{\it \roster \item There is an Abelian group
$G:=\{ \exp( \sum_n  t_{n} \partial_{t_n} )\ $ $
 |\  t_n \in \Cal V^\infty\}$
 acting on the moduli spaces $\BMeP$ and $\CMeP$,
whose orbits are identified with the KdVH flow.

\item There is a fixed normal  subgroup $G_0$ of $G$,
$G_0=\{ g_{t_1}\ | \ t_1 \in \Bbb R \} \approx \roman{U}(1)$; $G_0$
 trivially acts upon $\CMeP$:
$\gamma =g_{t_1} \gamma$ for $\gamma \in \CMeP$ and $g_{t_1} \in G_0$.

\item The group $G/G_0$ acts on $\CMeP$.
\endroster }
\endproclaim

Hence Proposition 3-11 (2) is proved.
We can express the equivalent class in $\CMeP$ by the group action
in the following proposition.

\tvskip

\proclaim{3-17. Proposition}
{\it
\roster
\item Fixing $\gamma \in \CMeP$, $G/G_0$ whose
element is given as $g_{t_2,t_3,\cdots}$ transitively acts upon
$\frak C[\gamma]$: For an arbitrary $\gamma' \in \frak
C[\gamma]$, we can find an element $g_{t_2,t_3,\cdots}$ of
the group $G/G_0$
such that $\gamma = g_{t_2,t_3,\cdots} \gamma'$.

\item For  an arbitrary  $\gamma \in \CMeP$, there exists
the KdVH
flow: $\CMeP$ can be decomposed,
$$   \CMeP = \coprod \frak C[\gamma].
$$

\item For $\gamma\in \BMeC$,
the energy functional $\Cal E[\gamma]$ is exactly conserved for
the KdVH flow.
\endroster }
\endproclaim

\demo{Proof}
(1) and (2) are obvious from the properties of group.
(3) is proved because the energy $\Cal E[\gamma]$
of the loop $\gamma$ given by Definition 2-18 is
identified with the conserved quantity of $H_0$.
\qed \enddemo

\tvskip
Hence Proposition 3-11 (7) is proved from Proposition 3-17 (3).
By Propositions 3-11, 3-15, 3-16 and 3-17, we completely proved
our main theorem 3-2.

\tvskip
As we have the classification of $\CMeP$, we will use it
and go on to investigate the moduli space $\CMeP$ in rest of
this paper because our purpose is to get some knowledge of
the moduli space $\CMeP$.

\tvskip
For later convenience, we will introduce a quotient space.
Due to Theorem 3-2 and Proposition 3-17, $\CMeP$ has natural
projections induced by the equivalent relation $\simKdVH$, {\it i.e.},
 $\pi_{\roman{KdVHf}}:
 \frak C[\gamma] \mapsto (\gamma)$, where $(\gamma)$ is a representative
element of $\frak C[\gamma]$.

\proclaim{3-18. Definition}
{\rm
\roster \item  We define a quotient space of the moduli space by,
$        \fMeP:= \pi_{\roman{KdVHf}} \CMeP
   :=\CMeP/\simKdVH.
$

\item The natural projection is denoted by
$\pi_\elas:\BMeP$
$\to \fMeP$.
\endroster }
\endproclaim

\tvskip
\proclaim{3-19. Remark}
{\rm
We will comment on Proposition 3-15 (4),
$$
	[\partial_s, \partial_{t_n}] = 0 \quad
         \text{ for } n >0.
$$
As the KdVH flow is very regular, we can regard
$\frak C[\gamma]\times S^1$ $\in \BMeP$ as a manifold.
Accordingly
 $\partial_{t_n}$ are regarded as a vector field. We will use it
 as a generator of a cohomology in \S 7.

\roster

\item It means the local length $d s$ preserves for the KdVH flow.

\item It can be
interpreted as Frobenius integrability conditions.

\item It is known as the compatibility condition or zero curvature
conditions known in the soliton physics.

}
\endproclaim

%\newpage
\tvskip
\centerline{\twobf \S 4. Algebro-Geometric Properties of the KdV flow I}
\centerline{\twobf --Algebraic Properties--}
\tvskip

As we proved Theorem 3-4, we will use the relation between the moduli
of a quantized elastica $\BMeP$ and the KdV flow in order to
give a finer classification, which is
based on the study {\it finite type flow} in $\BMeP$ and $\CMeP$,
in this section.
However as this classification comes from the algebraic investigation
of the KdV flow, we should replace
the base function space in the category of the smooth functions
by that of the formal power series
 in order to explain this classification,
though we need some subtle treatments.

This section is devoted for investigations of a commutative
differential ring, which were given by Mulase in [Mul],
  Burchnall and Chaundy in about seventy years ago [BC1, 2, Ba3],
and Mumford in [Mum1].
Our argument  basically
follows the arguments of Mulase for the
Schottky problem [Mul] and of Sato [SN, SS].
Following their theories, we will consider a part of
the moduli of a quantized elastica using
 the formal power series. Since
the part is dense in $\BMeP$ as mentioned in Theorem 4-2,
the replacement of the base field is not so critical.

Although investigation of $\gamma$ as a real one-dimensional curve
is our main subject,
we deal with a hyperelliptic curve as a complex one-dimensional curve
in the context of algebraic geometry in this section and next
section. Thus readers
should not confuse the terms {\lq\lq}curve" in the
categories of the differential geometry and the algebraic geometry.
We basically refer the complex algebraic curve {\it algebraic
curve}, {\it hyperelliptic curve} or {\it elliptic curve}
whereas we call such a real curve just {\it curve}.

\vskip 1.0 cm

Let us start this section with the following lemma.

\proclaim{4-1. Lemma}
{\it
If there is  a natural number $N$ such that
$\partial_{t_N}u$ is an eigen vector of the operator
$\Omega$ with an eigenvalue $k \in \CC$, {\it i.e.},
$$   k\partial_{t_N} u=\Omega \partial_{t_N}u ,
$$
 $\partial_{t_m}$ is a scalar multiplication of $\partial_{t_N}$
for all $m \ge N$.
Further by introducing $t_n'$ $n>N$ and setting
$\partial_{t_n'} := \partial_{t_n}-k^{n-N}\partial_{t_N}$,
the relation becomes $\partial_{t_n'} u \equiv 0$. }
\endproclaim

\demo {Proof} This proof is easily from  Definition 3-2
and Proposition 3-11 (5). \qed \enddemo

Lemma 4-1 means that some orbits in infinite dimensional vector space
$\Cal V^\infty$ are essentially  reduced to an orbit
consisting of finite $N$ dimensional vector space.
Let us refer this flow {\it finite flow} or {\it finite $N$-type flow}.

Here we will give our second main theorem:

\tvskip\proclaim{4-2. Theorem}
{\it \roster
\item We will write the set of the finite type flow by
$\BMeP_{,\finite}$
and the set of finite $g$-type flow by
$\BMeP_g$.
The moduli space of the elastica has decomposition,
$$
\BMeP_{, \finite}:=
\coprod_{g < \infty} \BMeP_g.
$$

\item $\BMeP_{,\finite}$ is dense in
$\BMeP$ with respect to the canonical topology
determined by the KdVH flow.

\endroster}
\endproclaim

\vskip 1.0 cm

In order to prove Theorem 4-2, we prepare
the knowledge of the KdV equation.  A more concrete statement
appears in Proposition 4-33.
Before that, we will recall the result of the Whitney for
the quasi-analytic system [Wh], which is easily proved by
the Weierstrass preparation theorem.

\proclaim{4-3. Proposition}{\rm [Wh]}
For a presheaf $\SCinf(\RR)$ of smooth functions
 over $\RR$ and a presheaf  $\SF(\RR)$ of
$\CC$-valued formal power series over $\RR$,
 we have a surjective presheaf morphism,
$$
	\eta: \SCinf(\RR) \to \SF(\RR),
$$
{\it i.e.},
for a germ $f \in \Gamma_p(\SCinf(\RR))$ and $t_1 \in \RR$ around $t_1^0$,
$$
   \eta(f) = \sum_i^\infty \left[\frac{d^i}{d s^i} f\right]_{t_1=t_1^0}
          (t_1-t_1^0)^i.
$$

\endproclaim

\tvskip

The map $\eta$ is not injective,
{\it e.g.}, due to a function $f$,
$f(s) = 0$ at $s=s_0$ and $f(s) = \exp( -1/(s-s_0)^2)$ at
 otherwise points.

Since $\eta$ is a local correspondence,
the map $\eta$ can be applied for the presheaves of
the $\Cinf$ functions and the formal power series over $S^1$,
or $\eta: \SCinf(S^1) \to  \SF(S^1)$.

On the other hand, for  an arbitrary element in $\SCinf(S^1)$
we can find a sequence in the presheaf of the formal power series
$\SF(S^1)$ which converges to the element
 using the same Weierstrass  preparation theorem.

By using these properties, we will replace the base ring
$\Cinf(S^1,\RR)$ with the formal power series in this section.
\tvskip

In order to express the system of the KdV equation,
we will mention the differential algebra and its division
algebra over a commutative ring $R$.
As we show in Definition 5-15 and
Proposition 5-18, the hyperelliptic $\wp$ function obeys
 the KdV
equation and  has a singularity of the second order.
Hence we might be ought to deal with Lourant expansion
ring $\CC[[t_1]][t_1^{-1}]$ as $R$.
However as we are concerned with one of the KdVH flows which
are  finite and real valued,
we deal only with a finite {\lq\lq}real" valued part of $\wp$
and avoid the singular points.
In other words, we employ a formal series ring $\CC[[t_1]]$
as the ring $R$.

In this section, $t_1$ is dealt with as a generic parameter but
can be regarded as a real (complex) number $t_1 \in \RR$ ($\CC)$. After
considering periodicity, we regard it as a point of $S^1 = \RR/\ZZ$
in later.
For a convenience, let $\partial_1 := \partial/\partial t_1$.

Here we assume that all algebra and subalgebra have unit as their
definitions in this article.

\proclaim{4-4. Definition}{\rm [Mul, S, SN, SS]}
{\rm
\roster
\item The differential ring $\DDf$  is defined by
$$  \DDf:=\{ \sum_{k\ge 0}^N a_k\partial_1^k \ |
    \ N < \infty, \  a_k\in \CC[[t_1]] \}.
$$

\item Let us identify commutative subalgebras
$B_1$ and $B_2$ in $\DDf$ if there exists
an invertible element $r\in\CC[[t_1]]^\times$ such that
$$
	B_1 = r B_2 r^{-1}.
$$
We define {\it standard representation } $B_s$ in
the equivalent class  $[B_1]$ by tuning $r \in \CC[[t_1]]^\times$
such that it contains,
$$
	 \partial_1^n + b_{n-2}\partial_1^{n-2} +\cdots + b_0
             \in B_s.
$$

\item The degree of a differential operator
and the projections $+$ and $-$ are defined by the same
as the case $\EEs$ in Definition 3-1.

\item The micro-differential ring   $\EEf$ to  $\DDf$ is defined by,
$$
    \EEf:=\{ \sum_{k= -\infty}^N a_k\partial_1^k \ |
   \ N < \infty, \  a_k \in \CC[[t_1]] \}.
$$

\item The constant coefficient subring of $\EEf$ is defined by,
$$
     \EEc:=\{ \sum_{k= -\infty}^N a_k\partial_1^k \ |
   \ N < \infty, \  a_k \in \CC \}.
$$

\item An invertible set $\WWf$ in $\DDf$ is defined by,
$$
	\WWf := \{ W \in \EEf \ | \ W = 1 + \sum_{i=1}^\infty
              a_i \partial_1^{-i}, a_k \in \CC[[t_1]] \},
         \quad \WWc := \WWf \cap \EEc.
$$

\endroster}
\endproclaim

For the readers who are not familiar with valuation,
we will review it.

\proclaim{4-5. Definition }{\rm [Ha]
$\Cal K$ is a topological field.
Let us call a topology space $\Cal E$ {\it left linear topological space}
if $\Cal E$ satisfies
\roster

\item $\Cal E$ is a $\Cal K$ linear space.

\item A map from $\Cal K \times \Cal E$ to $\Cal E$
($(\lambda, x) \mapsto \lambda x$) and
addition $x+y \in \Cal E$ are continues.

}
\endproclaim

\proclaim{4-6. Definition 4-6 }{\rm [Ha]
\roster

\item Let $\Cal K$ be a field.
A valuation of $\Cal K$ with values in $\ZZ$
is a map  $\val:\Cal K \to \ZZ$,
for all $x,y\in \Cal K$, $x,y\neq0$,
$$
    \val(xy) = \val(x)+\val(y), \quad
   \val(x+y)\ge \min(\val(x), \val(y)).
$$
and $\val(0)=\infty$.

\item The set $\Cal R:=\{x\in \Cal K\ | \ \val(x)\ge 0\ \}$ is
a local subring, called {\it valuation ring}.

\item The set $\frak m:=\{x\in \Cal K\ | \ \val(x)> 0\ \}$ is
called a {\it local ideal} of $\Cal R$

\item
Let the metric of $\Cal K$ be
$|x| :=\ee^{-\val(x)}$ for $x \in \Cal K$,
which is called {\it non-Archimedian metric}.

}

\endproclaim

For example, the valuation of a commutative ring $\CC[x]$
is given by its degree, {\it i.e.}, for $f(x) \in \CC[x]$,
$\val(x)=\deg_x(x)$. For a more general commutative ring,
we can find a local parameter by localization at a prime
ideal and its valuation is given by its degree of the local
parameter.

The valuation ring is a linear topological space due to
the non-Archimedian metric [Ha].
Similarly, we have the following proposition [SN],
which is naturally obtained.

\vskip 0.5 cm \proclaim{4-7. Proposition}{\rm [SN]

{\it\roster
\item  When we define $\EEf_m:=
\{ D \in \EEf \ | \ \roman{deg} D \le m \ \}$,
$\EEf$ has filter,
$$
   \EEf= \cup_{m} \EEf_m, \quad\{0\}= \cap_{m} \EEf_m, \quad
   \EEf_m \subset \EEf_{m+1} .
$$

\item $\EEf$ is a linear topological space with respect to this
filter.

\item $\EEf$ is an infinite dimensional algebra given by the
formal power sires whose element converges in the filter topology.

\item In $\EEf$, we can define  valuations in $\CC[[t]]$ and
 $\EEf$ as
$P = \sum_{i=-\infty}^\infty a_i \partial_1^i\in \EEf$
$$
\val(a) := \max\{ m\in \NN | \partial_1^m a \neq 0 \},
\quad
	\val(P) : = \inf\{ \val(a_i)-i \}.
$$

\endroster}
\endproclaim

\tvskip
Formally Proposition 4-7 is obvious from their definitions but we need
rigorous arguments to justify them mathematically, which is written in
[SW, SS, S2].
The differential operators appearing in the soliton theory and
in the following arguments converge in this topology.

\proclaim{4-8. Lemma}{\rm [Mul, SN, SS]}
\roster

\item
The adjoint map for $W\in \WWf$, $Ad(W): \EEf \to \EEf$
$(Ad(W)P = W P W^{-1})$,
defines the automorphism in $\EEf$.
$Ad(W)|_{\EEf_m}$ is invariant,
{\it i.e.}, $Ad(W)\EEf_m=\EEf_m$

\item For an operator $\tilde L \in \Cal L$, where
$$
\Cal L := \{ D \in \EEf \ | \
             D = \partial_s + \sum_{i=1}^\infty a_i \partial_s^{-i}\ \},
$$
we can find a unique $W \in \WWf$ modulo $\WWc$ such that
$$
       Ad(W) \tilde L = \partial_s,
$$
and then this relation induces the isomorphisms of
$$
\WWf/\WWc \approx \Cal L.
$$

\item
For every standard commutative subalgebra $\aaf \subset \DDf$,
there is a $\CC$-subalgebra $\aac \in \EEc$ such that
is $\CC$-isomorphic to $\aaf$ and
$$
	\aac \cap \EEc_- =\{0\}.
$$

\endroster

\endproclaim

\demo{Proof}
(1) is trivial.
(2) Let us find $W=\sum_{i=0}^\infty w_i \partial_1^{-i}$
such that $\tilde L=W \partial_1 W^{-1}$. Noting
$\tilde L = \partial_1 + \tilde L_-$, the relation is
reduced to $[\partial_1, W] = -\tilde  L_- W$,
{\it i.e.},
$$
  \partial_1 w_{k-1} = - \sum_{i+j+r=k, i\ge2}
 \pmatrix 1-i \\ r \endpmatrix u_i \partial_1^r w_j.
$$
 Then we can
recurrently determine $w_i$ from small numbers
since $\CC[[t]]$ has indefinite integrals.
When we find such $W_1$ and  $W_2$, {\it i.e.},
$\tilde L W_a-W_a \partial_1 = 0$, then
$W_1 \partial_1 W_1^{-1} W_2 - W_1 W_1^{-1}W_2 \partial_1=0$
or $[\partial_1, W_1^{-1} W_2] =0$. Hence
$ W_1^{-1} W_2\in \WWc$.
(3): Let us take a monic element of $\aaf$ such that its
form is
$$
L_n =  \partial_1^n + b_{n-2}\partial_1^{n-2} +\cdots + b_0.
$$
Then we have $\tilde L = L_n \in \Cal L$.
Let $S\in \WWf$ such that $\tilde L = S\partial_1 S^{-1}$.
Then $\aac := S^{-1}\aaf S$. For  an arbitrary
 $P\in \aaf$, $ S^{-1} P S$  belongs to
$\EEc$, {\it i.e.}, $[\partial_1, S^{-1} P S]=0$ because
 $[\partial_1, S^{-1} P S]=S^{-1}[S \partial_1 S^{-1}, P] S
=S^{-1}[L, P] S =0$ due to  the assumption $[P,L]=0$.
 Further the inner automorphism preserves
the order of the operator.\qed
\enddemo

\tvskip
As we will not prove here, it is known that
if we define a  left $\EEf$-module,
$$
          \VVf := \EEf/  \EEf t_1,
$$
the homomorphism from $\EEf$ to the endomorphism of  $\VVf$ is injective.
In other words,  the endomorphism is faithful
if it can be regarded as a representation.
$\VVf$ has a valuational topology $\val$ and becomes
a graded module.
 There the valuation and graded topology are
identified.
Further we have a natural $\EEc$-module isomorphism [SN],
$$
        \VVf \approx \EEc.
$$
Further we consider
an embedding of a submodule $\VVf_0$ into $\VVf$ with zero-index map
for a certain index, which can be regarded as a
Grassmannian manifold in a certain sense.
We note that the above isomorphism is not meaning of
$\EEf$-module.

In this article, we characterize such an embedding
by a finite subset of natural numbers $F$, which can be
regarded as the Weierstrass gap in the infinite
point of a corresponding algebraic curve.

Further we should note that the adjoint map $Ad$ is the
key of the Sato theory and in this section, we sometimes
call it {\it gauge} transformation.

\tvskip

\proclaim{4-9. Definition}{\rm [Mul]
A $\CC$-subalgebra $\aac$ in $\EEc$
is called a {\it rank one subalgebra} if it has
$\CC$-linear basis whose indices corresponds
to all of integer except a finite subset $F$,
{\it i.e.},
$$
\NN - F =\NN_{\aac}:=\{ n \in \NN \ | \
\exists P \in \aac \text{ such that } \ord(P)=n\ \}
$$
and
$$
	\aac \cap \EEc_- =\{0\}.
$$
}
\endproclaim

As $\aac$ is a $\CC$-algebra, there is a monic element
$P_n$ in $\aac$ of
order $n\in \NN-F$ with $P_0:=1$. Then $\{P_n|n\in \NN-F\}$
forms a $\CC$-linear basis of $\aac$. In other words arbitrary
$P\in \aac$ can be represented by $\CC$-linear combinations of
monic $P_n$ elements. In fact if the order of $P$ is $m$,
there exists $c\in \CC$ such that the order of
 $P-c P_m \in \aac$ must be less than $m$.
Such a recursion process
gives us the representations.

\proclaim{4-10. Lemma}{\rm [Mul]}
\roster

 Let $\aac\neq \CC$ be a rank one subalgebra, and $P$ and $Q$
be elements in $\aac$ whose orders are coprime.

\item $\dim_\CC(\aac/\CC[P,Q])< +\infty.$

\item $\aac$ is finite $\CC[P]$-module. There is a nontrivial
polynomial $f(x,y)\in \CC[x,y]$ such that $f(P,Q)=0$.

\item The transcendence degree of
$\aac$ over $\CC$ is one.

\item By regarding $\aac$ as a graded module with respect to
degree of differential operators:
$$
	\aac^{(n)}:=\{ P\in \aac\ | \ \ord(P)\le n\ \},
$$
$$
	\aac_n:= \aac^{(n)} \oplus \aac^{(n-1)}\cdot I \oplus
             \aac^{(n-2)}\cdot I^2
 \oplus\cdots \oplus \aac^{(0)} \cdot I^n,
$$
$$
  \gr \aac=\oplus_{n=0}^\infty \aac_n,
$$
we regard ${\roman{Proj}} (\gr \aac) $ as an algebraic curve
$C$.
Here $I$ is the identity of $\aac$.

\item Let $\HH^1(\aac)= \EEc/\aac\oplus \EEc_-$. We have
$$
	\HH^1(C, \Cal O_{C}^\times) = \HH^1(\aac),
$$
where $\Cal O$ is the sheaf of holomorphic functions on $C$ of
(4) and $\Cal O^\times$ is a multiplicative subset of $\Cal O$.

\endroster
\endproclaim

\demo{Proof}
Let  $\GCD(m,n)$ denote the greatest common divisor of two
non-negative integer $m$ and $n$. Since the rank of $\aac$ is unit,
we have the relations
$$
1 = \min\{ \GCD(\ord(P'),\ord(Q')) \ | \ P', Q' \in \aac \ \}.
$$
and the orders of $P$ and $Q$ are coprime. Hence $\CC[P,Q] \subset \aac$.
As $N_{\CC[P,Q]}=\NN$ and $\CC[P,Q]$ is $\CC$-linear vector space,
$N_\aac-N_{\CC[P,Q]}$ must be finite set.
Hence $\dim_\CC(\aac/\CC[P,Q])$ must finite.
On the other hand, since
$\NN-\{\ord\{P^m,Q^n\} \ | \ m, n \in \ZZ_{\ge 0}\}$ must be finite
set, $P$ and $Q$ satisfy an algebraic relation $f(P,Q)=0$.
Further the proofs of (4) and (5) are due to theory of
an ordinary commutative ring [Ha].
\qed
\enddemo

We note that $F$ in Definition 4-9 is related to
the Weierstrass gap at infinity point of the algebraic curve $C$.
\tvskip

After this point,
we will concentrate our attention only on the operator
$L=\partial_1^2 +u $, which is related to the  KdV equation:
$$
\Cal L_2 := \{ D \in \EEf \ | \ D = \partial_1^2 + u, \quad u \in
 \CC[[t_1]] \}.
$$
\tvskip
We give its related operators as examples.
\proclaim{4-11. Example}
$$\split    L^{1/2}&=\partial_1 + \frac{1}{2} u \partial_1^{-1}
                 -\frac{1}{4} (\partial_1 u)\partial_1^{-2}
               +\frac{1}{8} ((\partial_1^2 u)-u^2)\partial_1^{-3}\\
     &+\frac{1}{16} (6 u (\partial_1 u)-\partial_1^3 u)\partial_1^{-4}
       -\frac{1}{32}( -2 u^3 + 14 u (\partial_1^2 u)
   + 11 (\partial_1 u)^2 - (\partial_1^3 u))\partial_1^{-5} + \cdots,\\
   4 L^{3/2} &= 4 \partial_1^3 + 3\partial_1 u + 3 u \partial_1
   +\left( \frac{1}{2} \partial_1^2 u
   +\frac{3}{2} u^2\right) \partial_1^{-1} +\cdots.,\\
   16 L^{5/2} &= 16 \partial_1^5 +40u \partial_1^3
   + 60 (\partial_1 u) \partial_1^2
   +50 (\partial_1^2 u) \partial_1
   + 30 u^2\partial_1 15 (\partial_1^3 u) +30 u (\partial_1 u) \\
   &+\left(5 \left(u^3 + \frac{1}{2}(\partial_1 u)^2 \right)
         + \partial_1 f(u, \partial_1 u, \cdots)\right) \partial_1^{-1}
          +\cdots.\\
\endsplit$$
\endproclaim
Here $\partial_1 f(u, \partial_1 u, \cdots)$ is a functional of $u,
\partial_1 u, \cdots$.

\tvskip

Let us fix the operator $P = \partial_1^2$
of $\aac$ in Lemma 4-10 because we only consider
 $L=W\partial_1^2 W^{-1}$.
From the primitive number theory,
for an odd number $m$ and  an integer $n(>m)$,
we find $ a,b \in \ZZ_{\ge 0}$ such that
$$
	n = a m + 2 b , \quad ( a,b \in \ZZ_{\ge 0}).
         \tag 4.1
$$
When we fixed $\aac$ as a rank one subalgebra, the partner $Q$ of
$P\equiv\partial_1^2$
 in the Lemma 4-10 is an operator whose order is given by an
odd number $2g+1$. Thus $F$ in Definition 4-9 is given by a
smaller sequence of odd numbers, $\{1,3,5,7,9,\cdots,2g-1\}$.
Let us introduce a set of such subrings $\aac$ in $\EEc$.

\proclaim{4-12. Definition}{\rm
$$
   \split
      \AAc:=\{ \aac  \
              |& \ \text{ $\aac$ is a rank one subalgebra, }\\
              & \ \exists W \in \WWf\ \text{ such that
              $W\aac W^{-1}\in \DDf$ is a commutative subalgebra,}\\
              &\NN- \NN_{\aac}\subset
          \{1,3,\cdots,2g-1\} ,\ g < \infty\ \}.\\
\endsplit
$$}
\endproclaim

For the case of $g=1$,  $\aac\equiv W^{-1}\CC[L,[L^{1/2}]_+]W$.
Since $[L^{1/2}]_+\equiv \partial_1$, $[L, \partial_1]=0$
and thus $u$ must be $\CC$. In other words, $Q$ must be $\partial_1$,
$\CC[\partial_1^2,\partial_1]\equiv \CC[\partial_1]$.
For the case $g=1$, it becomes an ordinary polynomial ring.

\tvskip
\proclaim{4-13}
{\rm
We recall that an algebraic curve with a morphism to $\PP$
of order two is called hyperelliptic curve.
A hyperelliptic curve $C_g$  of genus $g$ $(g\ge 1)$,
including the case of elliptic curve,
is given by the homogeneous equation,
$$
   Y^2 Z^{2g-1} = h_g(X,W) := \lambda_0Z^{2g+1} +\lambda_1 X Z^{2g}
        +\lambda_2 X^2 Z^{2g-1}+\cdots +\lambda_{2g+1} X^{2g+1},
$$
where $\lambda_{2g+1}\equiv1$ and $\lambda_j$'s are
complex values. }
\endproclaim

\proclaim{4-14. Lemma}
Let $L =\partial_1^2+u\equiv W\partial_1^2 W^{-1}\in \Cal L_2$
to $W \in \WWf$.
\roster

\item $L^{n/2}=W\partial_1^n W^{-1}$.

\item $
 [{L^{2n}}_+,L]\equiv [{L^{2n}},L]\equiv0 $.

\item The set of the differential operators in $\DDf$
which commute with a given operator $L_2\in \DDf$ is itself
a commutative subalgebra of $\DDf$.

\item $\aac\in \AAc$ is $\CC[\partial_1^2]$-module and by
considering
$$
	\AAc = \sum \aac,
$$
$\AAc$ is also $\CC[\partial_1^2]$-module.

\item
For  an arbitrary
 $\aac\in \AAc$, we can find $Q_g\in \aac$
which satisfies an affine equation,
$$
  Q_g^2 = h_g(\partial_1^2,1),
$$
so that there is a $W\in \WWf/\WWc$ such that
$W \partial_1^2 W^{-1}=L$ and
$W Q W^{-1}$ are commutative in $\DDf$.
Further we have found a
 hyperelliptic curve $C=\roman{Proj}(\gr \aac )$
and
$$
	\HH^1(C, \Cal O_C) = \HH^1(\aac),
$$
which are generated by $<\partial_1, \partial_1^3,
 \cdots, \partial_1^{2g-1}>$.

\endroster
\endproclaim

\demo{Proof}
(1) and (2) can be shown by direct computations.  On (3), we consider
a commutative differential ring in $\DDf$ such that
$\Cal B:=\{P\in \DDf\ | \ [P,L]=0\ \}$.
Since $L^{1/2}=W \partial_1 W^{-1}$,
$[\partial_1, W^{-1} P W]=W^{-1}[L^{1/2}, P] W =0$ because of
the assumption. Hence $ W^{-1} P W$ is an element of $\EEc$ and
thus we can find $\aac \in \AAc$ such that $\Cal B = W^{-1} \aac W$.
Hence $\Cal B$ is a commutative ring.
Next  (4) is trivial.
From the definition Lemma 4-10, and Eq. (4.1), we reach (5).  \qed
\enddemo

\tvskip

Next we will consider the filter structure in $\AAc$
and its completion with respect to the filtration.

\proclaim{4-15. Proposition}
Let us define a filter,
$$
	F_g\AAc:=\{ \aac \in \AAc \ | \ \NN- \NN_{\aac}\subset
          \{1,3,\cdots,2g-1\} \}.
$$
This satisfies the following relations:

\roster
\item
$$
	F_g\AAc \subset F_{g+1}\AAc, \quad F_{n}\AAc\equiv0,
        \quad n<0.
$$

\item By letting $\AAc_g:=F_g\AAc/ F_{g-1}\AAc$,
there is a large gauge transformation between
$\aac_1, \aac_2 \in\AAc_g$, {\it i.e.}, there exists
$W \in \WWf$ such that
 $\aac_1 = W \aac_2 W^{-1}$.

\item The direct limit of the filtration gives
$$
\split
 \overline{\AAc} &:= \lim_{\to} F_g\AAc\\
                 &=\{\aac \in \EEc \ | \
                  \exists W\in \WWf \text{ such that }\\
           &{}\quad W \aac W^{-1} \in \DDf
             \text{ is a  subalgebra, }
          \NN- \NN_\aac \subset 2\NN-1 \ \}.\\
 \endsplit
$$
\endroster
\endproclaim

\demo{Proof} (1) and (3) are obvious.
(2) is due to the proof of Proposition 4-16 (2).
\qed \enddemo

For each
element $\aac\in \AAc_g$, we consider the
 correspondence in Lemma 4-10 (4), {\it i.e.},
$\roman{Proj}(\gr\aac_g)$. It turns out that
$\AAc_g$ is isomorphic to the set of the
hyperelliptic curves with genus $g$.

\proclaim{4-16. Proposition}
\roster

\item The set $\AAf$ of commutative subrings in $\DDf$
inherits from the above filtration of $\AAc$.

\item
For any elements $L_1$ and $L_2$  in $\Cal L_2$,
there is a gauge transformation $W \in \WWt/\WWc$ such
that
$$
	L_1 = W L_2 W^{-1}.
$$
\endroster

\endproclaim

\demo{Proof} (1) is trivial. (2): There is
an element $W_a \in \WWt/\WWc$ such that
$L_a=W_a \partial_1^2 W_a^{-1}$ for $(a=1,2)$.
Hence $L_2 =W_2 W_1^{-1} L_1 W_1 W_2^{-1}$. \qed
\enddemo

\tvskip
As we described the tools and their properties for the
differential ring over $\CC[[t_1]]$, we will extend its
base field to $\CC[[t_1, t_2, \cdots ]]$. However
before we will give the extension in Definition 4-19, we digress
and show a connection between $\Omega$ in Definition 3-2
and $L$ in $\Cal L_2$.
Following the arguments in [D], we firstly prepare a lemma.

\proclaim{4-17. Lemma}{\rm [D]}
The {\lq\lq}resolvent" operator for $L=\partial_1^2 +u $,
$$
	T^{(\pm)} :=\left[ \frac{1}{2 z^2} \sum_{r = -\infty}^\infty
             (\pm z)^r L^{-r/2} \right]_-,
$$
 has the following properties:

\roster

\item $ (T^{(+)} + T^{(-)} ) = ( L - z^2)^{-1}$.

\item $[T^{(\pm)}(L-z^2) ]_- = [(L-z^2)T^{(\pm)}]_- = 0$.

\item When we define a map for a $X\in \EEs$, called Adler map [D],
$$
	\frak h(X):= [(L-z^2) X]_+ (L-z^2)
                  - (L-z^2) [X (L-z^2)]_+ ,
$$
we have the relation, $\frak h(T^{(\pm)})=0$.

\item $T^{(\pm)}$ has a formal expansion,
$$
	T^{(\pm)}= \sum_{r=1}^\infty S^{(\pm)}_r \partial_1^{-r}.
$$

\endroster

\endproclaim

\demo{Proof} (1) is trivial.
(2) is given by the relation,
$$
\split
[(L-z^2) T^{(\pm)}]_-&=\frac{1}{2}\left[(L-z^2) \left[
 \sum_{r = -\infty}^\infty
             (\pm z)^r L^{-r/2} \right]_-\right]_-\\
&=\frac{1}{2}\left[(L-z^2)
 \sum_{r = -\infty}^\infty
             (\pm z)^r L^{-r/2} \right]_-\\
&=\frac{1}{2}\left[
 \sum_{r = -\infty}^\infty
             \left((\pm z)^r L -z^2 (\pm z)^r
              \right) L^{-r/2} \right]_- .\\
\endsplit
$$
It is clear that it vanishes.
(3) is proved due to the property of the Adler map,
$$
	\frak h(X)\equiv - [(L-z^2) X]_- (L-z^2)
             + (L-z^2) [X (L-z^2)]_- .
$$
(4) is obvious from the definition of the resolvent.
\qed \enddemo

Due to the lemma, we gave the connection.

\proclaim{4-18. Proposition}{\rm [D]}
\roster
\item
$$  [2^{(n-1)}[{L^{(2n-1)/2}}]_+,L]=\Omega_1^n\partial_1 u,
$$
where
$\Omega_1:= \partial_1^2 + 2u + 2 \partial_1 u \partial_1^{-1}$.

\item By letting $\overline h_n
= \res \dfrac{2^{2n}}{2n+1} L^{(2n-1)/2}$,
$(${\lq\lq}$\res$" means the coefficient of $\partial_1^{-1})$,
we have
$$
	\partial_1 \overline h_n
= \Omega_1 \partial_1 \overline h_{n-1}.
$$

\endroster
\endproclaim

\demo{\Proof}
Due to the condition $\frak h(T^{(\pm)})=0$,
we can determine the first two coefficients $S_1$ and $S_2$ as
$$
	\partial_s^3 S^{(\pm)}_1 + 2 ( \partial_s u) S^{(\pm)}_1 +
                      4( u + z^2) \partial_s S^{(\pm)}_1 =0,
$$ $$
	\partial_s^2 S^{(\pm)}_2 =
- \frac{1}{2} \partial_s S^{(\pm)}_1.
$$
Let us consider the following operator,
$$
	(T^{(+)}- T^{(-)})
           =\left[ \sum_{r = -\infty}^\infty
             z^{2r+1} L^{-(2r+1)/2} \right]_-.
$$
The left hand side in the relation,
$$
    [ L^{(2r+1)/2}]_+ L - L[ L^{(2r+1)/2}]_+
     =L[ L^{(2r+1)/2}]_- -  [ L^{(2r+1)/2}]_- L,
$$
appears as a coefficient of
  $z^{2r -1 }$ in the series
$ \frak h(T^{(+)}- T^{(-)}) $ with respect to $z$.
Thus we are concerned with $S_r := S^{(+)}_r - S^{(-)}_r $,
which must have the expansion,
$$
           S_r = \sum_{i=-\infty}^\infty S_r^{(i)} z^{2i +1}.
$$
Comparing the coefficients in $z^{2r-1}$, we obtain,
$$
	4\partial_1 S_1^{(i+1)}=(\partial_s^3
          +   2  (\partial_s u) + 4 u  \partial_1) S^{(i)}_1
          =\Omega_1 \partial_1 S^{(i)}_1.
$$
We have the relation,
$$
     [[L^{(2r+1)/2}]_+, L]
      = \frac{1}{4}\partial_1 S_1^{2r+1}=
      \frac{1}{4^r}\Omega_1^r \partial_1 S^{(1)}_1,
$$
with $S_1^{1}=-u/2$.
Then we let $\overline h_n$  identified with $S_n$ by tuning
its coefficient. \qed
\tvskip

As we finished the digression, we extend
$\CC[[t_1]]$ to  $\CC[[t_1, t_2, \cdots]]$.
In the extension of the valuation over
 $\CC[[t_1]]$ to that of  $\CC[[t_1, t_2, \cdots]]$,
let the degree of $t_i^n$ be $(2i-1)n$.

\proclaim{4-19. Definition}{\rm [Mul, SN]
\roster
\item The differential ring $\DDt$, and its related set and ring are
defined by,
$$
\DDt:=\{ \sum_{k\ge 0}^N a_k\partial_1^k \ |
    \ N < \infty, \  a_k\in \CC[[t_1, t_2, \cdots]] \},
$$
$$
    \EEt:=\{ \sum_{k= -\infty}^N a_k\partial_1^k \ |
   \ N < \infty, \  a_k \in \CC[[t_1,t_2, \cdots]] \},\quad
     \EEt = \DDt + \EEt_-,
$$
$$
\Cal L_2^t := \{ D \in \EEt \ | \ D = \partial_1^2 + u, \quad u \in
 \CC[[t_1,t_2,t_3,\cdots]] \}.
$$

\item By letting $\val(t_i^n):=(2i-1)n$, we extend the valuation of
$\DDt$ and $\EEt$, which are also called valuations of $\DDt$ and $\EEt$.

\item
$$
	\WWt := \{ W \in \EEt \ | \ W = 1 + \sum_{i=1}^\infty w_i
            \partial_1^{-i} \
             \}.
$$

\item
$$
   \hat\DDt := \{\ P=\sum_{i=0}^\infty a_i \partial_1^i\in \DDt\ | \
            \exists N \in \NN, \ \val(a_i)>i - N\
              \text{ for } \forall i \gg 0\
                \},
$$
$$
	\hat\EEt := \{\ P=\sum_{i=-\infty}^\infty
               a_i \partial_1^i\in \EEt\ | \
            \exists N \in \NN, \ \val(a_i)>i - N\
                \text{ for } \forall i \gg 0\
                   \}.
$$

\endroster}
\endproclaim
\tvskip
We note that   $\DDt$, $\EEt$ and so on,  have
natural embeddings of $\DDf$, $\EEf$ and so on,
{\it e.g.},
$$
	\DDf\ni P(t_1) \mapsto P(t_1, 0, 0, \cdots ) \in \DDt.
$$
Using the embeddings, we regard $\DDf$ as a subring of
$\DDt$ as following.

\proclaim{4-20. Definition}
\roster

\item The moduli space of the KdV equations is defined by
$$   \BMKdV:=\{ u\in \CC[[t_1,t_2,\cdots]]
      \ |\  \partial_{t_n} u - \Omega_1^{n-1}\partial_1 u =0
\text{ for }\ \forall n \  \} ,
\quad   \CMKdV:= \BMKdV/(t_1) ,
$$
where
$\Omega_1:= \partial_1^2 + 2u + 2 \partial_1 u \partial_1^{-1}$.
Here $y$ is an element of the vector space generated by
$t_1, t_2, \cdots$.
\item
If $\tilde L \in $ $
\Cal L^t := \{ D \in \EEt \ | \
             D = \partial_s + \sum_{i=1}^\infty a_i \partial_s^{-i}\ \}
$ and $P \in \DDt$ satisfy $[ P, \tilde L]\in
\EEt_-$, the equation
$[\partial_y - P, \tilde L]=0$ is called {\it Lax equation}
and $(P,L)$ is {\it Lax pair}.

\endroster}
\endproclaim

\proclaim{4-21. Proposition}{\rm [Mul] }
Let $L:= \partial_1^2 - u \in \Cal L^t_2$.
{\it \roster

\item $ [\partial_{t_n}-2^{2(n-1)}{[L^{(2n-1)/2}}]_+,L]=0$
is the Lax equation.

\item For  an arbitrary  $P\in \DDt$ of
the Lax pair $(P,L)$, $P$ can be expressed by
$$
	P = \sum^{n}_{j=1} c_j [L^{j/2}]_+,
$$
where $c_j \in \CC[[t_2,t_3,\cdots]]$.

\item If and only if $u$ satisfies $[\partial_y - P, L^{1/2}]=0$,
 $[\partial_y - P, L]=0$.

\item The equation
$   [\partial_{t_n}-2^{2(n-1)}[{L^{(2n-1)/2}}]_+,L]=0$
gives the $n$-th KdV equation,
$\partial_{t_n}u-\Omega_1^{n-1}\partial_1 u =0$, and thus
we have a bijection
$$
\BMKdV\approx\{ L \in \Cal L_2^t \ | \
[\partial_{t_n}-2^{2(n-1)}[{L^{(2n-1)/2}}]_+,L]=0, \ n >1 \ \}.
$$
Here $\approx$ is given by the correspondence between $u$ and
$L = \partial_1^2 + u$.

\endroster }

\endproclaim

\demo{\Proof}
First we consider (3).
Let $L=W\partial_1^2W^{-1}$. $[\partial_y - P, W\partial_1 W^{-1}]=0$
gives $[W(\partial_y - P)W^{-1}, \partial_1]=0$ and then
we obtain
$[W(\partial_y - P)W^{-1}, \partial_1^2]=0$ and
$[\partial_y - P, L]=0$.
For an operator $Q \in \aac$,
$[Q,\partial_1^2]=0$ means
$(\partial_1^2 Q)+2(\partial_1 Q)\partial_1=0$,
{\it i.e.},
$(\partial_1^2 Q)=0$ and $(\partial_1 Q)=0$.
Hence
$[W(\partial_y - P)W^{-1}, \partial_1^2]=0$ also means
$[\partial_y - P, L^{1/2}]=0$.
(1) It is know that $[[L^{j/2}]_+, L^{1/2}]\in \EEt_-$.
Due to (3), (1) is proved.
Next we consider (2).
$[L^{j/2}]_+$ is a monic operator. Hence if order of $P$ is $n$,
there exists $c \in \CC$ such that the order of
$P-c [L^{n/2}]_+\in \DDt$ is $n-1$. By induction, we have the
results in (2).
Proposition 4-18 (1) leads us to (4). \qed
\enddemo

Here we will translate the relations in terms of
geometrical language.
Due to Proposition 4-21 (4), we also denote
the right hand side there by $\BMKdV$.

\proclaim{4-22. Lemma}{\rm [Mul, SN, SS]}
Let $L:= \partial_1^2 - u = W^{-1} \partial_1^2 W \in \Cal L_2^t$,
$$
d L := \partial_1 L d t_1 + \partial_2 L d t_2
+ \partial_3L d t_3 + \cdots,
$$ $$
d W := \partial_1 W d t_1 + \partial_2 W d t_2
+ \partial_3 W d t_3 + \cdots,
$$
 $$
dZ := 2 L^{1/2} d t_1 + 4 L^{3/2} d t_2
+ 8 L^{5/2}  d t_3 + \cdots,
$$ $$
dZ_+ := 2 [L^{1/2}]_+ d t_1+ 4 [L^{3/2}]_+ d t_2
+ 8 [L^{5/2}]_+  d t_3 + \cdots, \quad Z = Z_+ + Z_-.
$$
{\it \roster

\item
The Lax equation becomes
$$
	d L = [Z_+, L], \quad	d L = -[Z_-, L].
$$

\item
$$
	d Z_+ = \frac{1}{2} [ Z_+, Z_+], \quad
	d Z_- = -\frac{1}{2} [ Z_-, Z_-].
$$

\item $d L = [ dW \cdot W^{-1}, L]$.

\item $W^{-1} d W - Z_+
\in \DDc dt_1+ \DDc dt_2 + \DDc dt_3+ \cdots$
or by using the gauge freedom,
$$
	d W = Z_+ W, \quad dW =- Z_- W.
$$
\endroster}

\endproclaim

\demo{Proof}
(1) is trivial.
(2): Noting $d^2 L\equiv 0$, $[L, d Z_+ - Z_+ Z_+]\equiv 0$
and then we obtain (2).
(3): From $d (W W^{-1})\equiv 0$, $d W^{-1}=-W^{-1}dW W^{-1}$.
Hence $dL=d(W \partial_1^2 W^{-1})$ becomes the right hand side.
(4): Using (2) and (3), $[d W W^{-1} - Z_+, L]=0$, and we obtain
 $[W^{-1}d W - W^{-1} Z_+ W, \partial_1]=0$.
It implies (4).
\qed \enddemo

\tvskip

Here we note that the conditions $dZ_+ = [Z_+,Z_+]$ and so on are
the Frobenius integrability conditions. Due to the conditions,
the orbit as a dynamical system can be  uniquely determined.
Conclusively we have the following proposition on the
orbit of the KdV equations. As its proof is a little bit
complicate, we will give only the result.

\proclaim{4-23. Proposition}{\rm [Mul]}
For
$L(0)=S(0)\partial^2 S(0)^{-1}$,
$$
	U(t)=exp( t_1\partial_1+ t_2\partial_1^3
                   +t_3 \partial_1^5 +\cdots ) S(0)^{-1}\in
      \hat \EEt,
$$
$U(t) = S(t)^{-1} Y$ for $S(t) \in \frak G$ and $Y \in \hat \DDt$,
we have the time development,
$$
	L(t) =S(t) \partial_1^2 S(t)^{-1}.
$$
\endproclaim

\tvskip

\proclaim{4-24. Definition}{\rm [Mul]
\roster
\item For $L \in \Cal L_2^t$, if the map
$$
	T_0(R_{t,n}) \ni \frac{\partial}{\partial y }
          \mapsto \left.
           \frac{\partial L^{1/2}}{\partial y}\right|_{y=0}
          \in \EEt_-
$$
is injective, we say that $R_{t,n}$ is {\it effective}.
Here we write $R_{t,n}$ as the orbit space generated by
$t_1, t_2, \cdots, t_n$ and $T_0(R_{t,n})$ as
its tangent space at the origin $0 \in R_{t,n}$.

\item
If  for  $L \in \Cal L_2^t$  $R_{t,n}$ is not
effective for $n>g$ but $n\le g$ is effective,
we say that $L=\partial_1^2+u$ or $u$ is
 {\it finite $g$ type solution}
 of the KdV equation.

\endroster}
\endproclaim

\proclaim{4-25. Lemma}{\rm [Mul, S, SN]}
{\it \roster
\item
If there is  a natural number $N$ such that $\partial_{t_N}u$
is an eigen vector of the operator
$\Omega$ with an eigenvalue $k \in \CC$, {\it i.e.},
$$   k\partial_{t_N} u=\Omega \partial_{t_N}u ,
$$
$\partial_{t_m}$ is
scalar multiplication of $\partial_{t_N}$
for $m \ge N$. If not, we refer that $t_m$ is {\it effective}.

\item For the finite $g$ solution of $L$ and for $n> g$, we have
the commutation relation,
$$
[2^{2(g+1)}[{L^{(2g+1)/2}}]_+,L]\equiv0,
$$
by construction $t_n$ in
terms of  a linear combination in $\CC<t_1,t_2,\cdots,t_g>$,

\item Let $L \in \Cal L_2^t$
such that $[2^{2(g+1)}[{L^{(2g+1)/2}}]_+,L]\equiv0$ and
$$
[\partial_{t_j} -2^{2(j-1)}[{L^{(2j-1)/2}}]_+,L]=0,
\quad \text{for } j<g,
$$
is effective. Then we have a
commutative subring  $\aat := \CC[L, [L^{(2g+1)/2}]_+] \subset \DDt$
such that $W \in \WWt$, $\aac = W^{-1} \aat W \in \AAc$,
and an isomorphism as $\CC$-vector space,
$$
	\HH^1(\aac) \approx
\CC< dt_1, dt_2, \cdots, dt_g>.
$$

\endroster }
\endproclaim

\demo{\Proof}
(1) is essentially the same as Lemma 4-1 and (2) is obvious
from (1). So we will concentrate our attention on (3).
The integrability conditions makes the conditions in $\DDt$
reduced to those in $\DDf$ as an initial state.
Due to Lemma 4-14 (5), (3) is proved. \qed
\enddemo

\tvskip

\proclaim{4-26. Definition}
{\rm
\roster
\item The
filter with respect to the effective differential equations
is defined by
$$
\split
      F_g\BMKdV&:= \{L \in \Cal L_2^t \ | \
  [\partial_{n}- 2^{2(n-1)}[L^{(2n-1)/2}]_+,L]=0
 \text{ is not effctive for } n > g\} \\
      &= \{L \in \Cal L_2^t \ | \
  [[L^{(2g-1)/2}]_+,L]\equiv0\ \} \\
\endsplit
$$
and
$$
	F_g\BMKdV \subset F_{g+1}\BMKdV, \quad
	F_n\BMKdV =\emptyset, \text{ for } n<0.
$$

\item A set of finite $g$ type solutions of the KdV equation
is denoted by
 $$
\BMKdV_g:=F_{g}\BMKdV\setminus F_{g-1}\BMKdV.
$$

\endroster
}
\endproclaim

\tvskip
Due to Lemma 4-25 (3), $F_g\BMKdV$ corresponds to $F_g \AAc$ and
the correspondence becomes a bijection by considering their appropriate
quotient spaces.

As the system of the KdV equations is a dynamical system,
there is a $g$-dimensional orbit in each solution space
in $\BMKdV_g$ by neglecting its periodicity.
We can regard it as a fiber bundle,
$$
\CD       \text{orbit}  @>>> \BMKdV_g  \\@.  @VV{\pi_\KdV^{g}}V \\
                   @.     \fMKdV_g   .
\endCD
$$
For each orbit space $ {\pi_\KdV^{g}}^{-1}(p)$
at a point $p$ in $\fMKdV_g$,
there is a commutative ring $\aac \in \AAc_g$ such that
$$
         T^*     {\pi_\KdV^{g}}^{-1}(p) \approx \HH^1(\aac).
$$
For later convenience, we also define a space
$\fMKdV_{\roman{finite}}:=\coprod_{g} \fMKdV_g$.

Next we will consider $\BMKdV$ itself.
$(F_g,\AAc)$ has direct limit due to Proposition 4-15.
Let us consider the set of subrings in $\DDt$
$$
	\frak B:=\{ L \in \Cal L_2^t \ |
          \ \exists W \in \WWt, \exists \aat \in \DDt
             \text{ and }
          \ \exists \aac \in \overline\AAc \text{ such that }
                      \aat = W \aac W^{-1} \text{ and }
             L = W \partial_1^2 W^{-1}\  \}.
$$
Since solving the KdV equations are an initial value problem,
for  an arbitrary  initial state $u \in \CC[[t_1]]$
 we can find the time-development obeying the KdV
equations. Thus we have $\frak B \subset \BMKdV$.
On the other hand, from the definition, we can find
$\CC[\partial_1^2] \in \overline\AAc$ which gives
$\NN_{\CC[\partial_1^2]} = 2 \NN$.
Further for  an arbitrary
$L \in \BMKdV$, there is a gauge transformation,
$W \in \WWt$ such that $W^{-1} L W=\partial_1^2$ due to
Lemma 4-8 (2).
Hence $\frak B \supset \BMKdV$ and then
$\frak B \equiv\BMKdV$. Such a consideration
is justified by the direct limit and graded topology of $\DDt$
or $\EEc$.

Thus $\BMKdV$ has naturally the topology induced from the linear
topology of the micro-differential operator in Proposition 4-7
and the filter of $\CC[\partial_s^2]$-module in Proposition 4-15,
even though $\BMKdV$ itself is not vector space.

\proclaim{4-27. Proposition}
{\it\roster

\item
$\BMKdV$
 is a filter space.

\item The set of finite $g$ type solutions of the KdV equation
is denoted by
 $\BMKdV_g$ $:=$ $F_{g}\BMKdV\setminus F_{g-1}\BMKdV$ and
the set of finite type solutions of the KdV equation
 is denoted by $\BMKdV_\finite$.
Then we have decomposition,
$$
 \BMKdV_{\finite}= \coprod_{g < \infty} \BMKdV_g.
$$

\item
In the sense of Proposition 4-15 (3), it converges.
$$
	\BMKdV = \overline{\cup_{g=1}^\infty F_g\BMKdV}.
$$
\endroster}
\endproclaim

For a point of $\fMKdV_g$, they have non-trivial (effective)
differential equations,
$$
[\partial_{n}- 2^{2(n-1)}{L^{(2n-1)/2}}_+,L]=0,
\quad (n=1,\cdots,g).
$$
For the orbital as a dynamical system,
we find a natural volume form
$<dt_1,dt_2,\cdots,dt_g>_{\CC}$.

\proclaim{4-28. Lemma}
\roster
\item
For $L_1, L_2 \in \BMKdV_g$, there is $W \in \WWt$
such that
$$
	W L_1 W^{-1} =L_2.
$$

\item
For a subset of $\WWf$ $(g>0)$,
$$
	\WWf_g:=\{ W \in \WWf\ | \
        W L W^{-1} \in \BMKdV_g, \text{ for }L \in \BMKdV_g\ \},
$$
the projection $\pi_{KdV}^g $ along the orbit space in
the quotient space of $\BMKdV_g$ by the action of $\WWf_g$
is given by,
$$
 \pi_{KdV}^g (\BMKdV_g/\WWf_g) \sim \roman{pt}.
$$

\item For
$$
	\WWf_{0,1}:=\{ W \in \WWf\ | \
        W L W^{-1} \in \BMKdV_0\cup\BMKdV_1,
          \text{ for }L \in \BMKdV_0\cup\BMKdV_1\  \},
$$
the following relation holds,
$$
 \pi_{KdV} (\BMKdV_0\cup\BMKdV_1 /\WWf_{0,1}) \sim \roman{pt}.
$$

\endroster
\endproclaim

\demo{Proof} (1) is essentially the same as Proposition 4-16 (2).
The action of $\WWf_g$ to $\BMKdV_g$ is transitive due to (1)
and thus (2) is obtained.\qed

\tvskip
We note that $\BMKdV_0$ should be regarded as a compactification
of the base field $\CC$, {\it i.e.}, $\BMKdV_0\equiv \PP$
itself.

\tvskip

Now let us come back to the elastica problem.
Firstly we note that
the elastica problem is defined over the real functions.
Hence we should restrict the above result to a real analytic
problem.
In other words, we choose
 a natural complex structure $J$ $(J^2=-1)$ in the orbits space
$<dt_1, dt_2, \cdots, dt_g>_\CC$ and constraint
it by $<dt_1, dt_2, \cdots, dt_g>_\RR$
using the fact that finite $g$-type flow is a
finite $g$-type solution of the KdV equation.
Further the orbit satisfies the reality condition
$|\partial_1 \gamma|=1$, which characterizes a certain
type of hyperelliptic curves.

Secondly we should notice the difference of the categories of
the previous chapter and this chapter.
However  as the $g$-type flow $u$ is expressed by meromorphic
functions over a hyperelliptic curve of genus $g$,
 elements of the {\it finite real} flow
 exist in the category of the formal power series.
Hence the investigation of the finite flow
does not depend on the difference.

Further the arc-length $s$ corresponds to $t_1$ in
the above argument but we are consider only the
closed one. Hence firstly $t_1$ must be an element
of $S^1=\RR/\ZZ$. Even though $u(s)\equiv \{\gamma, s\}_{\SD}$
is periodic, $\gamma$ is not in general.
We should restrict the space of the solution space
of the KdV equation so that $\gamma(0)=\gamma(2\pi)$ or
$\gamma(0) = \gamma(\infty)$.

We will define a projective structure in $\BMeP_g$ by
$\pi_\elas^g: \BMeP_g \to \fMeP_g$ so that for a point $p$ in
$\fMeP_g$, ${\pi_\elas^g}^{-1}(p)$ is the real number orbit,
and let
$\fMeP_\finite=\coprod_{g}\fMeP_g$ as did in Definition 3-18.

We  summary these results in the following proposition.

\tvskip
\proclaim{4-29. Proposition}
{\it  There are
natural injections
$$
i_{\KdV}:\BMeP_{\text{finite}} \hookrightarrow
 \BMKdV_{\text{finite}},
\quad
\iota_{\KdV}:\fMeP_{\text{finite}}
 \hookrightarrow \fMKdV_{\text{finite}}
$$
which satisfy

\roster
\item
$$
\iota_{\KdV} \circ \pi_{\elas} =  \pi_{\KdV}\circ i_{\KdV},
$$

\item
$$    \BMKdV\setminus \BMeP \neq \emptyset .
$$
\endroster
}\endproclaim

\tvskip

Using the above results,
there is a filtration in $\BMeP$ such that
$$
	F_g\BMeP :=\BMeP \cap F_g \BMKdV,
$$
which satisfies
$$
	F_g\BMeP \subset F_{g+1}\BMeP, \quad
        \BMeP_g =F_g\BMeP/ F_{g-1}\BMeP.
$$
 We have written  just $\BMeP$ as $i_{\KdV}(\Bbb
M_{\elas}^{\PP})$ and $\fMeP$ as $\iota_{\KdV}(\frak
M_{\elas}^{\PP})$ for brevity.

Next we will consider the real orbits or the {\lq\lq}time" development
of each finite $g$ type flow in $\BMeP$ (instead of $\BMKdV$).
Let us recall the fact that rational points
in $[0,1)$, {\it i.e.}, $\QQ/\ZZ$, are measure zero in $[0,1)$.
Further it is known that for a torus $\CC/(\ZZ + \sqrt{-1}\ZZ)$,
a real direct line (orbit) stemmed from the origin with
an angle $\theta$ does not stand upon the origin again
if $\theta \not\in \tan^{-1} (\QQ/\ZZ)$.
Similarly  in general, the real number {\lq\lq}time" development
of the finite $g$ type solution is not periodic in {\lq\lq}time"
$t_i$ ($i>1$), in the $g$-dimensional torus $\Cal J_g$
which is called quasi-periodic solutions.
Hence we conclude that such an orbit is homeomorphic to $\RR^{g-1}$
in this sense and show the following proposition.

\proclaim{4-30. Proposition}
{
For each $\roman{pt}\in \BMeP_g$, we have a restricted
action of $\WWf_g$ and thus the following results are satisfied:
\roster

\item
$$
	\pi_\KdV^g|_{\BMeP_g}(\BMeP_g/[\WWf_g|_{\BMeP_g}])= \roman{pt}
$$

\item
$$
	\BMeP_g/[\WWf_g|_{\fMeP_g}]\approx S^1\times \RR^{g-1}, \quad
        \CMeP_g/[\WWf_g|_{\fMeP_g}]\approx \RR^{g-1},
        \quad \text{ for } g>1.
$$

\item
$$
 \BMeP_0\cup\BMeP_1/\WWf_{0,1}|_{\BMeP_0\cup\BMeP_1} \approx S^1, \quad
    (\BMeP_0\cup\BMeP_1/\WWf_{0,1}|_{\BMeP_0\cup\BMeP_1})/\Isom(S^1)
            \approx \roman{pt}.
$$

}
\endproclaim

\tvskip

We will recover the base ring with smooth functions.
In other words,
we  show that the completion in Proposition 4-27
can be extended to $\EEs$ because the convergence is
determined only
by the topology of order of the differential operator
as shown in the following lemma.

\vskip 0.5 cm \proclaim{4-31. Lemma}
{\it\roster
\item  When we define $\EEs_m:=
\{ D \in \EEs \ | \ \roman{deg} D \le m \ \}$,
$\EEs$ has a filter topology,
$$
   \EEs= \cup_{m} \EEs_m, \quad\{0\}= \cap_{m} \EEs_m, \quad
   \EEs_m \subset \EEs_{m+1} .
$$

\item $\EEs$ is a linear topological space with respect to this
filter topology.

\endroster
\endproclaim

\tvskip
Due to Lemma 4-31 and note  below the Proposition 4-3,
 we have the following proposition.

\proclaim{4-32. Proposition}
{\it
Let us define the moduli space of the KdV equations over the ring
of the smooth functions:
$$   {\BMKdV}^\infty:=\{ u\in \Cinf(\Cal V^\infty)
      \ |\  \partial_{t_n} u - \Omega_1^{n-1}u =0
\text{ for }\ \forall n \  \} ,
\quad   \CMKdV^\infty:= \BMKdV^\infty/(t_1) .
$$
Then

\roster
\item ${\BMKdV}$ is dense in ${\BMKdV}^\infty$.

\item ${\BMKdV}_\finite$ is a subset of ${\BMKdV}^\infty$.

\endroster

}
\endproclaim

\demo{Proof} Due to the Weierstrass  preparation theorem,
for  an arbitrary germ in $\SCinf(\RR)$,
there is a sequence in $\SF(\RR)$ conversing it.
Integrability due to Proposition 3-15
asserts that the difference does not
enlarge for the time development.
Hence (1) is proved. (2) is obvious \qed \enddemo

Hence we have the finial statement in this section.

\proclaim{4-33. Proposition}
 $\BMeP \subset \BMKdV^\infty$  has the filter topology
induced from $\BMKdV^\infty$.

\roster
\item There is a natural decomposition,
$$
 \Bbb M_{\elas,\finite}^\PP= \coprod_{g < \infty} \BMeP_g.
$$

\item
With respect to the induced topology and in the sense of
Propositions 4-27 (3) and 4-32 (1),
we have
$$
	\BMeP = \overline{\cup_{g=1}^\infty F_g\BMeP_g}.
$$

\endroster
\endproclaim

%\newpage
\tvskip
\centerline{\twobf \S 5.
 Algebro-Geometric Properties of the KdV flow II}
\centerline{\twobf--Analytic Expressions--}
\tvskip

In this section, we will give a more concrete argument.
For example, we will give another proof of Theorem 4-2
and Proposition 4-33
at Proposition 5-22.
This section is based on the inverse scattering method,
Krichever's scheme and Baker's  approach.

\vskip 0.5 cm

The studies of the KdV equation have a long history. There were so many
researchers contributing them, {\it e.g.}, Miura, Gardner,
 Greene, Kruskal,
Lax,  and so on [D, DJ]. Owing to their studies, we will give,
here, another aspect of the KdV equation without proofs,
 which is called the inverse scattering method.

\proclaim{5-1. Proposition}{\rm [BBEIM, D, DJ, Kr]}
{\it \roster \item
 For a solution $u$ of the $n$-th KdV equation, there  is
a complex valued smooth function $\psi_{\bx}$ over
$\Bbb R \times \{t_n\}$, which is a
universal covering of $S^1$ $(S^1 = \Bbb R/2\pi \Bbb Z)$ and
$\{t_n\} \subset \RR$,
$$     \psi_{\bx} \in \Cinf(\Bbb R \times \{t_n\}, \CC),
$$
as an eigen vector of the eigenvalue problem over $\Bbb R$ $(S^1 = \Bbb
R/ 2\pi \Bbb Z)$,
$$   L=\partial_1^2 + u ,\quad  -L\psi_{\bx} =\bx \psi_{\bx} ,
$$ and as a solution of  $$
   (\partial_{t_n}-2^{2(n-1)}{L^{(2n-1)/2}}_+) \psi_{\bx} =0.
$$
The deformation of $u$ with respect to $t_n$ which preserves the eigen
value $\bx$ is equivalent with that  $u$ is a solution of
 the $n$-th KdV-equation and
vice-versa. Here the extension of the domain of $u$ over $S^1$ to $\Bbb
R$ is naturally defined as $u(t_1)=u(t_1+ 2\pi)$. These equations are
equivalent to the Lax equations in Definition 4-20.

\endroster
\endproclaim

\tvskip
\proclaim{5-2. Remark}
{\rm
\roster
\item The eigenvalue problem $- L\psi_{\bx} =\bx \psi_{\bx}$ can be
regarded as a quantization  of the {\lq\lq}classical" equation
$$   (-\partial_1^2 - \frac{1}{2}\{\gamma,s\}_{\SD}(t_1) )\psi(t_1) =0 .
$$ Indeed, $ (\partial_\tau - L)\Psi =0$,  ($\Psi = \exp(\tau\bx) \psi$)
appears when we quantize $\psi(t_1)$ by means of the
 path integration
 [R,Mat0].

\item For finite type solutions of the KdV hierarchy,
the Lax equations and the
compatibility condition are essentially reduced to finite relations.
 Due to Lemma 4-1 and Definition 4-24, the equations with respect to
$t_m$ $m> N$ are trivial one for $N$-type solution.

\endroster}
\endproclaim

\tvskip
For a while, we will assume that  $u$ is real.
Let  $\Spect(-L)$ denote a set of $\bx$.
 Due to hermitian properties of $-L$, $\Spect(-L)$ is a subset
of real number bounded from below. The function $\psi_{\bx}(t_1)$ is
regarded as a section of line bundle over $\Spect(-L)$.

For bases $y_0$ and $y_1$ of the solution space of
$-L\psi_{\bx} =\bx\psi_{\bx}$,
($ \psi_{\bx} = a y_0 + b y_1$, for $a$, $b \in \CC$),
$$   y_0(0,\bx)=1, \quad y_1(0,\bx)=0, \quad
   \partial_1 y_0(0,\bx)=0, \quad\partial_1 y_1(0,\bx)=1,
$$
we have monodoromy matrix  defined as
$$     M(\bx) :=\pmatrix y_0(\pi,\bx) & y_1(\pi,\bx) \\
         \partial_1 y_0(\pi,\bx) & \partial_1 y_1(\pi,\bx) \endpmatrix ,
$$
whose determinant is unity.
If the eigenvalue of this matrix $\rho$ is in the unit circle in $\Bbb
C$ ($|\rho|=1$), the solution $\psi_{\bx}$ is called stable and exist as
a global section over the line bundle over $s \in \Bbb R$. Unless, it is
called unstable and it means that there is no global section over $s \in
\Bbb R$ even though we can find local solutions of $-L\psi_{\bx}
=\bx\psi_{\bx}$. We sometimes refer the unstable state
{\lq\lq}gap state" or {\lq\lq}forbidden state". The determinant whether it is
stable or unstable is done by the characteristic equation,
$$   \rho^2 - \Delta_u \rho +1 = 0 ,
$$
where $\Delta_u:=\tr M$. If its discriminant $\Delta_u^2 -4$ is
non-positive, corresponding $\bx$ becomes stable.

Since $\Delta_u^2-4$ is an analytic function over $\Spect(-L)-\{ \infty\}$
and has ordered zero points $\bx_1,\bx_2 \cdots$,
it has infinite product expression:
$$   (\Delta_u^2 -4) = c \prod_{j=0}^\infty  (\bx -\bx_j),
$$
where $c$ is a constant in $\bx$. This fact is correct even for the case
that $u$ is complex valued and
thus we will return to the general $u$ form here.

\proclaim{5-3. Proposition}{\rm [MM]}
{\it  For $-L\psi_{\bx} =\bx\psi_{\bx}$ with
smooth $u(t_1)$ over $\Bbb R$,
the discriminant $\Delta$ is characterized by infinite
$\bx_j$ and  can be rewritten as,
$$
   (\Delta_u^2 -4) = \left(\prod_{j=0, \text{single zeros}}
     (\bx -\bx_j) \right)h(\bx)^2 ,
$$
where $h(\bx)=\sqrt{c} \prod_{j', \text{double zeros}}^\infty
(\bx-\bx_j)$ is the part of double zeros.}
\endproclaim

\tvskip
For large $\bx$, $-L$ asymptotically behaves like
$-\partial_1^2$ for bounded
$u$ and thus  the asymptotic behavior of $\Delta$ can be investigated.
Since the ground state corresponds to a single zero of
$\Delta^2 - 4$ and
other each gap has two single zeros of $\Delta^2 - 4$, the number of
single zeros of $\Delta^2 - 4$ must be odd.
Here we will consider a case with  finite single zero points $2g+1$:
$$
 \frac{(\Delta_u^2 -4)}{h(\bx)^2} = \prod_{j=1}^{2g+1}  (\bx -\bx_j) .
$$
We refer such a case as finite-gap-state. It should be noted that
$\psi_{\bx}$ has natural involution  $\pi :\Spect(-L) \longrightarrow
\Spect(-L)$ ($\pi : \bby \longrightarrow -\bby$, $\pi : \infty = \infty$)
where $\bby = \sqrt{\Delta_u^2 -4}/h(\bx)$.
 Due to analyticity, we can extend $\Spect(-L)$ to complex.  As for
$u\equiv 0$ case, $\Spect(-\partial_1^2)$ is complexfied to $\PP$ (even
though we need more precise arguments), the energy spectrum
$\Spect(-L)$ is, in general, reduced to a hyperelliptic curve $C_g$ due to its
two-folding property. In fact for $\bby =
\sqrt{\Delta_u^2 -4}/h(\bx)$, this relation means a hyperelliptic curve
defined in  4-13.

\tvskip
In this section, we will fix a hyperelliptic curve $C_g$ with
genus $g$ given by an affine curve,
$$
	\bby^2 =h_g(\bx,1)= (\bx-c_1)\cdots (\bx- c_{2g+1}).
$$
In other words, we deal with a commutative ring
$\CC[\bx,\bby]/(\bby^2-h_g(\bx,1))\cup \{\infty\}$.
We should note that
for a hyperelliptic curve $C_g$, there exists a differential
operator $-L$ with $u$ such that its spectrum Spect$(-L)$ gives
the hyperelliptic curve  isomorphic to $C_g$.

\proclaim{5-4. Proposition}
{\it
Let the moduli space  of hyperelliptic
curves of genus $g$ be denoted by $\frak M_{\hyp,g}$.
Then $\frak M_{\hyp,g}$  is $(2g-1)$ dimensional space.  }
\endproclaim

\demo {Proof} A point in the moduli space $\frak M_{\hyp,g}$ is
characterized by $2g+1$ zero points of $h_g(x,1)$
in the above definition and $\infty$ point.
However in these variables, there are several symmetries which express
the same compact Riemannian surface. First one is
translational symmetry $c_j \to c_j + \alpha_0$, $\alpha_0 \in \Bbb C$.
Second one is dilatation $c_j \to c_j \alpha_1$ $\alpha_1 \in \Bbb C$.
Third one is $(\bx,\bby) \to (1/\bx,\bby \prod_j c_j/\bx^{(2g+1)/2})$,
which reduces
$c_j \to 1/c_j$. Hence the remainder degree of freedom is $2g-1$.
\enddemo

\tvskip
\proclaim{5-5. Remark}
{\rm
We will mention $\frak M_{\hyp,g}$ here.
We consider a smooth curve in $\frak M_{\hyp,g}$
 which is not degenerated;
$c_i \neq c_j$ if $i\neq j$ and
all of $c_j$ are finite value of $\Bbb C$.
Let us find the largest distance $|c_j-c_k|$ of pair $(c_j,c_k)$
in  $\{c_j\}$
as  an arbitrary  $|c_j-c_k|$ does not vanish
because the curve is not degenerated.
Let us rename them as $(c_1, c_{2g+1})$ and define
$$
(\alpha_1, \cdots, \alpha_{2g-1}):=
((c_2-c_1)/(c_{2g+1}-c_1), \cdots,(c_{2g}-c_1)/(c_{2g+1}-c_1))
  \in \CC^{2g-1}.
$$
Since $1 - \alpha_j = (c_{j+1}-c_{2g+1})/(c_1-c_{2g+1})$
 and $|c_1-c_{2g+1}|$ is the largest distance,
the region of each $\alpha_j$ must be constrained
as $|\alpha_j|\le 1$ and $|1-\alpha_j|\le 1$.
Next we will  order $\alpha$ following the law,
\roster
\item if $Re (\alpha_i) < Re (\alpha_j)$, $i<j$.

\item if $Re (\alpha_i) =Re (\alpha_j)$ and
$Im (\alpha_i) < Im (\alpha_j)$, $i<j$.

\endroster
Hyperelliptic curves of genus $g$ are determined as
two-fold coverings of
$\PP^1$ ramified at $0$, $1$, $\infty$ and $2g-1$ additional points as
the above order.

However it is difficult to deal with deformation from
 non-degenerate hyperelliptic curves to
degenerate curves [HM, IUN, Mum0].
}
\endproclaim

\tvskip
\proclaim{5-6. Definition} [BBEIM, IUN, Kr, Mum0-2]
{\rm \roster
\item
Let
$$
\HH_1(C_g, \Bbb Z)
  =\bigoplus_{j=1}^g\Bbb Z\alpha_{j}
   \oplus\bigoplus_{j=1}^g\Bbb Z\beta_{j}.
$$
denote the homology of an algebraic curve $C_g$.

%\item We define the Jacobi
%variety $\hat{\Cal J_g}$ associated with
%$C_g\in \frak M_{\hyp,g}$ by the exact sequence,
%$$     0 \longrightarrow \HH_1(C_g, \Bbb Z) \longrightarrow
%       \hat{\Cal J_g} \longrightarrow
%       \HH_1(C_g, \Cal O^\times) \longrightarrow 0 .
%$$

\item We introduce the periodic matrix of the
 curve $C_g$, in terms of the normalized first
kind one-form $\omega_{i}$ over $C_g$:
$$
  1=\left[\int_{\alpha_{j}}\hat \omega_{i}\right], \quad
  \Bbb T =\left[\int_{\beta_{j}}\hat \omega_{i}\right], \quad
  \hat{\pmb{\Omega_1}}=\left[\matrix 1 \\ \Bbb T \endmatrix\right].
$$

\item For fixing $\Bbb T$, we define the theta function
$ \theta:\CC^g \longrightarrow \CC$ by,
$$
     \theta(z)   : =\theta(z| \Bbb T)
     :=\sum_{n \in \Bbb Z^g} \exp \left[2\pi i\left\{
     \dfrac 12 \ ^t\negthinspace n \Bbb T n
     + \ ^t\negthinspace n z\right\}\right].
$$
\endroster
}
\endproclaim

\proclaim{5-7. Proportions}{\rm [BBEIM, IUM, Kr, Mum0-2]}
{\it
\roster \item
By defining the Abel map for $g$-th symmetric product
of the curve $C_g$,
$$       \hat w:\roman{Sym}^g( C_g) \longrightarrow \CC^g, \quad
      \left( \hat w_k(Q):=
\sum_{i=1}^g \int_\infty^{Q_i}s \hat \omega_k \right),
$$
the Jacobi variety $\hat{\Cal J_g}$ is realized as a complex torus,
$$   \hat{\Cal J_g} = \CC^g /\hat{ \pmb{\Lambda}} .
$$
Here $\hat{ \pmb{\Lambda}}$ is a lattice generated by
$\hat{\pmb{\Omega}}$.  For the Abelian group of the
divisor of the line bundle over a hyperelliptic curve
$C_g$, which is called Picard group
$\Pic^0(X)$, the Abel theorem is expressed by $\Pic^0(X) \approx
\CC^g/\hat{ \pmb{\Lambda}}$.

\item
The theta function has monodoromy properties
$$    \theta( z + e_k )= \theta(z), \quad\theta( z + \tau_k ) =
       \ee^{ - 2\pi i z_k + \pi i \tau_{kk}} \theta(z) .
$$

\item The Riemann theorem gives that
$$   \theta(\hat w(Q)-\sum_{i=1}^g \hat w(P_g) + K) \not \equiv0 ,
$$
where $K$ is a constant called Riemann constant if and only if
$P_g$'s are general points on $C_g$.

\endroster
}
\endproclaim

\tvskip
As $\Bbb M_{\elas,g}^{\PP}$ and $\Bbb M_{\KdV,g}$ have
the natural projections,
we will introduce  the universal family of hyperelliptic
 curves of
  $\frak M_{\hyp,g}$
  induced from  $\pi_\hyp: \Cal J_g \mapsto \text{pt}\in
  \frak M_{\hyp,g}$.

\tvskip
\proclaim{5-8. Proposition}{\rm (Krichever, Mulase)[Kr, Mul, Mum1]}
{\it \roster
\item A finite $g$ type solution of the KdV equation is given
by a meromorphic function over the Jacobi variety $\Cal J_g$ of a
 hyperelliptic curves $C_g$.

\item There is a natural bijection between the moduli spaces
of hyperelliptic curves $\frak
M_{\hyp,g}$ and $\frak M_{\KdV,g}$,
$$  \frak M_{\hyp,g} \approx \frak M_{\KdV,g} .
$$\endroster}
\endproclaim

\tvskip

As (2) comes from the previous section,
  we will mention its idea of (1) as follows [Kr, SW].
Krichever started with $\psi_x$, a solution of
$(-\partial_1^2 - u+x^2)\psi_x=0$, which is called the Baker-Akhiezer
function. His approach is very natural in the soliton theory and can be
generalized from the case of the KdV hierarchy, which is related to
hyperelliptic curves, to that in the KP hierarchy related to
more general compact Riemannian surfaces.

\proclaim{5-9. Lemma}{\rm [Kr]}
{\it \roster
\item
For a solution of the KdV equation whose $\Spect(-L)$ is associated with
the hyperelliptic curve $C_g$,
we parameterize the eigenvalue $-x^2$ for $L\psi_x =x^2 \psi_x$.
Then $1/x$ is a local parameter of $\infty$ of $C_g$.

\item $\psi_x$ is meromorphic on $C_g - \infty$ and at the point
$\infty$ it has an essential singularity
$$   \psi_x = \ee^{s x} \psi_W,  \quad
         \psi_W:=( 1+ \sum_{i=1}^\infty a_i(t_1) x^{-i}).
$$
Here this expansion gives us the recursive relation
$   - 2 \partial_1 a_i = - L a_{i-1} $with $a_0 = 1$.
\endroster }
\endproclaim

\demo{Proof}
(1): For a sufficiently large $|x|$,
this equation can be approximated by
$(-\partial_1^2+x^2)\psi_x\sim 0$.
Thus we can regard $\psi_x \sim \exp( s x)$. In other
 words for a local coordinate $z=1/x$ around $\infty \in \Spect(-L)$,
 $\psi_x \sim \exp( -s/z )(1 + \Cal O(z))$: $1/x^2=1/\bx$ is a local
 coordinate around  $\infty \in \Spect(-L)$.
 (2)  can be obtained by straightforward computations.
\enddemo

\tvskip
Using this Lemma 5-9, we  follow the Krichever's construction of the
finite $g$ type solution.
As we gave the Jacobi varieties and theta functions of
hyperelliptic curve $C_g$ in 4-13 and Proposition 5-7,
 we  introduce a normalized Abelian differential of the second kind,
$\hat \eta_{P,i}$,
$$   \hat \eta_{P,n} = d \left(\frac{1}{t^{n-1}} + \Cal O(1)\right),
$$
around $P$ using a local parameter $t$ ($t(P)=0$) with the normalization
$$
  \int_{\alpha_j} \hat \eta_{P,n}=0 , \quad \text{ for } j=1,\cdots,g.
$$

As we have prepare to express the Baker-Akhiezer function, we
consider the deformation equation,
$$ (\partial_{t_n}-2^{2(n-1)} L^{(2n-1)/2}_+)\psi_x =0.
$$
Since $z=1/x$ is a local parameter around $\infty$ and
around there $ L^{(2n-1)/2}_+ \sim \partial_1^{(2n-1)}$,
we introduce
$$   \hat \eta_{\infty,n} = d (x^{2n-1} + \Cal O(1)) ,
$$
and  consider the function
$$   \Cal E(t,Q)= \exp\left( \sum_{\alpha, j }
   2^{2(n-1)}t_{\alpha, j} \int^Q \eta_{P_\alpha,i} \right) .
$$

Around $\infty$, $\Cal E(t,Q) \sim \exp( \sum_{n=1}^\infty 2^{2(n-1)}t_n
x^{2n-1} )$ and
$\partial_{t_n} \Cal E(t,Q) \sim 2^{2(n-1)} x^{2n-1}\Cal
E(t,Q)$. Due to Lemma 4-8 and 5-9 and by letting
$L=W(s,\partial_1)\partial_1^2W(s,\partial_1)^{-1}$
in the sense of Lemma 4-14, we obtain the relations
$\psi_x=W(s,\partial_1) \Cal E(t,Q)
 + \Cal O\left(\dfrac{1}{x}\right)$ and
$$   L^{n/2} \ W(s,\partial_1)\  \Cal E(t,Q)= W(s,\partial_1)\
   \partial_1^n \Cal E(t,Q)
   + \Cal O\left(\frac{1}{x}\right).
$$
From the Lax equations in Proposition 5-1, $\psi_x$ is
 expressed by $\psi_x/\Cal E =
(\psi_x/\Cal E)(x t_1, 4x^3 t_2, 8x^5 t_3, \cdots)+
\Cal O(\frac{1}{x})$.

On the other hand, even though $\Cal E(t,Q)$ is satisfied with the
dispersion relation around $\infty$ and has no monodoromy around
$\alpha_j$'s, it has monodoromy around $\beta_j$
$$   \exp( 2 \pi i U_j) := \exp( \sum_{j,\alpha}
         2^{2(j-1)} t_j H_{\ \alpha,j}^i )
$$ where $$
    H_{\ \alpha,j}^i= \frac{1}{2\pi i}
    \int_{\beta_i} d  \hat \eta_{P_\alpha,j} .
$$

Noting this  monodoromy of the theta function in
Proposition 5-7, we can find a
single value
function over $\CC^g$, which is known as Baker-Akhiezer function;
$$   \psi_x = \Cal E(t,Q)
   \frac{   \theta( w(Q) +\sum_{\alpha,j}2^{2(j-1)}
   t_{\alpha,j} H_{\alpha,j} -
   \sum_{i=1}^g w(P_g) + K)}{
   \theta( w(Q)    -\sum_{i=1}^g w(P_g) + K)} .
$$
This is a solution of the Lax equations in Proposition 5-1.  We can find
a finite type solution of the KdV equation by using the zero mode
using Proposition 2-8. $\psi_x$ is determined by an analysis on
the functions
over  $\CC^g$ related to the Jacobi variety $\hat \Cal J_g$.
 As the map form $C_g$'s to
the Jacobi variety $\hat \Cal J_g$ is known as Abel map,
finding inverse map
from functions over $\hat \Cal J_g$ to functions over $C_g$'s is known
as Jacobi inverse problem. Krichever's scheme should be regarded
as the Jacobi inverse method and can be applied even to
a generalized Jacobi variety.
It shows the existence of an injection form
$\frak M_\hyp$ to $\frak M_\KdV$,

\tvskip

\tvskip\proclaim{5-10. Remark}
{\rm
For a finite type solution of the KdV hierarchy $u$,
we have the hyperelliptic
curve $C_g$ as a spectrum of $-L$ to $u$. Then the above arguments give
the following results:

\roster\item
The orbit of the equations of the KdV hierarchy is realized in a direct
line in the Jacobi variety $J_g$ of $C_g$.

\item Any finite $g$ solution $u$
 is given as a solution of the Jacobi inverse problem of the
Jacobi variety $J_g$.
\endroster}
\endproclaim

\tvskip
As far as we will deal
with only hyperelliptic curves and the KdV hierarchy, we can give more
concrete arguments based upon Baker's original argument
 [Ba1, Ba2].

\tvskip

\proclaim{5-11. Definition}{\rm [Ba1, Ba2, BEL]
We introduce the family of the differential forms:
 \roster
\item
$$   \omega_1 = \frac{ d \bx}{2\bby}, \quad
      \omega_2 =  \frac{\bx d \bx}{2\bby}, \quad \cdots \quad
     \omega_g =\frac{\bx^{g-1} d \bx}{2 \bby}.
$$

\item$$
     \eta_{j}=\dfrac{1}{2\bby}\sum_{k=j}^{2g-j}(k+1-j)
      \lambda_{k+1+j} \bx^k d\bx ,
     \quad (j=1, \cdots, g) .
$$\endroster}
\endproclaim

\tvskip
\proclaim{5-12. Lemma}{\rm [Ba1, Ba2, BEL]}
{\it \roster
\item $\omega$'s are the basis of the holomorphic function valued
cohomology of  hyperelliptic curve $C_g$,
which give unnormalized periods:
$$    \pmb{\Omega}'=\left[\int_{\alpha_{j}}\omega_{i}\right], \quad
      \pmb{\Omega}''=\left[\int_{\beta_{j}}\omega_{i}\right], \quad
    \pmb{\Omega}=\left[\matrix \pmb{\Omega}' \\ \pmb{\Omega}''
     \endmatrix\right].
$$

\item They are related to the normalized ones:
$$    \ ^t\left[\matrix \hw_{1}  \cdots & \hw_{g}\endmatrix\right]
       :={\pmb{\Omega}'}^{-1}  \ ^t\left[\matrix \omega_{1} & \cdots
   \omega_{g}\endmatrix\right] ,\quad
   \Bbb T:={\pmb{\Omega}'}^{-1}\pmb{\Omega}''.
$$

\item $\eta$'s are the unnormalized one-form of the second kind over
$C_g$ and then the complete hyperelliptic integral of
the second kinds is given  as
$$      H':=\left[\dint_{\alpha_{j}}\eta_{i}\right], \quad
         H'':=\left[\dint_{\beta_{j}}\eta_{i}\right] .
$$
\endroster}
\endproclaim
Here the contours in the integral are, for example,
given in p.3.83 in [Mum2].

\demo{Proof} We check holomophicity of the forms in  (1) and (3).
A zero point of $\bby=0$, or a root $c_j$ of $f(\bx)=0$, corresponds to
a point  $(c_j,0)$ of the curve $C_g$. We  use a local coordinate
$z^2:=(\bx-c_j)$ and
 $\bx^m d\bx/(2 \bby) \sim (z^2 + c_j)^m dz + \cdots$.
On the other hand, around $\infty$ point,
 let us choose local coordinate $1/x$ as $1/x^2 = 1/\bx$ and
 then  $\bx^m d\bx/(2 \bby) \sim  (1/x)^{2g-2m+2} dx+ \cdots$.
Hence $\omega$ is holomorphic all over the curve $C_g$ while
$\eta$ is holomorphic except $\infty$ point.\qed \enddemo

\tvskip

\proclaim{5-13. Definition}
{\rm\roster
\item The unnormalized Jacobi variety $\Cal J_g$ is defined
by a complex torus,
$$   \Cal J_g = \Bbb C^g / \pmb{\Lambda} ,
$$
where $\pmb{\Lambda}$ is a lattice generated by $\pmb{\Omega}$.

\item We defined the theta function by,
$$
\theta\negthinspace\left[\matrix a \\ b \endmatrix\right](z)
   =\theta\negthinspace\left[\matrix a \\ b \endmatrix\right](z; \Bbb T)
   =\sum_{n \in \Bbb Z^g} \exp \left[2\pi i\left\{
    \dfrac 12 \ ^t\negthinspace (n+a)\Bbb T(n+a)
    + \ ^t\negthinspace (n+a)(z+b)\right\}\right].
$$
\endroster}
\endproclaim

\tvskip
\proclaim{5-14. Proposition}
{\it The Riemann constant of the hyperelliptic curve
$C_g$ is given as
$$
K=\sum_{j=1}^g\int_{\infty}^{{\mrm A}_j}\hw = \delta' + \delta''\Bbb T
$$
where $ \delta' =\ ^t\left[\matrix \dfrac g2 & \dfrac{g-1}2 & \cdots
      & \dfrac 12\endmatrix\right],
   \quad \delta''=\ ^t\left[\matrix \dfrac12 & \cdots & \dfrac12
   \endmatrix\right].$
 }
\endproclaim

\demo {Proof} This proof is in p.3.82 in [Mum2]. \enddemo

\tvskip
Using the Abel map $\roman{Sym}^g( C_g) \longrightarrow \Bbb C^g$,
we define the coordinate in $\Bbb C^g$,
$$
   \ft_j :=\sum_{i=1}^g \int^{(\bby_i,\bx_i)} \omega_j .
$$
Here we  note that $\ft_j$ behaves $(1/x)^{2(g-j)+1}$ around
$\infty$ point if we use the parameter $x^2 =\bx$.

\vskip 0.5 cm
\proclaim{5-15. Definition} {\rm ($\wp$-function, Baker)[Ba1,Ba2]
 \roster
\item Using the coordinate $\ft_j$, the $\sigma$-function, which is
a holomorphic function over $\Bbb C^g$, is defined by
$$
\sigma(\ft):=\sigma(\ft;C_g):
  =\ \roman{exp}(-\dfrac{1}{2}\ ^t\ \ft
  H'{\pmb{\Omega}'}^{-1}\ft)
  \vartheta\negthinspace
  \left[\matrix \delta'' \\ \delta' \endmatrix\right]
  ({\pmb{\Omega}'}^{-1}\ft ;\Bbb T) .
$$

\item
In terms of the $\sigma$-function, the
hyperelliptic $\wp$-function over the hyperelliptic curve
$C_g$ is defined by
$$   \wp_{ij}(\ft):=-\dfrac{\partial^2}{\partial \ft_i\partial \ft_j}
   \log \sigma(u)  =\dfrac{\sigma_i(\ft)\sigma_j(\ft)-\sigma_{ij}
   (\ft)\sigma(\ft)}{\sigma(\ft)^2}.
$$

\endroster}
\endproclaim

As $\sigma$-function is an entire function over $\CC^g$ and has
single zero at $g-1$ dimensional subvariety of $\CC^g$
which is called theta-divisor, the hyperelliptic $\wp_{ij}$ has
the second order singularity and function of $\Cal J_g$.

\tvskip

\proclaim{5-16. Remark }
{\rm
It is worth while noting that from Definition 5-15, the hyperelliptic
$\wp$-function can be concretely computed for a
given hyperelliptic curve $C_g$. The summation in the definition
of $\theta$ function rapidly converges due to effects of $\Bbb T$
  and others are integrations of primary functions.
Further it is known that
$\wp_{gi}$ is an elementary symmetric function,
{\it i.e.},
for $F(\bx) = (\bx-\bx_1) (\bx- \bx_2) \cdots (\bx-\bx_g)$,
[Ba1, BEL],
$$
	F(\bx) = \bx^g - \sum_{i=1}^g \wp_{g i}\bx^{i-1}.
$$
Accordingly, by numerical approach,
we can compute a value of the hyperelliptic $\wp$
function as Euler determined a
value of the elliptic integral to know the shape of a
classical elastica by numerical method [E, L, T1, 2].
This approach was discovered by Baker about one hundred years ago
 [Ba1, 2, Mat7-10].

We emphasize that it completely differs from
Krichever's approach based upon Baker-Akhiezer theorem explained in \S 4.
Krichever's arguments might not give us
practical algorithms to fix parameters
 of general hyperelliptic function except  solutions expressed
by elliptic or hyperbolic functions. (Due to its
abstract, it is a good strategy to construct soliton theory.)

On the other hand, Baker's original method determines  concrete function
forms of corresponding $\wp$ functions,
for any algebraically given hyperelliptic
curves (even for degenerate curves in $\frak M_{\hyp,g}$).
We can expand $\wp$-function
around a general point and know its parameter dependence.

 Since this Baker's construction in \cite{Ba2}
might be no longer in recent researchers' memory
 as long as I know,
we believe that this review of Baker's work has meaning.
We believe that it is very useful for the analysis
in physics [BEL, Ma7-10].

}
\endproclaim

\vskip 0.5 cm
In [Ba2] Baker found that the $\wp$-functions obey
the following differential equations, which contain
the KdV hierarchy.

\proclaim{5-17. Example }(genus = 3){\rm [Ba2]}
Let us express $\wp_{ijk}:=\partial \wp_{ij}(\ft) /\partial \ft_k$ and
$\wp_{ijkl}:=\partial^2 \wp_{ij}(\ft) /\partial \ft_k \partial \ft_l$.
The hyperelliptic $\wp$-function obeys the relations
$$   \allowdisplaybreaks \align
     ( 1)\quad& \wp_{3333}-6\wp_{33}^2
      =  2\lambda_5 \lambda_7
       + 4\lambda_6 \wp_{33} + 4 \lambda_7 \wp_{32},\\
     ( 2)\quad& \wp_{3332}-6\wp_{33} \wp_{32}
      =  4\lambda_6 \wp_{32} + 2 \lambda_7 (3\wp_{31}-\wp_{22}),\\
     ( 3)\quad& \wp_{3331}-6\wp_{31}\wp_{33}
      =  4\lambda_6 \wp_{31} - 2 \lambda_7 \wp_{21},\\
     ( 4)\quad& \wp_{3322}-4\wp_{32}^2-2\wp_{33} \wp_{22}
    =  2\lambda_5 \wp_{32} + 4 \lambda_6 \wp_{31} - 2
                              \lambda_7 \wp_{21},\\
     ( 5)\quad& \wp_{3321}-2\wp_{33} \wp_{21}-4\wp_{32} \wp_{31}
      =  2\lambda_5 \wp_{31},\\
     ( 6)\quad& \wp_{3311}-4\wp_{31}^2-2\wp_{33} \wp_{11}
      = 2 \Delta_\wp, \\
     ( 7)\quad& \wp_{3222}-6\wp_{32}\wp_{22}
      = -4\lambda_2 \lambda_7
        -2\lambda_3 \wp_{33} + 4 \lambda_4 \wp_{32}
        +4\lambda_5 \wp_{31} - 6 \lambda_7 \wp_{11},\\
     ( 8)\quad& \wp_{3221}-4\wp_{32} \wp_{21}-2\wp_{31} \wp_{22}
      =- 2\lambda_1 \lambda_7+ 4 \lambda_4 \wp_{31} - 2  \Delta_\wp, \\
     ( 9)\quad& \wp_{3211}-4\wp_{31} \wp_{21}-2\wp_{32} \wp_{11}
      =- 4\lambda_0 \lambda_7+ 2 \lambda_3 \wp_{31},\\
     (10)\quad& \wp_{3111}-6\wp_{31} \wp_{11}
      =  4\lambda_0 \wp_{33} - 2 \lambda_1 \wp_{32} + 4
                  \lambda_2 \wp_{31}, \\
     (11)\quad& \wp_{2222}-6\wp_{22}^2 \\
      =-8&\lambda_2 \lambda_6+ 2 \lambda_3 \lambda_5
                      - 6 \lambda_1 \lambda_7
       -12\lambda_2 \wp_{33} + 4 \lambda_3 \wp_{32}
                        + 4 \lambda_4 \wp_{22}
       + 4\lambda_5 \wp_{21} -12 \lambda_6 \wp_{11} +12 \Delta_\wp,\\
     (12)\quad& \wp_{2221}-6\wp_{22} \wp_{21}
      =- 4\lambda_1 \lambda_6- 8 \lambda_0 \lambda_7
       - 6\lambda_1 \wp_{33} + 4 \lambda_3 \wp_{31}
                  + 4 \lambda_4 \wp_{21}
       - 2\lambda_5 \wp_{11}, \\
     (13)\quad& \wp_{2211}-4\wp_{21}^2        -2\wp_{22} \wp_{11}
      =- 8\lambda_0 \lambda_6
       - 8\lambda_0 \wp_{33}
             - 2 \lambda_1 \wp_{32} + 4 \lambda_2 \wp_{31}
       + 2\lambda_3 \wp_{21}, \\
     (14)\quad& \wp_{2111}-6\wp_{21} \wp_{11}
      =- 2\lambda_0 \lambda_5
       - 8\lambda_0 \wp_{32} + 2 \lambda_1(3\wp_{31} - \wp_{22})
       + 4\lambda_2 \wp_{21}, \\
     (15)\quad& \wp_{1111}-6\wp_{11}^2
      =- 4\lambda_0 \lambda_4+ 2 \lambda_1 \lambda_3
       + 4\lambda_0(4\wp_{31}- 3 \wp_{22})
       + 4\lambda_1 \wp_{21} + 4 \lambda_2 \wp_{11}.
     \endalign
$$ where $$
     \Delta_\wp  = \wp_{32} \wp_{21} - \wp_{31} \wp_{22}
     + \wp_{31}^2 - \wp_{33} \wp_{11}.
$$
\endproclaim

\proclaim{5-18. Proposition}
{\it For $u=-2(\wp_{gg}-\lambda_{2g}/3 )$ and
$u(s,t_2,t_3)=u(\frak t_g, \dfrac{\ft_{g-1}}{2^2},
\dfrac{\ft_{g-2}}{2^4} + \dfrac{3}{2^4\lambda_{2g}} \ft_{g-1})$
obeys the first and the second KdV equations.}
\endproclaim

\demo {Proof} Let us consider $g=3$ case.
If we regarded as $u=-2(\wp_{33}-\lambda_6/3 )$,
 it is obvious that (1) in Example 3-14 becomes the KdV equation
noting $\lambda_7 = 1$. By setting $2 \partial_{\frak
t_3} \times (2) + \partial_{\frak t_2}\times (1)$ and $\partial_{t_3}=
16\partial_{\ft_1}+ \dfrac{16\lambda_2}{3}
 \partial_{\ft_2}$, we obtain the
second KdV equation. From arguments of Baker [Ba1, Ba2], even for $g>3$ the
relations (1) and (2) maintain for $g$ case. \qed \enddemo

\tvskip\proclaim{5-19. Remark}{\rm
\roster
\item
By above arguments,
for given hyperelliptic curve $\bby^2 =f(\bx)$, we can
construct a solution of the first and second KdV equations.
 Further the
compatibility of Lax system gives more general argument for the other
equations in the KdV hierarchy.
Then it implies that we explicitly showed
the existence of an injective  map
$$   \frak M_{\hyp,g} \longrightarrow \frak M_{\KdV,g} .
$$
This correspondence is valid even for degenerate curves.

s\item Our development of the quantized elastica after submitting
this article is in [Mat7-10]. In [Mat7, 10], we showed more
explicit function forms of quantized elastica over $\CC$.

\item
After submitting this article, we knew the works of
Buchstaber, Enolskii, Leykin and related people
[reference in [BEL]]
as an extension of parts of Baker's studies.

\endroster}
\endproclaim

\tvskip
Now let us give another proof of Theorem 4-2 (2) and Proposition
4-33 (2).

\proclaim{5-20 Proposition}
Proposition 4-2 (2), 4-27 (3), and 4-33 (2)
 can be regarded as an approximation theory
based upon the Weierstrass  preparation theorem.
\endproclaim

\demo{\Proof}
Let us recall the moduli space of the
KdV equations whose base ring is
smooth functions and definition is in Proposition 4-32,
$$   {\BMKdV}^\infty=\{ u\in \Cinf(\Cal V^\infty)
      \ |\  \partial_{t_n} u - \Omega_1^{n-1}u =0
\text{ for }\ \forall n \  \} ,
\quad   \CMKdV= \BMKdV/(t_1) ,
$$
For  an arbitrary  $u \in{\BMKdV}^\infty$, there is a unique
spectrum $\Spect(L:=-\partial_s^2 -u)$
 up to its orbits, by solving the eigenvalue
equation $(-\partial_s^2 -u)\psi_x = \bx \psi_x$.
We assume that the spectrum
does not have finite gap $\{ (c_1,c_2),$ $ (c_3,c_4),$ $ \cdots, $
$(c_{2g-1},c_{2g}),$ $ \cdots\}$ and then the corresponding
characteristic equation becomes transcendental equation
$\bby^2=f(x)$, where $f(x)$ is the
transcendental function with zeros $(c_j)_{j=0,1,\cdots}$.
Since  $u$ is a  smooth function over $S^1$, $|u|$ is bounded the
above. Hence around $\infty$ of the $\Spect(L)$, $L \sim - \partial_s^2$
and $\Spect(L)$ at $\infty$ is patched by the affine space $\CC$;
the width of gap converges to zero for
$\bx \to \infty$. Thus we can approximate $\Spect(L)$ by finite gap
spectrum $\Spect(L_g):=$ $\{ (c_1,c_2),$ $ (c_3,c_4), $ $ \cdots,$
$(c_{2g-1},c_{2g}),$ $ (c_{2g+1},\infty)\}$. The approximated potential
$u_g$ is given by the $\wp$ function of the hyperelliptic function
$\bby^2 = h_g(\bx,1)$ whose zero points are $(c_j)_{j=1,2,\cdots,2g+1}$.
By using Weierstrass  preparation theorem and taking
appropriate $g$,
we can approximate $f(\bx)$ by $h_g(\bx,1)$ for desired.

Hence up to the KdVH flow, $u_g$ approaches to $u$ for $g$ approaches
to $\infty$ from its construction. (For  an arbitrary
 finite $g$, $u_g$ is unique up to
its orbits). Thus for  an arbitrary $u$ in ${\BMKdV}^\infty$,
there is a series of points $u_g$ belonging to $F_g{\BMKdV}$
such as $u_g \to u$ for $g \to
\infty$ up to orbit.
(We note that the finite type solutions does not depend upon
the base rings $\Cinf$ or formal power series.)
Hence we have
$$
{\BMKdV}^\infty = \overline{\cup_g F_g{\BMKdV}}.
$$

Since $\BMeP$ is a subset of ${\BMKdV}^\infty$
and for  an arbitrary  curve $\gamma \in \CMeP$,
$u:=\{\gamma,s\}_{\SD}$ has a unique value,
the above statement is valid. \qed
\enddemo

\tvskip

\proclaim{5-21 Example}{\rm [E, T1, 2, Mat2]
As an element $\gamma$ in $\BMeP$ must satisfy
$\gamma(s+L)=\gamma(s)$ in $\PP$ and a reality condition
$|\partial_s \gamma |=1$.
Even though the hyperelliptic function $\wp$ is
a meromorphic function over $\Cal J_g$, we can find
a trajectory or real line in $\Cal J_g$ which avoids
the singularities and satisfies the reality and closed
conditions. In other words, we will find
$\BMeP$ as a subset of $\BMKdV$.
We give  examples of the $\gamma\in \PP$ in terms of
the local chart around the origin.
\roster

\item genus $g=0$ case: a circle with radius 1.

\item genus $g=1$ cases: $\CMeP_1$ consists of two points:

2-1) Jacobi elliptic modulus $l=1$ case
$$
         \gamma(s)=s - \frac{2}{\alpha} \left(\tanh(\alpha s)
                 -\sqrt{-1}\sinh(\alpha s ) \right).
$$

2-2) Jacobi elliptic modulus $l=0.908911\cdots$, which gives
the eight shape loop in a complex plane $\CC$ [E].

\endroster
\endproclaim

Here we note that in [Mum3],
Mumford  gave simple and deep expression of
the shape of elastica, which shows the depth, importance
and beauty of this problem. There he showed
 how the reality condition $|\partial_s \gamma|=1$
restricts the moduli of  elliptic curves.

%\newpage
\tvskip
\centerline{\twobf \S 6. Cohomology of a Loop Space }
\tvskip

As we mentioned in Introduction, in this section,
we will digress from
our analysis of the moduli of a quantized elastica and review arguments
of a loop space over $S^2$ in the category of topological spaces $\Top$ whose
morphism is a continuous map (isomorphism is homeomorphism,
monomorphism is injective continuous map and so on). Studies on a loop
space in $\Top$ are well-established and its cohomological properties
are well-known as in the textbook of Bott and Tu [BT]. We can recognize
the moduli space of a quantized elastica in $\PP$ as a loop space in the
category of the differential geometry $\DGeom$. When we replace
smooth functions with continuous functions and $\PP$ with $S^2$
respectively,  it is expected that the moduli space
of a quantized elastica in $\PP$ is related to that in $\Top$.
In this section, we will review
a loop space in $\Top$ and show its cohomological properties.

\proclaim{6-1 Definition}{\rm [BT]}
{\rm
$E$ and $X$ are topological space and $X$ has a good cover $\frak U$.
A map $\pi : E \longrightarrow X$ is called a
{\it fibering} if it satisfies the
covering homotopy properties: for given  a map $f:Y \longrightarrow E$
from an arbitrary
 topological space $Y$ into $E$ and homotopy
$\overline f_t$ of $\overline f = \pi \circ f$ in $X$
($Y \times [0,1] \longrightarrow X$, $f_0 := f$), there
is a homotopy $f_t$ of $f$ which covers
$\overline f_t$; ($Y \times [0,1]
\longrightarrow E$ such that $\overline f_t:=\pi \circ f_t$).}
\endproclaim

\tvskip\proclaim{6-2. Definition}{\rm [BT]}
{\rm\roster
\item The {\it path space} of $X$ is defined to be the space $P(X)$
consisting of all the paths in $X$ with initial point $*$:
$$
  P(X):=\{ \text{maps } \mu :[0,1] \longrightarrow X\
                                | \ \mu(0)=* \in X\ \}.
$$

\item The {\it loop space} over $X$ with a fixed point is defined by,
$$   \Omega X = \{ \mu :[0,1] \longrightarrow X|\ \mu(0) = \mu(1)
   = * \in X\}.
$$\endroster}
\endproclaim

\tvskip
In the category of topological spaces $\Top$, $\PP$ and  $S^2$ are
identified by homeomorphism as its morphism.
Thus we will give properties of the loop space over $S^2$ in
$\Top$ as follows.

\proclaim{6-3. Theorem}{\rm [BT]}
{\it\roster \item $P(S^2)$ is  a fibering
whose fiber is $\Omega(S^2)$:
$$\CD        \Omega(S^2)  @>>> P(S^2) \\@.  @VVV \\
                   @.    S^2   \\
\endCD$$

\item Its cohomology is torsionless and given by
$$   \HH^q( \Omega S^2,\Bbb Z) = \Bbb Z \quad
         \text{  for } q \in \Bbb Z_{\ge0}.
$$
as a module and its algebraic properties are given by
$$    \HH^*( \Omega S^2, \Bbb Z)
          = \frak E(x) \otimes_{\Bbb Z} \frak Z_\gamma(e) ,
$$
where $x$ and $e$ generators of $\HH^{1}( \Omega S^2,\Bbb Z)$ and
$\HH^{2}( \Omega S^2,\Bbb Z)$ respectively (
${\roman {dim}}\ x =1$ and ${\roman {dim}}\ e=2$). Here $\frak
E(x)$ is the exterior algebra $\Bbb Z[x]/(x^2)$ and $\frak Z_\gamma(e)$
is the divided polynomial algebra whose base is
$(1,e,e^2/2,e^3/3!,\cdots)$.
 In other words, the generator of
$ \HH^{2k+1}( \Omega S^2,\Bbb Z)$ is $x \cdot e^k/k! $ and that of
$ \HH^{2k}( \Omega S^2,\Bbb Z)$ is $e^k/k!$. \endroster}
\endproclaim

\tvskip
In order to prove Theorem 6-3, we prepare two well-known results in
algebraic topology and triangle category without proofs [BT].

\tvskip\proclaim{6-4. Proposition}{\rm [BT]}
{\it For given a double complex $K=\oplus_{q,p \ge 0} K^{p,q}$,
there is a spectral sequence $\{  \roman{E}_r, d_r\}$
converging to the total
cohomology $ \HH_D(K)$ such that each $ \roman{E}_r$
has a bigrading with
$$  d_r :  \roman{E}_r^{p,q} \longrightarrow  \roman{E}_r^{p+r,q-r+1}
$$ and $$
    \roman{E}_1^{p,q} =  \HH_d^{p,q}(K), \quad
    \roman{E}_2^{p,q} =  \HH_\delta^{p,q} \HH_d(K),
$$
where $d$ and $\delta$ are derivative:
$d:K^{p,q} \longrightarrow K^{p+1,q}$
and $\delta :K^{p,q} \longrightarrow K^{p,q+1}$, $D= d+ (-)^p \delta$. }
\endproclaim

\tvskip
We will consider the double complex for a fibering
$\pi: E \longrightarrow M$,
$$
   K^{p,q}:= C^p(\pi^{-1} \frak U, \Omega^q).
$$
Here $\frak U$ is a ramification of $M$ and $\Omega^q$ is a $q$-form
along the fiber.

\tvskip\proclaim{6-5. Proposition}{\rm (Leray-Hirsch theorem) [BT]}
{\it $\pi:E \longrightarrow X$ is a fibering with fiber $F$ over simply
connected topological space which has a good cover,
$$    \roman{E}^{p,q}_2 = \HH^p(X,   \HH^q(F,A)) ,
$$
where $A$ is a commutative ring.
If $ \HH^q(F,A)$ is a finitely generated $A$-module,
$$    \roman{E}_2:=  \HH^*(X;A) \otimes  \HH^*(F;A) .
$$}
\endproclaim

\tvskip
\tvskip\proclaim{Proof of Theorem 6-3}{\rm [BT]}
Since $P(X)$ is contractive,
$$    \HH^q(P(X)) = \left\{\matrix  \Bbb Z &\text{ for } q=1\\
   0 &\text{ otherwise } \endmatrix \right.
$$
and the spectral sequence must converge to $ \HH^p(P(X))$,
$ \roman{E}_2$ must give isomorphism except 0-dimension.
$$            \roman{E}_2:
        \quad \matrix {}_5 &\vdots & \vdots &\vdots &\vdots &\ddots \\
       {}_4 &\Bbb Z & 0 &\Bbb Z &0 &\cdots \\
       {}_3 &\Bbb Z & 0 &\Bbb Z &0 &\cdots \\
       {}_2 & \Bbb Z & 0 &\Bbb Z &0 &\cdots \\
       {}_1 &\Bbb Z & 0 &\Bbb Z &0 &\cdots \\
       {}_0 & \Bbb Z & 0 &\Bbb Z &0 &\cdots \\
       \ & {}_0      & {}_1 & {}_2 & {}_3 & \cdots \endmatrix
$$

Next we will consider the algebraic properties.
From Proposition 6-5, $ \roman{E}_2$  is the
tensor product $ \HH^*(\Omega S^2) \otimes  \HH^*(S^2)$.
Let $v$ be a two-form of $S^2$.
Then if $ \HH^1(\Omega S^2)$ is denoted as $\Bbb Z x$,
$ \roman{E}_2^{0,1}$ is expressed by $\Bbb Z x \otimes 1$.
The derivative  $d_2$
in $ \roman{E}_2$, which is isomorphism, acts on $x\otimes1$
as $d_2(x \otimes 1)=(1
\otimes v)$. Since $d_2(x^2 \otimes 1)=$ $(d_2 x \otimes 1) \cdot x
\otimes 1- x  \otimes 1\cdot d_2 x \otimes 1 $ $ = (1 \otimes v)(x
\otimes 1) - (x \otimes 1) (1 \otimes v)=0$, we have $x^2=0$ because
$d_2$ is isomorphism. Thus $d_2^{-1}(x \otimes v)$ is expressed by
another generator $e$ in $ \HH^2(\Omega S^2)$, which is algebraically
independent of $x$. $d_2 ( e \otimes 1) = (x \otimes v)$. Since
$d(ex\otimes 1) =e\otimes v$, $ex$ is a generator in dimension 3.
Similarly $d_2(e^2 \otimes 1)= 2ex \otimes v$ means that $e^2/2$ is a
generator in dimension 4. In other words, we have a table such that,
$$       \roman{E}_2:
   \quad \matrix {}_5 &\vdots & \vdots &\vdots &\vdots &\ddots \\
    {}_4 & e^2/2 \otimes 1 & 0  &0 &0 &\cdots \\
    {}_3 & ex \otimes 1 & 0 & ex \otimes v &0 &\cdots \\
    {}_2 & e \otimes 1& 0 & e \otimes v &0 &\cdots \\
    {}_1 & x\otimes 1 & 0 &x\otimes v &0 &\cdots \\
    {}_0 & 1 & 0 & 1\otimes v &0 &\cdots \\
    \ & {}_0      & {}_1 & {}_2 & {}_3 & \cdots \endmatrix
$$
Hence Theorem 6-3 is proved. \qed
\endproclaim

\tvskip

\proclaim{6-6. Remark}{\rm [Br]
Thought we showed the result on the loop space defined in
Definition 6-2.
However there are several studies on another loop space
$$
\{ \gamma: S^1 \hookrightarrow S^2\ | \ \text{ smooth immersion} \},
$$
and its cohomology, which differs from the result in Theorem 6-3
[Br]. This loop has a freedom of choice of starting points
of $S^1$ in $S^2$.
However in this article, we are concerned with a loop space
with fixed point as we mentioned in Remark 2-15 and Definition 2-2.
Accordingly we mentioned only the result.

}
\endproclaim

%\newpage
\tvskip
\centerline{\twobf \S 7. Topological Properties of Moduli $\Cal
M_{\elas}^{\PP}$ }
\tvskip

As in previous section,
 we reviewed the cohomological properties of a loop
space in $\Top$, in this section we will argue its relation to our loop
space in $\DGeom$ or the moduli space of a quantized elastica
again.
We believe that
such considerations are important for the quantization of
 an elastica and the
statistical mechanics of polymer physics [KL, Mat1, Mat2, Mat3].

The loop spaces in both {\bf  Top
} and $\DGeom$ are infinite dimensional spaces when we regard
them as manifolds in an appropriate sense.
Even though it is not known that
de Rham's theorem can be applicable to such an infinite dimensional
manifold, it is expected that
cohomological sequences should correspond to each other.

Precisely speaking, as we will show later, the closed condition
and the reality condition $|\partial_s\gamma|=1$ in the
moduli space $\CMeP$ makes its topological properties difficult.
Thus we must tune the 0-dimension of the cohomology related to $\BMeP$.
Then we will reach our second main Theorem 7-4, which implies
that cohomology of $\BMKdV$ reproduces Theorem 6-3 with $\RR$
coefficients.

\tvskip

Since the loop space in $\Top$ is given with the fixed point,
there is no translation freedom for the loop in $S^2$, which corresponds
to our situation of quotient of $\EE^0(\CC)$ in Definition 2-2.
Further there is no freedom of change of the origin of the loop
in $\Top$.
Hence we must compare $\Omega S^2$ with $\CMeP$ rather than $\BMeP$.

Further $\CMeP_0 \approx \CMeP_1/\WWt_1 \approx {\pt}$,
which should be regarded the same class because both these are
zero dimension.
The $\Cal F\CMeP$ might be natural sequence:
$$
\Cal F\CMeP\ : \
\emptyset \to  F_1{\CMeP} \hookrightarrow  F_2{\CMeP}
\hookrightarrow
\cdots \hookrightarrow   F_{g-1}{\CMeP}
 \hookrightarrow  F_g{\CMeP}\hookrightarrow
F_{g+1}{\CMeP}
\hookrightarrow  \cdots.
$$

Noting $ {\CMeP}_g:=F_g{\CMeP}/ F_{g+1}{\CMeP}$,
as we are concerned only with its topological properties,
let us consider the related complex of vector spaces,
$$
\split
\Cal G\CMeP\ : \
\emptyset {\overset\delta\to\longrightarrow}
 F_1{\CMeP}/\WWt_{0,1}
{\overset\delta\to\longrightarrow}{\CMeP}_{2}/\WWt_{2}
{\overset\delta\to\longrightarrow}
&\cdots
{\overset\delta\to\longrightarrow}  {\CMeP}_{g-1}/\WWt_{g-1}\\
&{\overset\delta\to\longrightarrow} {\CMeP}_g/\WWt_{g}
{\overset\delta\to\longrightarrow}
 {\CMeP}_{g+1}/\WWt_{g+1}
{\overset\delta\to\longrightarrow} \cdots,
\endsplit
$$
with trivial map $\delta=0$ and $\delta^2 =0$.

As each ${\CMeP}_g/\WWt_{g}$ is a finite dimensional vector space
$\RR^{g-1}$ thanks to Proposition 4-30,
we have de Rham complex  $\DCMeP_g$ ($g>1$),
$$
\DCMeP_g\ : \
0 \to\Omega^0({\CMeP}_g/\WWt_{g}) {\overset d\to\longrightarrow}
\Omega^1({\CMeP}_g/\WWt_{g}) {\overset d\to\longrightarrow}
\Omega^2({\CMeP}_g/\WWt_{g}) {\overset d\to\longrightarrow}
\cdots.
$$
and $\DCMeP_1$
$$
\DCMeP_1\ : \
0\to\Omega^0(F_1{\CMeP}/\WWt_{0,1}) {\overset d\to\longrightarrow}
\Omega^1(F_1{\CMeP}/\WWt_{0,1}) {\overset d\to\longrightarrow}
\Omega^2(F_1{\CMeP}/\WWt_{0,1}) {\overset d\to\longrightarrow}
\cdots,
$$
where $\Omega^p(M)$ is the set of $p$-forms over $M$.

\proclaim{7-1. Proposition}
Let us consider a double complex $\CCMeP$ with
the derivative $D=d+(-)^{g} \delta$,
$$
0\to \DCMeP_1 \to \DCMeP_{2} \to
\cdots \to  \DCMeP_{g-1}
 \to \DCMeP_g\to
 \DCMeP_{g+1} \to \cdots.
$$
Then its cohomology,
$$
	\HH^p(\CCMeP):=\oplus_g \HH^{p-g+1}(\DCMeP_g),
$$
is given by $\HH^0(\CCMeP):=\RR$ and
$$
	\HH^p(\CCMeP)=
        \RR dt_2\wedge dt_3 \wedge \cdots \wedge dt_{p+1},
         \quad p>0.
$$
\endproclaim

\demo{\Proof}
First we note
$$
	\CMeP_g/\WWt\approx\RR^{g-1}, \text{ for } g\ge 1,
         \quad \CMeP_1/\WWt\approx\pt.
$$
Since we have for $n\ge0$ [BT],
$$
	\HH^p(\RR^n) =\RR \quad \text{for } p =0.
$$
Due to Poincar\'e duality, we have
$$
	\HH^p(\RR^n)=\HH^{n-p}_c(\RR^n),
$$
if we write the compact support function valued cohomology
by $\Hc^{p}$ [BT].  The generator is expressed by,
$$
 dt_2\wedge dt_3 \wedge \cdots \wedge dt_g,
$$
with a compact support function over there.
\qed\enddemo

\tvskip
First from Proposition 3-11 (5), let us interpret
$\Omega:\partial_{t_n}\mapsto
\partial_{t_{n+1}}$ as an endomorphism of  tangent space of
Jacobi varieties $T_*J_g$ of a hyperelliptic
curve related to a point $\gamma$ in $\BMeP$.
 Since the Jacobi variety is a quotient space of $\CC^g$,
its tangent space (and also its cotangent space) can be identified
with $\CC^g$: $T^*J_g \approx T_* J_g \approx \CC^g$.
Of course, we are concerned only with its real part $\RR^g$.
 Then using the canonical duality in the real part $\RR^g$,
$$   < \partial_{t_n}, d t_m> = \delta_{n,m},
$$
we can introduce an endomorphism $\Omega^{-1*}$ and $\Omega^*$
of $\BMeP$,
$$   \Omega^*: d t_n \mapsto   d t_{n-1}=\Omega^* d t_n ,
    \quad \Omega^{-1*}: d t_n \mapsto   d t_{n+1}=\Omega^{-1*} d t_n ,
$$ where
$ < \Omega \partial_{t_n}, d t_m >=<\partial_{t_n},\Omega^*  d t_m>$.

\vskip 0.5 cm
\proclaim{7-2. Definition}
{\rm Let us define an endomorphism $\epsilon$  of $\BMeP$ by,
$$
  \epsilon:= d t_2 \Omega^{-1*},
$$
where $\Omega^{-1*} $ is regarded as a right action operator,
$\epsilon^q = d t_2 \Omega^{-1*}(\wedge
\epsilon^{q-1})$ for $q>1$ and $\Omega^{-1*} \cdot
 1 := 1$. }
\endproclaim

Then we have the properties of $\epsilon$ as follows.
\proclaim{7-3. Lemma}
{\it\roster
\item We have the relation $ \epsilon^q \cdot 1 = d t_2 \wedge  d t_3
\wedge \cdots\wedge d t_{q+1}$.

\item $\epsilon$ can be realized by $\tilde\epsilon$,
$$   \tilde\epsilon :=\sigma \sum_{k=0} \epsilon_k, \quad
       \epsilon_0 := d  t_2, \quad
     \epsilon_k :=d t_{k+2} \wedge (d t_{k+1}  i_{ \partial_{t_{k+1}}})
        \quad (k>0),
$$
where $\sigma$ is a permutation operator $
\pmatrix 1 & 2& 3& \cdots& q-1 & q \\
                q & q-1& q-2& \cdots& 2 & 1 \endpmatrix$
and $i_{ \partial_{t_k}}$ is an inner product $*$ operator; $i_{
\partial_{t_k}}\cdot d t_l =<  \partial_{t_k}, d t_l>=\delta^k_l$.

\item There is a ring isomorphism,
$$
\varphi_0: \RR \otimes_\RR\frak
E(x) \otimes_\RR \frak Z_\gamma(e) \to \RR[[ \epsilon^{2},d t_2]]
$$
 by
$$
\varphi_0: (e,x) \to (\epsilon^{2},d t_2),
$$
where the product in $\RR[[ \epsilon^{2},d t_2]]$
 is defined by
$$
\epsilon * d t_2= d t_2 *\epsilon:= \epsilon\cdot d t_2,
\quad
 \epsilon * \epsilon= \epsilon^2, \quad d t_2* d t_2= d t_2\wedge d t_2 =s0.
$$

\endroster}
\endproclaim

\demo {Proof}  (1): for example $\epsilon^2 \cdot1 $ $=dt_2 \Omega^{-1*}
(\wedge dt_2 \Omega^{-1*})\cdot 1$ $=dt_2\wedge d t_3$ and this
can be extended to general case. (2): noting  $\epsilon_k^2 = 0,$ $
(k \ge 0)$, straightforward computations gives the results. (3): noting
Theorem 6-3, it is obvious.\qed \enddemo

\tvskip

 Here we will note that
$\epsilon : \HH_c^g(\RR^g) \to \HH_c^g(\RR^{g+1})$
generates the sequence $\RR^g \hookrightarrow \RR^{g+1}$,
and thus $\epsilon^m$
could be regarded as a generator of the filter topology of
$\BMeP$ and $\BMKdV$.  Thus it means that
we can evaluate the moduli space of a quantized elastica $\BMeP$
 using the induced topology and
$\epsilon$ as in Proposition 4-26.

Finally we  reach our third main theorem.
\proclaim{7-4. Theorem}
{\it
 By setting $e = \epsilon^2$, $x =  d t_2$, the cohomology
$\HH^q( \CCMeP)$,
is a ring isomorphic to $\HH^q(\Omega S^2,\Bbb R)$,
$$
 \phi : \HH^*(\CCMeP)\ \
{\widetilde\longrightarrow}\ \ \HH^*(\Omega S^2,\Bbb R).
$$
}
\endproclaim

\tvskip
\proclaim{7-5. Remark}{\rm
\roster
\item
The closed condition $\gamma(s+L)=\gamma(s)$ for some $L$
and the reality condition $|\partial_s\gamma|=1$
are  too strong.
For example due to the condition,
$\CMeP_0$ and $\CMeP_1$ consist only of disjoint points
as mentioned in Examples 5-21 and [Mat2].
Thus if we assign real vector bases each point,
these cohomology might be $\HH^p(\CMeP_0)=\RR\delta_{0,p}$
 and $\HH^p(\CMeP_1)=\RR\oplus\RR \delta_{0,p}$.
These phenomena come from a {\lq\lq}elasticity" in the category
 $\DGeom$ but we wish to consider the topological properties of
the loop space in $\DGeom$.
Thus we have replaced $\CMeP$ with $\CCMeP$ by loosing
strongness of the condition and make its topology weak;
it implies a replacement to fewer open sets.
This replacement comes from modulo computations
in the gauge transformation by $\WWt$ in the  KdV equations,
using the natural immersion $i_\KdV:\BMeP \to \BMKdV$.

However for a sufficiently large $g$ case, the closed condition
and the reality condition might not have serious effects. Then the
quotient by gauge transformation can be also guaranteed by
the fact that each moduli space of compact Riemannian surface
of genus $g$ is simply connected [McLac].
Accordingly we consider that the replacement is not worse.

The isomorphism
$\phi$ could be regarded as a functor between triangle
categories  of loop spaces in $\DGeom$ and $\Top$
and a quasi-isomorphism between
 $\CCMeP$ and $\Omega S^2$ [Br].
(These objects in the triangle categories are vector spaces given by
$\epsilon^n$ and $(x^a,e^m)$ respectively.
The morphisms are multiplications as their ring structures.)

\item From the definition,
$\epsilon^m$ can be regarded as a map from
$\HH^q(\CCMeP)$
to $\HH^{q+m}(\CCMeP)$.
We should regard that this map comes from the properties of
vertex operator, which change the genus of curves [DJKM]
and $\epsilon^m \cdot 1$ is
 interpreted as a topological base of $\CCMeP$.

\item
The operator $\epsilon$ induces the complexes $\Cal F\CMeP$ and
$\Cal G\CMeP$.
This essentially exhibits the topology of Sato theory
because in Sato theory [S, SN, SS],
the existence of the gauge transformation
$\WWt$ is a key factor.
Theorem 7-4  means that its topology is as strong as that of a
 loop space
in $\Top$. It implies that the topology of Sato theory is too weak to
lead us to express fine structure of the moduli space
as Harris and Morrison pointed out in [HM, p.44-5];
they stated that the geometrical approach in [Mul]
does not influence the study of the moduli space of algebraic curves
including $\fMhyp$ [HM].
In fact as mentioned in [HM, Mum0],
the moduli space of $\fMhyp$ is, in general,
 very complicate but our approach is not so difficult.
Accordingly we wish to obtain stronger topology to express the
moduli space.
We hope that the day comes  that the studies on
quantized elastica are connected with those of $\fMhyp$
as Euler did for the case of genus one [E, T1, 2].

\item
As we will comment in Remark 8-9,
the correspondence between loop spaces in
$\DGeom$ and $\Top$ can be extended to
higher dimensional loop spaces by considering recent result of
a quantized elastica in $\Bbb R^3$ [Mat4].

\endroster}
\endproclaim

%\newpage \tvskip
\centerline{\twobf \S 8. Discussion }
\tvskip
\proclaim{ 8-1}{\rm
Although we have correspondence between homological properties of
 $\Omega S^2$ in $\Top$ and those of $\Cal
M_{\elas,g}^{\PP}$ in  $\DGeom$,
there is open problem for a correspondence
of homotopy group between them, {\it e.g.},
$$   \pi_{q-1}(\Omega S^2) = \pi_{q} (S^2 )  \quad (q \ge 2) ,
$$ $$
   \pi_{q-1}(\Omega S^2) \times \Bbb Q
   =\left\{\matrix \Bbb Q & \text{ for } q=1,2\\
                   0 & \text{ otherwise } \endmatrix \right. .
$$
}
\endproclaim

\vskip 0.5 cm\proclaim{ 8-2}{\rm [Mat2]
We will consider $\gamma \in \BMeC$ in this remark.
By defining
$$
   v=\left( \frac{\partial_s^2 \gamma}
                {2\sqrt{-1} \partial_s \gamma} \right) ,
$$
this problem is related to the quantization of an elastica in $\CC$,
$$
 Z[\beta] = \int_{\CMeC} D \gamma
 \exp( -\beta \int_{S^1} v^2 \dd s ) .
$$
For $\beta>0$, the domain of $E=\int_{S^1} v^2 $ can be extended to
$\infty$-point and we will define
$$
 \overline{\CMeC} = \{ \gamma: S^1 \longrightarrow \CC \ |\
   \gamma \text{ is continuous,}\quad |\partial_s \gamma(s)|=1
    \}/\sim .
$$
In other words, as we assign the energy of $\gamma$ with wild shape
to $\infty$-point of $E$,
it does not contribute the partition function $Z$.
Then we can regard the partition function as
$$
  Z : \overline{\CMeC}\times \RR_{\ge0} \longrightarrow \Bbb R .
$$
The integral region in $Z$ is recognized as $\overline{\CMeC}$.
Due to our Theorem 3-4, we have a natural projection operator $\Pi_E$:
$$
   \Pi_E :
     \Cal M_{\elas}^{\CC} \longrightarrow  \Cal M_{\elas,E}^{\CC},
   \quad \Pi_E^2 = \Pi_E .
$$
We  have a spectral decomposition,
$$   1_{\Cal M_{\elas}^{\CC}} = \int d E \Pi_E .
$$

Hence the partition function becomes
$$   Z[\beta] =
 \int d  E \  \roman{Vol}( \Cal M_{\elas,E}^{\CC}) e^{-\beta E},
$$
where Vol$( \Cal M_{\elas,E}^{\CC})$ means the volume of $(  \Cal
M_{\elas,E}^{\CC})$.

Here we will comment on a question why we can use the concept of the
orbits of {\lq\lq}kinematic" system even though in the
noncommutative algebra,
one sometimes encounters nonsense of concept of orbit,
{\it e.g.}, Kronecker
foliation [C]. Even in quantized problem, we can go on to
use the concept of orbit and commutative geometry even though the
dimension of the orbit space need not be finite.

Let extend to the domain of $\beta \in \Bbb R_{\ge 0}$ to
$\Bbb R_{\ge 0}+\infty$. Note that as the inverse image,
$$   \Cal M_{\elas,\cls}^{\CC} = Z^{-1} ( Z(\infty) ) ,
$$
the classical moduli space of the harmonic map of the elastica
depending upon
the boundary condition is naturally immersed in our moduli space
$\overline{\CMeC}$. In other words, our analysis
naturally contains Euler's perspective of the classical elastica
[E, T1, 2, L].
}
\endproclaim

\vskip 0.5 cm\proclaim{ 8-3}{\rm
Due to the projection operator,
we can define the order in the moduli
space $\CMeP$. Noting that the energy
$E$ is real in $\Cal
M_\elas^{\CC}$, let
$$
  \Cal M_{\elas,<E}^{\CC} := \coprod_{E' < E} \Cal M_{\elas,E'}^{\CC} .
$$
For $E_1 < E_2$, we have
$$     \Cal M_{\elas,<E_1}^{\CC} \subset
     \Cal M_{\elas,<E_2}^{\CC}.
$$
Then the moduli space $ \Cal M_{\elas}^{\CC}$ is an ordered space.
}
\endproclaim

\vskip 0.5 cm\proclaim{ 8-4}{\rm
The operator $\epsilon$ in Lemma 7-3 can be regarded as a creation
operator in the quantum field theory.
The vacuum state is regarded as $1$.
We can define the dual space of
$\Cal V^\infty$; $<e^m , e_n> = \delta^m_n$
where $e_n = d t_n$ and $e^m = \partial_{t_m}$.

Further by noting  $\epsilon$ modulo $\epsilon^2$, we can reconstruct
$\Cal C\CMeP$ in Lemma 4-31. On the other
hand, we can introduce the micro-differential operator $e^m$
 $(m \ZZ)$ as the base of $\Cal C\CMeP$ as in the Definition 3-1 and
Proposition 7-1. Then as the dual of $\Cal C\CMeP$, we can define
$e_m$ $(m <0)$ and the vacuum of
this field operator in the quantum field
theory has affine structure as physicists think.
}
\endproclaim

\vskip 0.5 cm\proclaim{ 8-5}{\rm
In the differential operator ring, $\DDs$, the integral
$\int_{S^1} \partial_s u = 0$
means that since the integral is linear map,
its kernel belongs to  $\DDs/\partial_s \DDs$.

Using the Definition 3-7 and Proposition 3-11,
 let us define,
$$
     h := \sum_j h_j d t^j, \quad \delta :=\sum_j d t^j \partial_{t_j},
       \quad \frak a= u d s .
$$
we have the transformation in $(\DDs/\partial_s \DDs)$:
$$   \delta \frak a = \tilde \Omega h, \quad
   \tilde \Omega := d s \partial_s \frac{\delta }{\delta u},
$$ $$   \delta * h = 0. $$
This relation  is called Becchi-Rouet-Stora (BRS) relation
 [LO, Mat2].
}
\endproclaim

\vskip 0.5 cm\proclaim{ 8-6}{\rm
We will introduce a dilatation flow
$$   \partial_t \psi_x = t \partial_s \psi_x .
$$

The intersection between this flow and the KdV flow
is governed by the Painlev\'e equation of the first kind,
$$   s = 3 u^2 + \partial_s^2 u.
$$
This statement can be proved as follows. Since the KdV flow in
Remark 3-8 is given by
$B_1 = u$ while this flow $B_1=t$. Hence $u=t$ and the KdV flow becomes
$$   \partial_t u=1 = \partial_s( 3 u^2 + \partial_s^2 u),
$$
and we obtain the Painlev\'e equation of the first kind [In, Mat2].
}\endproclaim

\vskip 0.5 cm\proclaim{ 8-7}{\rm
Since the Schwarz derivative $u$ is invariant for $\PSL(\CC)$ and
$\PSL(\CC)$ transitively acts upon $\PP$, we can regard $\Cal
M_{\elas}^{\PP}$ as
$$   \Omega\SL(\CC) :=\{\gamma
   :S^1\hookrightarrow  \PSL(\CC)\ | \ \gamma(0)=1 \}.
$$
Because $\gamma(s)= g_s \gamma(0)$ for $g\in
\PSL(\CC)$, we have the condition $g(0)= g(2 \pi)$.
As Witten pointed out, for a loop space we can naturally construct its
tangent space as a loop space of the tangent space of the target space
[Wi].
In other words, we can naturally define a loop algebra
$\Omega\roman{sl}_2(\CC)$.  In the loop algebra, we have only the
condition $g^{-1} d g(0) =g^{-1} d g(2 \pi)$ using $g \in
\SL(\CC)$, which is not  stronger condition than the condition
$g(0) =g(2 \pi)$. Since there is  a
smooth map from $S^1$ to $S^1$ as $\roman{Diff}(S^1)$, we obtain an
expression of the loop algebra,
$$    \roman{Diff}(S^1) \otimes \roman{sl}_2(\CC)\oplus  \CC,
$$
which acts upon $\BMKdV^\infty$ in Proposition 4-32
with the weaker condition.

In fact, the KdV flow has bi-hamiltonian structure and 2-cocycle
$$   \omega_{\underline \Omega}(X,Y)
   :=\omega(\underline \Omega X, Y)
   + \omega(X, \underline \Omega Y).
$$
Using ordinary functional derivative (Gatuex derivative $\delta
u(y)/\delta u(x)= \delta(x-y)$), we can write down the (second) Poisson
relation,
$$    \{u(s),u(s')\}=\underline \Omega \delta(s -s') ,
$$
where $\delta(s)$ is the Dirac $\delta$-function.

Let
$$    l_n := \frac{1}{2\pi}
            \int \dd s\ u_\kappa e^{\ii s n}              .
$$
denote its Fourier component.
Then it obeys the semi-classical Virasoro algebra,
$$       \{l_n,l_m\}=(n-m)l_{n+m}+n(n^2-1)\delta_{n+m,0}           .
$$
where the second term the unit central charge.
We have the Virasoro algebra.

Using the topological relation $\CC^* \sim S^1$, the problem of
conformal field theory is reduced to that of the loop algebra. Thus our
relation can be also interpreted in the regime of
the conformal field theory.
Thus it is clear that our problem is  related to the two dimensional
quantum gravity [HM].
}
\endproclaim

\vskip 0.5 cm\proclaim{ 8-8}{\rm
It is known that for
$H_0 := \int u d s$, the second Poisson structure of $H_0$
reproduces  the KdV equation; when the second Poisson
bracket is defined as
$$
   \{X,Y\}_{\underline \Omega}=
   \omega_{\underline \Omega}(X,Y),
$$
 $\partial_t u = \{u, H_0\}_{\underline \Omega}$
is $\partial_t u + 6 u \partial_s u + \partial_s^3 u=0$.

If we will used the Hamiltonian $H_n$ of the higher dimensional KdV as
the energy functional of the system, we will have another decomposition,
$$   \Cal M_{\elas}^{\PP, (n)} = \coprod \Cal M_{\elas,E}^{\PP, (n)}
$$
$$       \Cal M_{\elas,E}^{\PP, (n)}
           :=\left\{ \gamma_t \in \CMeP  |\ H_n-E=0\ \right\} .
$$
The space is determined by the $n(>1)$-th KdV hierarchy,
}
\endproclaim

\vskip 0.5 cm
\proclaim{ 8-9}{\rm [BT, Mat4]
According to the results in [BT], we have the relation
$$
   \text{H}^q( \Omega S^n, \Bbb Z ) = \Bbb Z  \quad
   \text{ for } q=0 \ \text{ modulo } n-1.
$$
As we mentioned in Remark 7-5,
it is expected that the moduli space of a quantized
elastica in $S^{n}$ has similar cohomological properties.
In fact, one of
these authors calculated the quantized elastica in $\Bbb R^n$
and obtained the same structure of the moduli space of a
quantized elastica in $\Bbb R^n$ [Mat4].
}\endproclaim

\vskip 0.5 cm\proclaim{ 8-10}{\rm
We wish to know the volume of each $\Cal M_{\elas,E}^{\CC}$. However
this problem is not easy. In fact
as pointed out in [HM], the soliton theory might not affect to get
any information of the structure of $\Cal M_{\elas,E}^{\CC}$.

In other words, our Theorem 7-4 means that the filter topology in
 the
soliton theory is too week and is equivalent with the
topological properties of the loop space. It might have no
 effect on the study of geometrical future
of moduli space of hyperelliptic curve.
Thus we believe that we must go beyond
the ordinary soliton theories
 to another theoretical world for the study of
moduli space of a quantized elastica
as Euler investigated the elliptic functions
by studying the shape of classical elasticas [E, T1, 2 ,L, We].
}\endproclaim

\vskip 0.5 cm
\proclaim{ 8-11}{\rm
First we will note
 the relations for $\PP$, $\CC$ and upper half complex plane $\Bbb H$;
$$   \matrix
   \PP    & :&\PSL(\CC) &:&\dfrac{a \gamma + b}{c \gamma + d}\\
   \CC    &: &\         &:&a \gamma + b\\
   \Bbb H & :& \PSL(\RR)&:&\dfrac{a \gamma + b}{c \gamma + d}
   \endmatrix
$$
We showed that loops on $\PP$ are related to the KdV flow and that loops
on $\CC$ are related to the MKdV flow. Next we should consider loops on
$\Bbb H$.
}
\endproclaim

\tvskip
\proclaim{ 8-12}{\rm
One of  solutions of
$$   ( - \partial_s^2 -\frac{1}{2}\{\gamma,s\}_\SD ) \psi =0
$$
is given by $1/\sqrt{\partial_s \gamma}$. The coordinate transformation
for the Diff$(S^1)$ leads us to redefine $\psi$ as the invariant form
$\sqrt{d s /d \gamma}$.
This reminds us of the prime form and the Dirac field which
has a half weight as same as the theta function [Mum2, Mat7, 9].

In fact, for a curve in $\CC \subset \PP$, there is a natural topology of
$\gamma$ induced form  the distance in $\CC$, which is given by the
Frenet-Serret relation:
$$   \pmatrix \partial_s & k/2 \\
      -k/2 & \partial_s \endpmatrix
   \pmatrix1/\sqrt{\partial_s \gamma}\\
   i/\sqrt{\partial_s \gamma}
   \endpmatrix =0 .
$$
This operator is regarded as the Dirac operator.
The Dirac operator could
be regarded as a translator from the category of analysis
 to the category of geometry.
Hence as we are dealing with the topology of the Dirac operator,
we might have a stronger topology of the curve.

We can extend this structure to a conformal surface in $\Bbb R^3$ as the
generalized Weierstrass relation [Kno, KL, Mat5, 6].

We note that this Dirac operator
(and the Schr{\"o}dinger operator in Proposition
2-8) defined upon the loop space differs from the
Dirac operator of Witten in [Wi] because
Witten's one is related to the conformal field theory and the
ordinary string which is determined by intrinsic properties whereas ours
are related to the extrinsic Polyakov string [KL, Mat5, 6].
}\endproclaim

\tvskip
\proclaim{ 8-13}{\rm
As we noticed in  8-2, the partition function $Z$ can be expressed
by
$$
	Z = \int d E \roman{Vol}( \Cal M_{\elas,E}^\CC ) e^{-\beta E},
$$
where $\roman{Vol}( \Cal M_{\elas,E}^\CC )$ is formally represented by
$$
     \roman{Vol}( \Cal M_{\elas,E}^\CC )=\sum_g
     \int_{\frak M_{\elas,E,g}^{\CC}} d \roman{vol}(J)
\int_J d t_2 d t_3 \cdots d t_g,
$$
where $d \roman{vol}(J)$ is the volume form around a point $J$ in
$\frak M_{\elas,E,g}^{\CC}$
and $\frak M_{\elas,E,g}^{\CC}:=\frak M_{\elas,E}^{\CC}\cap\frak
M_{\elas,g}^{\CC}$.
Then we will leave integral over $t_2$,
in the above expression and obtain the
time $t_2$ depending partition function,
$$
Z[t_2] = \int d E \sum_g\int_{\frak M_{\elas,E,g}^{\CC}}
d \roman{vol}(J) \int_J d t_3 \cdots d t_g e^{-\beta E}.
$$
Similarly we obtain $Z[t_2,t_3,\cdots,t_g]$,
which is a generating function [R].
 Then we  can expect that it might obey the KdV
equation or related equation. This situation might be related to
with Witten's conjecture and Kontsevich's theorem [HM].
}\endproclaim

%\newpage

\end
\Refs
\widestnumber\key{BBEIM}
\ref \key AM \by R.~Abraham and J.~E.~Marsden \book Foundations of
 Mechanics second ed.\publ Addison-Wesley \yr1985
\publaddr Reading\endref
\ref \key Ba1 \by H.~F.~Baker \book Abelian Functions \publ Cambridge
  \yr1897 \publaddr Cambridge \endref
\ref \key Ba2 \bysame \jour Acta Math. \vol 27 \pages 135-156
 \yr 1903 \paper On a system of differential equations leading
 to periodic functions \endref
\ref \key Ba3 \bysame \paper
Note the foregoing paper
{\lq\lq}Commutative Ordinary Differential Operators" by
 J.~L.~Burchnall and T.~W.~Chaundy
  \jour Proc. Royal Society London (A)
  \yr    1928 \vol 118 \pages 584-593
\endref
\ref \key BBEIM \by E.~D.~Belokolos, A.~I.~Bobenko, V.~Z.~Enol'skii,
  A.~R.~Its and V.~B.~Matveev
\book Algebro-Geometric Approach to Nonlinear
 Integrable Equations \publ Springer \yr 1994 \publaddr New York \endref
\ref \key   {BC1}
  \by    J.~L.~Burchnall and T.~W.~Chaundy
  \paper Commutative Ordinary Differential Operators
  \jour Proc. Royal Society London (A)
  \yr    1928 \vol 118 \pages 557-583
\endref
\ref\key   {BC2}
  \by    J.~L.~Burchnall and T.~W.~Chaundy
  \paper Commutative Ordinary Differential Operators II
  \jour Proc. Royal Society London (A)
  \yr    1931 \vol 134 \pages 471-485
\endref
\ref
  \key   {BEL}
  \by    V.H. Buchstaber, V.Z.  Enolskii, and D.V. Leykin
  \paper Klein Function, Hyperelliptic Jacobians and
         Applications
  \jour Rev. Math. \& Math. Phys.
  \yr    1997 \vol 10 \pages 3-120
\endref
\ref \key Br \by J-L. Brylinski \book Loop Spaces Characteristic
  Classes and Geometric Quantization \publ Birkh\"auser \yr1992
  \publaddr Boston \endref
\ref \key BT \by R.~Bott and L.~W.~Tu \book Differential Form in
  Algebraic Topology \publ  Springer \yr1982 \publaddr New York \endref
\ref \key C \by A.~Connes \book Noncommutative Geometry
 \publ Academic Press \yr1994 \publaddr
   Singapore \endref
\ref \key D \by L.~A.~Dickery \book Soliton Equations and Hamiltonian
Systems \publ  World Scientific \yr1991 \publaddr Singapore \endref

\ref \key DJKM \by E. Date, M. Jimbo, M. Kashiwara  and T. Miwa \inbook
Nonlinear Integrable Systems -Classical Thoery and Quantum Thoery- \ed
M. Jimbo and T. Miwa  \publ World Scientific\yr 1983 \publaddr
Singapore \endref
\ref \key E \by L. Euler
\yr 1744 \book Methodus inveniendi lineas curvas
maximi minimive proprietate gaudentes \publ Lausanne
 \endref

\ref \key DJ \by P.~G.~Drazin  and R.~S.~Johnson
\yr 1989 \book Solitons: an introduction
 \publaddr Cambridge \publ Cambridge University Press \endref
\ref \key GP1  \by R.~E.~Goldstein and D.~M.~Petrich \jour Phys. Rev.
Lett.\vol  67  \yr 1991 \page 3203-3206  \paper The Korteweg-de Vries
hierarchy as dynamics of closed curves in the plane \endref
\ref \key GP2 \bysame \jour Phys. Rev.
Lett.\vol  67  \yr 1992 \page 555-558 \paper Solitons, Euler's equation,
and vortex patch dynamics \endref
\ref \key G \by M.~A.~Guest \book Harmonic Maps, Loop Groups,
and Integrable Systems (London Math. Soc. Student Text 38)
\publ  Cambridge Univ. Press \yr1997 \publaddr
   Cambridge \endref
\ref \key Ha \by R. Hartshorne \book Algebraic Geometry
\publ Springer \yr1977 \publaddr
 Berlin \endref
\ref \key HM \by J.~Harris and I.~Morrison,
   \book Moduli of Curves
  \publ  Springer \yr1998 \publaddr New York \endref
\ref \key IUN \by S.~Iitaka, K.~Ueno and Y.~Namikawa
   \book Sprits of Deescartes and Algebraic Geometry (in Japanese)
  \publ  Nihon-Hyouron-Sha \yr1980 \publaddr Tokyo \endref
\ref \key In \by E.~L.~Ince \book Ordinary Differential Equations
\publ Dover  \yr 1956 \publaddr New York \endref
\ref \key KV \by A.~L.~Kholodenko and T.~A.~Vilgis
\jour Phys. Rep. \vol 298 \yr1998
\page 251-370
\paper Some geometrical and topological problems in polymer physics
\endref
\ref \key Kr \by I.~M.~Krichever \jour Russian Math. Surverys \vol 32
  \yr 1977 \page 185-213 \paper Methods of Algebraic Geomtery in the
  Theory of Non-linear Equations \endref
\ref \key Kno \by B.~G.~Knopelchenko
\jour Studies in Appl. Math. \yr 1996
   \vol 96 \page 9-51 \paper Induced Surfaces and Their
   Integrable Dynamics \endref

\ref \key KL \by B.~G.~Knopelchenko and G.~Landlfi \yr 1998
  \jour math.DG/9804144
\paper Generalized Weierstrass representation for
  surface in multidimensional Riemann spaces  \endref

 \ref \key LP \by J.~Langer and R.~Perline
 \paper Poisson Geometry of the Filament Equation
 \jour  J.~Nonlinear Sci.  \vol 1 \yr 1991 \page 71-91 \endref

\ref \key L \by A.~E.~H.~Love
\book A Treatise on the Mathematical Theory
   of Elasticity \publ Cambridge Univ. Press \yr 1927
   \publaddr Cambridge \endref
\ref \key LO \by J.~M.~Leinass and K.~Olaussen
 \jour  Phys. Lett.  \yr 1982
   \vol 108B \page 199-202 \paper  Ghosts  and Geometry \endref
\ref \key McLac \by C.~Maclachlan \jour Proc.A.M.S.  \yr 1971
   \vol 29 \page 85-86 \paper Modulus space is simply-connected \endref
\ref \key McLau \by C.~MacLaughlin \jour Pac. J. Math  \yr 1992
   \vol 155:1 \page 143-156
   \paper Orientation and string structures on loop spaces \endref
\ref \key MM \by H. P. McKean and P. van Moerbeke
\paper The spectrum of Hill's equation
\jour Inventions math. \vol 30 \pages 217-274 \yr 1975
\endref
\ref \key Mat0 \by S.~Matsutani
       \paper The Physical Realization of the
       Jimbo-Miwa Theory of the Modified Korteweg-de Vries
  Equation on a Thin Elastic Rod: Fermionic Theory
   \jour Int.  J. Mod. Phys. A \vol 10 \yr 1995 \pages 3091-3107 \endref

\ref \key Mat1 \bysame \jour  Int. J. Mod. Phys. A  \yr 1995
   \vol 22 \page 3109-3123 \paper Geomtrical Construction of The Hirota
   Bilinear Form of the Modified Korteweg-de Vries Equation on a Thin
   Elastic Rod: Bosonic Classical Theory \endref
\ref \key Mat2 \bysame  \jour  J.Phys.A  \yr 1998 \vol 31
   \page 2705-25 \paper Statistical mechanics of elastica on plane:
   origin of MKdV hierarchy \endref
\ref \key Mat3 \bysame \jour  J.Phys.A  \yr 1998 \vol 31
   \page 3595-3606 \paper On Density of State of Quantized Willmore
   Surface:A Way to a Quantized Extrinsic String in $R^3$ \endref
\ref \key Mat4 \bysame  \jour  J. Geom. Phys.  \yr 1999
   \paper Statistical mechanics of elastica in $\Bbb R^3$ \vol 29 \pages
   243-259 \endref
\ref \key Mat5 \bysame
   \paper Dirac Operator of a Conformal Surface Immersed in $\Bbb R^4$:
   Further Generalized Weierstrass Relation
   \jour Rev. Math. Phys. \vol 12 \yr 2000\pages 431-444 \endref
\ref \key Mat6 \bysame \jour
   Rev. Math. Phys.\yr 1999 \paper Immersion Anomaly of Dirac Operator
    on Surface in $\Bbb R^3$ \vol 11 \pages 171-186\endref
\ref\key Mat7
  \bysame
      \paper         Closed Loop Solitons and Sigma Functions:
         Classical and Quantized Elasticas with Genera One and Two
          \jour J. Geom. Phys.
          \yr 2001 \vol 39 \pages 50-61\endref
\ref\key Mat8 \bysame
         \paper Hyperelliptic Solutions of KdV and KP equations:
        Reevaluation of Baker's Study on Hyperelliptic Sigma
        Functions
          \jour  J. Phys. A \yr 2001\vol 34 \pages 4721-4732
        \endref
\ref\key Mat9 \bysame
         \paper Hyperelliptic Loop Solitons with  Genus $g$:
           Investigations of a Quantized Elastica
          \jour J. Geom. Phys. \yr 2002 \vol 43 \pages 146-162
        \endref
\ref\key Mat10 \bysame
         \paper Explicit Hyperelliptic Solutions of
Modified Korteweg-de Vries Equation:
Essentials of Miura Transformation
          \jour   J. Phys. A. Math \& Gen
        \yr 2002 \vol 35 \pages 4321-4333
        \endref
\ref \key Mul \by M.~Mulase
\jour  J. Diff. Geom. \yr 1984 \page 403-430
    \paper Cohomological Structure in Soliton Equations and Jacobian
    Varieties \endref
\ref \key Mum0 \by D.~Mumford \book Curves and Their Jacobians
    \publ  Univ.\ of Michigan \yr1975 \publaddr Michigan \endref
\ref \key Mum1 \bysame \jour Intl. Symp. on Algebraic Geomtery
     \paper An Algebro-Geomtric Construction of Commuting Operators and
    of Solutions to the Toda Latticd Equation,
 Korteweg-de Vries Equation
  and Related Non-linear Equation \yr1977\publaddr Kyoto \pages 115-153
      \endref


\ref \key Mum2 \bysame \book Tata Lectures on Theta, vol II
     \publ  Birkh\"auser \yr1983-84 \publaddr Boston \endref
\ref \key Mum3 \bysame \paper Elastica and Computer Vision
     \book {\lq\lq}Elastica and Computer Vision"  in
          Algebraic Geometry and its Applications
      \eds C.~Bajaj, \publ Springer-Verlag \publaddr Berlin \yr 1993
      \pages 507-518 \endref
\ref \key Po \by H. Poincar\'e \book Papers on Fuchsian Functions \transl
     J. Stillwel \publ Springer \yr1985  \endref
\ref \key Ped \by F.~Pedit
\paper KdV flows on the Riemann sphere,
{\rm a talk at the
meeting on {\lq\lq}Study on Integrability in Differential Geometry",
 Lecture on Tokyo Metropritan University, Jan. 8-10, 1998}  \endref
\ref \key R \by P.~Ramond  \book Field Theory, A Modern Primer
     \yr1981 \publ Benjamin/Cummings \publaddr  Massachusetts \endref

\ref \key S \by M. Sato  \paper $D$-Modules and Nonlinear System
   \jour Adv. Stud. in Pure Math. \vol 19 \yr 1989
\pages 417-434 \endref
\ref \key SN \by M. Sato and M. Noumi
  \book Soliton equation and universal
Grassmannian manifold (in Japanese) \publ Shophia univ. \publaddr Tokyo
\yr 1984 \endref
\ref \key Se \by G.~Segal  \paper The Geometry of the KdV equation
  \book Topological Methods in Quantum Field Theory  \eds W. Nahm et.al.
     \yr 1990 \pages 96-106 \publ World Scienctific \publaddr Singapore
     \endref
\ref \key SS \by M. Sato and Y. Sato \paper Soliton  equations as
dynamical systems on infinite dimensional Grassmann manifold \inbook
Nonlinear Partial Differentail Equations in Applied Science
\ed  H. Fujita,
P. D. Lax and G. Strang \publ Kinokuniya/North-Holland\yr 1983
\publaddr Tokyo \endref
\ref \key SW \by G.~Segal and G.~Wilson
\paper Loop groups and equations
    of KdV type \jour IHES \vol 61 \yr 1985 \pages 5-65 \endref
\ref \key T1 \by C.~Truesdell \jour Bull. Amer. Math. Soc. \vol 9 \yr
    1983 \page 293-310 \paper The influence of elasticity on analysis:
    the classic heritage \endref
\ref \key T2 \bysame
\book Leonhrdi Euleri Opera Omnia
ser. Secunda XI; The Rational Mechanics
of flexible or elastic bodies 1638-1788
\publ Birkhauser Verlag \yr 1960
    \publaddr Berlin \endref
\ref \key Wi \by E.~Witten \book Elliptic Curves and Modular Forms
    in Algebraic Topology, Proceedings Princeton 1986 \paper
    The index of the Dirac Operator in Loop Space \eds P.~S.~Landweber
    \publ Springer \publaddr Berlin \yr 1986 \endref
\ref \key We \by A.~Weil
\book Number Theory: an approach through history;
From Haammurapi to Legendre
\publ Birkh\"auser \publaddr Cambridge \yr 1983
\endref
\ref \key Wh \by H.Whitney
  \jour Trans. Amer. Math. Soc. \vol 36 \yr
    1934 \page 63-89
\paper Analytic extensions of differentiable functions
            defined in closed sets\endref
\endRefs
